\documentstyle[amssymb,12pt]{amsart}%\input amssym.def
\newtheorem{prop}{Proposition}[section]
\newtheorem{prop:def}{Proposition-Definition}[section]

\newtheorem{lemma}{Lemma}[section]
\newtheorem{thm}{Theorem}[section]
\newtheorem{cor}{Corollary}[section]

\theoremstyle{remark}
\newtheorem{remark}{Remark}

\begin{document}

\newcommand{\nc}{\newcommand} \nc{\on}{\operatorname}

\nc{\pa}{\partial}

\nc{\cA}{{\cal A}} \nc{\cB}{{\cal B}}\nc{\cC}{{\cal C}} 
\nc{\cE}{{\cal E}} \nc{\cF}{{\cal F}}\nc{\cG}{{\cal G}}
\nc{\cH}{{\cal H}} \nc{\cI}{{\cal I}} \nc{\cJ}{{\cal J}}
\nc{\cK}{{\cal K}} \nc{\cL}{{\cal L}} \nc{\cM}{{\cal M}}
\nc{\cO}{{\cal O}} \nc{\cR}{{\cal R}} \nc{\cS}{{\cal S}}   
\nc{\cT}{{\cal T}} \nc{\cV}{{\cal V}}\nc{\cW}{{\cal W}}  
\nc{\cX}{{\cal X}}

\nc{\sh}{\on{sh}}\nc{\Id}{\on{Id}}\nc{\Diff}{\on{Diff}}
\nc{\Perm}{\on{Perm}}\nc{\conc}{\on{conc}}
\nc{\ad}{\on{ad}}\nc{\Der}{\on{Der}}\nc{\End}{\on{End}}
\nc{\no}{\on{no\ }} \nc{\res}{\on{res}}\nc{\ddiv}{\on{div}}
\nc{\Sh}{\on{Sh}} \nc{\card}{\on{card}}\nc{\dimm}{\on{dim}}
\nc{\Sym}{\on{Sym}} \nc{\Jac}{\on{Jac}}\nc{\Ker}{\on{Ker}}
\nc{\Vect}{\on{Vect}} \nc{\Spec}{\on{Spec}}\nc{\Cl}{\on{Cl}}
\nc{\Imm}{\on{Im}}\nc{\limm}{\lim}\nc{\Ad}{\on{Ad}}
\nc{\ev}{\on{ev}} \nc{\Hol}{\on{Hol}}\nc{\Det}{\on{Det}}
\nc{\Bun}{\on{Bun}}\nc{\diag}{\on{diag}}\nc{\pr}{\on{pr}} 
\nc{\Span}{\on{Span}}\nc{\Comp}{\on{Comp}}\nc{\Part}{\on{Part}}
\nc{\tensor}{\on{tensor}}\nc{\ind}{\on{ind}}\nc{\id}{\on{id}}
\nc{\ins}{\on{ins}}

\nc{\al}{\alpha}\nc{\g}{\gamma}\nc{\de}{\delta}
\nc{\eps}{\epsilon}\nc{\la}{{\lambda}}
\nc{\si}{\sigma}\nc{\z}{\zeta}

\nc{\La}{\Lambda}

\nc{\ve}{\varepsilon} \nc{\vp}{\varphi} 

\nc{\AAA}{{\mathbb A}}\nc{\CC}{{\mathbb C}}\nc{\ZZ}{{\mathbb Z}} 
\nc{\QQ}{{\mathbb Q}} \nc{\NN}{{\mathbb N}}\nc{\VV}{{\mathbb V}} 
\nc{\KK}{{\mathbb K}} 

\nc{\ff}{{\mathbf f}}\nc{\bg}{{\mathbf g}}
\nc{\ii}{{\mathbf i}}\nc{\kk}{{\mathbf k}}
\nc{\bl}{{\mathbf l}}\nc{\zz}{{\mathbf z}} 
\nc{\pp}{{\mathbf p}}\nc{\qq}{{\mathbf q}} 
\nc{\Assoc}{{\mathbf Assoc}}

% RACINET \def\SSh{{\mathop{\scriptstyle\bigsqcup\!\hspace{-2.5pt}\bigsqcup}}}
% \def\varSSh{{\mathop{\scriptstyle\!\vspace{+4pt}\coprod\!\hspace{-4.4pt}\coprod}}} % big sha

\def\SSh{{\mathop{\scriptstyle\amalg\!\hspace{-1.5pt}\amalg}}}
\def\Ts{{\mathop{\scriptstyle\amalg}}}

% VARIANTS \def\SSh{{\mathop{\amalg\!\hspace{-2.5pt}\amalg}}}
% \def\Ts{{\mathop{\amalg}}}

\nc{\ub}{{\underline{b}}}
\nc{\uk}{{\underline{k}}} \nc{\ul}{{\underline}}
\nc{\un}{{\underline{n}}} \nc{\um}{{\underline{m}}}
\nc{\up}{{\underline{p}}}\nc{\uq}{{\underline{q}}}
\nc{\us}{{\underline{s}}}\nc{\ut}{{\underline{t}}}
\nc{\uw}{{\underline{w}}}
\nc{\uz}{{\underline{z}}}
\nc{\ual}{{\underline{\alpha}}}\nc{\ualpha}{{\underline{\alpha}}}
\nc{\ugamma}{{\underline{\gamma}}}
\nc{\ula}{{\underline{\lambda}}}\nc{\umu}{{\underline{\mu}}}
\nc{\unu}{{\underline{\nu}}}\nc{\usigma}{{\underline{\sigma}}}
\nc{\utau}{{\underline{\tau}}}
\nc{\uN}{{\underline{N}}}\nc{\uM}{{\underline{M}}}
\nc{\uK}{{\underline{K}}}

\nc{\A}{{\mathfrak a}} \nc{\B}{{\mathfrak b}} \nc{\G}{{\mathfrak g}}
\nc{\D}{{\mathfrak d}} \nc{\HH}{{\mathfrak h}}  \nc{\iii}{{\mathfrak
i}}  \nc{\K}{{\mathfrak k}}  \nc{\mm}{{\mathfrak m}} 
\nc{\N}{{\mathfrak
n}}\nc{\ttt}{{\mathfrak{t}}}  \nc{\U}{{\mathfrak u}}\nc{\V}{{\mathfrak
v}}

\nc{\SL}{{\mathfrak{sl}}}

\nc{\SG}{{\mathfrak S}}

\nc{\wt}{\widetilde} \nc{\wh}{\widehat}
\nc{\bn}{\begin{equation}}\nc{\en}{\end{equation}} \nc{\td}{\tilde}

% ****** GISPIC **********
%
%** by GISLI MASON *******
%
%**for commutative diagrams
%

\newcommand{\ldar}[1]{\begin{picture}(10,50)(-5,-25)
\put(0,25){\vector(0,-1){50}}
\put(5,0){\mbox{$#1$}} 
\end{picture}}

\newcommand{\lrar}[1]{\begin{picture}(50,10)(-25,-5)
\put(-25,0){\vector(1,0){50}}
\put(0,5){\makebox(0,0)[b]{\mbox{$#1$}}}
\end{picture}}

\newcommand{\luar}[1]{\begin{picture}(10,50)(-5,-25)
\put(0,-25){\vector(0,1){50}}
\put(5,0){\mbox{$#1$}}
\end{picture}}

\title[Quantization of Lie bialgebras and shuffle algebras] 
{Quantization of Lie bialgebras and \\ shuffle algebras of Lie algebras}

\author{B. Enriquez} 

\address{D\'epartement de Math\'ematiques et Applications, 
Ecole Normale Sup\'erieure, UMR du CNRS, 75005 Paris, France}

\date{September 2000}

\begin{abstract}
To any field $\KK$ of characteristic zero, we associate a set 
$\SSh(\KK)$. Elements of $\SSh(\KK)$ are equivalence classes of 
families of Lie polynomials subject to associativity relations. 
We construct an injection and a retraction between $\SSh(\KK)$ and the set 
of quantization functors of Lie bialgebras over $\KK$.  
This construction involves the following steps. 1) To each
element $\varpi$ of $\SSh(\KK)$, we associate a functor
$\A\mapsto\Sh^\varpi(\A)$ from the category of Lie bialgebras to that of Hopf
algebras; $\Sh^\varpi(\A)$ contains $U\A$. 
2) When $\A$ and $\B$ are Lie algebras, and $r_{\A\B}
\in\A\otimes\B$, we construct an element $\cR^\varpi(r_{\A\B})$
of $\Sh^\varpi(\A)\otimes\Sh^\varpi(\B)$
satisfying quasitriangularity identities; in particular, 
$\cR^\varpi(r_{\A\B})$ defines a Hopf algebra morphism from 
$\Sh^\varpi(\A)^*$ to $\Sh^\varpi(\B)$. 3) When $\A = \B$ and 
$r_\A\in\A\otimes\A$  is a solution of CYBE, we construct
a series $\rho^\varpi(r_\A)$ such that $\cR^\varpi(\rho^\varpi(r_\A))$
is a solution of QYBE. The expression of $\rho^\varpi(r_\A)$ in terms of 
$r_\A$ involves Lie polynomials, and we show that this expression 
is unique at a universal level. This step relies on vanishing 
statements for cohomologies arising from universal algebras for 
the solutions of CYBE. 4) We define the quantization of a Lie bialgebra $\G$ as 
the image of the morphism defined by $\cR^\varpi(\rho^\varpi(r))$, where
$r\in\G\otimes\G^*$ is the canonical element attached to $\G$. 
\end{abstract}

\maketitle

\subsection*{Introduction}

According to Drinfeld, a quantum group is a formal deformation
of the universal enveloping algebra of a Lie algebra $\G$. 
The semiclassical structure associated with such a deformation
is a Lie bialgebra structure on $\G$. The quantization 
problem of Lie bialgebras, as posed by Drinfeld in \cite{QG,unsolved}, 
is to construct a functor from the category of Lie bialgebras 
to that of quantum groups, whose composition with the 
``semiclassical limit'' functor is the identity. Such an 
object is called a quantization functor. When the structure 
constants of the quantum group can be expressed by polynomial 
formulas in terms of those of the Lie bialgebra, the quantization 
functor is called universal.

An approach to the construction of universal quantization 
functors was proposed by 
N.\ Reshetikhin (\cite{Resh}). Later, it was 
solved by P.\ Etingof and D.\ Kazhdan (\cite{EK}). 
Their quantization procedure involves associators. 
An associator in an element of an abstract algebra, subject 
to certain conditions. An example of an associator is Drinfeld's
``Knizhnik-Zamolodchikov associator'', which involves special values
of multiple zeta functions. Drinfeld also proved the existence of 
nontrivial associators defined over $\QQ$. 
Let $\KK$ be a field of 
characteristic zero, and let us denote by $\Assoc(\KK)$ the set of 
associators defined over $\KK$. The main result of \cite{EK} is the
construction of a map  
$EK : \Assoc(\KK) \to \{$universal quantization functors of Lie bialgebras 
over $\KK\}$. 

Our purpose in this paper is to study further the set of quantization
functors of Lie bialgebras. Our main result may be stated
as follows. We introduce a set $\SSh(\KK)$ of equivalence 
classes of Lie polynomials, satisfying certain associativity 
equations, and we define two maps $\al_\KK : \SSh(\KK) \to 
\{$universal quantization functors of Lie bialgebras 
over $\KK\}$ and $\beta_\KK : \{$universal 
quantization functors of Lie bialgebras 
over $\KK\}\to \SSh(\KK)$, such that $\beta_\KK\circ\al_\KK = 
\id_{\SSh(\KK)}$. In particular, the composition $\beta_\KK\circ EK$
yields a natural map $\Assoc(\KK)\to\SSh(\KK)$. 
The image of $\al_\KK$ is contained in the set of 
universal quantization functors
with functorial behavior with respect to the operations of 
dualizing and taking the double of Lie bialgebras. We set up a 
bijection between this set of quantization functors with the 
product of $\SSh(\KK)$ with a universal group $\cG_0$. 

Our approach may be viewed as close to the original 
approach of Reshetikhin. To describe the latter, it is useful to 
recall how quantum Kac-Moody algebras were initially constructed in 
\cite{Drinfeld,Jimbo}. 

Let  $\A$ be a Kac-Moody Lie algebra, and let $\delta : \A\to
\wedge^2\A$ be the cocycle defining its standard (Drinfeld-Sklyanin)
Lie bialgebra structure.  Then there exist two opposite Borel subalgebras
$\B_+$ and $\B_-$ of  $\A$, which are Lie subbialgebras of $\A$. The
quantization of $\B_\pm$  is constructed as follows. Let $h_i,x_i^\pm$
be the generators of $\B_\pm$.  Then $\delta$ is expressed simply on
these generators, by $\delta(h_i) = 0$,  $\delta(x_i^\pm) = \pm d_i
h_i\wedge x_i^\pm$. One then defines the coproduct  $\Delta$ on these
generators by $\Delta(h_i) = h_i\otimes 1 + 1\otimes h_i$,
$\Delta(x_i^\pm) = x_i^\pm \otimes e^{\pm d_i h_i}+ 1\otimes x_i^\pm$;
up to  now the only condition on $\Delta$ is that it should be
coassociative and  it should deform $\delta$.  This means that we have
defined Hopf algebra structures on the tensor algebras of $(\oplus_i \CC
h_i) \oplus  (\oplus_i \CC x_i^\pm)$. One then finds relations between
the generators $h_i$ and $x_i^\pm$,  which deform the classical relations, and are
skew-primitive with respect to $\Delta$. One may then show that the
resulting algebras are flat deformations of $U\B_\pm$,  and that the
relations are the generators of the radical of a Hopf pairing between 
$T\big((\oplus_i \CC h_i) \oplus  (\oplus_i \CC x_i^+)\big)$  and 
$T\big((\oplus_i
\CC h_i) \oplus  (\oplus_i \CC x_i^-)\big)$  (see e.g.\ \cite{Lusztig}).  
$U_\hbar\A$ is then constructed as a subalgebra of the quantum double of 
$U_\hbar\B_+$, because this quantum double contains two copies of the Cartan 
subalgebra of $\A$. The procedure is the same in the case of the 
quantum current algebras (``new realizations algebras'').

The approach of Reshetikhin to the quantization problem was to  imitate these steps in 
the general case. Let $(\G,\delta_\G)$ be a Lie bialgebra and 
let $(\G^*,\delta_{\G^*})$ be the dual Lie bialgebra. 
Then the tensor algebras $T(\G)$ and $T(\G^*)$ of $\G$
and $\G^*$ are cocommutative Hopf algebras, and $\delta_\G$
and $\delta_{\G^*}$ extend to Hopf co-Poisson structures on 
 $T(\G)$ and $T(\G^*)$. The first step is then to construct
coproducts $\Delta_\G$ and $\Delta_{\G^*}$ on these tensor 
algebras, deforming $\delta_\G$ and $\delta_{\G^*}$. 
The second step is to construct a Hopf algebra 
pairing between $(T(\G),\Delta_{\G})$ and $(T(\G^*),\Delta_{\G^*})$
and define the quantizations of $\G$ and $\G^*$ as the 
quotients of $T(\G)$ and $T(\G^*)$ by the radicals of this
pairing.  The pairing between $T(\G)$
and $T(\G^*)$ should be chosen in such a 
way that the quotients are 
flat deformations of $U\G$ and $U(\G^*)$.  
  
In this paper, we study this program placing ourselves 
in the dual framework. 
For any vector space $V$, the dual to the Hopf algebra 
structure on its tensor algebra $T(V)$ is a Hopf algebra
structure defined on the shuffle algebra $\Sh(V^*)$ of its 
dual. The duals
of $(T(\G),\delta_\G)$ and $(T(\G^*),\delta_{\G^*})$ are  
then commutative Hopf-Poisson structures on $\Sh(\G^*)$ and $\Sh(\G)$.
Moreover, the map $\G\mapsto \Sh(\G)$ is a functor from the 
category of Lie algebras to that of Hopf-Poisson algebras. 
The dual translation of the first step of the above program is 
to construct a functorial quantization of these 
Hopf-Poisson structures.

In Section \ref{sect:sh:g}, we introduce the main new object of this
paper. This is a functor $\KK\mapsto \SSh(\KK)$ from the category of 
fields to that of sets. Elements of $\SSh(\KK)$ are equivalence 
classes of families of Lie polynomials, subject to certain 
relations (see Section \ref{sect:ssh}). We construct a bijection 
between $\SSh(\KK)$ and the set of universal quantization functors of the 
Hopf-Poisson structure of $\Sh(\G)$ (Proposition 
\ref{prop:quant:functors}). 
In Section \ref{sect:cort}, we explain the connection 
between the quantizations of $\Sh(\G)$ provided by elements of 
$\SSh(\KK)$ and the PBW quantization of $S\G$ 
(\cite{Cortinas,Kath}). 
(The shuffle algebras appeared many times in the theory of 
Hopf algebras and quantum groups (see 
\cite{Nichols,Andr,Varchenko,Schauenburg,Rosso}), but the
construction of $\Sh(\G)$ using elements of $\SSh(\KK)$ 
seems to be new.) 
If $\varpi$ is an element of $\SSh(\KK)$, we denote by $\G\mapsto 
\Sh^\varpi(\G)$ the corresponding quantization functor. 
 
When $(A,m_A,\Delta_A)$ and $(B,m_B,\Delta_B)$ are any finite-dimensional Hopf
algebras, a Hopf pairing between them is the same as an element 
$\cR_{AB}$ in $A^*\otimes B^*$, such that  $(\Delta_{A^*}\otimes
\id)(\cR_{AB}) = \cR_{AB}^{(13)}\cR_{AB}^{(23)}$  and $(\id \otimes
\Delta_{B^*})(\cR_{AB}) = \cR_{AB}^{(12)}\cR_{AB}^{(13)}$  (where
$\Delta_{A^*}$  and $\Delta_{B^*}$ are the dual maps to the 
multiplication maps of $A$ and $B$). In that case, the map 
$\ell : B^*\to A$ defined by $\ell(\xi) = \langle \cR_{AB},\id\otimes
\xi\rangle$ defines a Hopf algebra morphism from $B^*$ to $A$. 
With some restrictions, the same
is true in our situation. Our next step is  then to associate to any
$\varpi\in\SSh(\KK)$, any pair $(\A,\B)$ of Lie 
algebras and any element $r_{\A\B}$ of
$\A\otimes \B$, a family of elements $\cR^{\varpi}_n(r_{\A\B})$ in 
$\Sh^\varpi(\A)\otimes\Sh^{\varpi}(\B)$, such that 
if $t$ is any formal parameter, $\cR^{\varpi}(r_{\A\B})
= \sum_{n\geq 0} t^n \cR^{\varpi}_n(r_{\A\B})$ satisfies these identities 
(Proposition \ref{prop:QT}). 

In the rest of the construction, we fix an element $\varpi$ in 
$\SSh(\KK)$. 
Assume that  $\A = \B$ and set $r_\A = r_{\A\B}$. 
It is natural to ask when 
$\cR^\varpi(r_{\A})$ satisfies the quantum Yang-Baxter equation
(QYBE). To state our next main result, we need some notation. 
Denote by $FL_n$ the part of the free Lie algebra in $n$ generators, 
homogeneous of degree one in each generator. The symmetric
group $\SG_n$ acts by simultaneous permutations of generators of both 
factors of $FL_n\otimes FL_n$. We set $F_n = 
(FL_n \otimes FL_n)_{\SG_n}$. Then there is a unique linear 
map $\kappa^{(ab)}_n(r_\A) : F_n\to \A\otimes\A$, sending the 
class of $P(x_1,\ldots,x_n)\otimes Q(x_1,\ldots,x_n)$ to 
$\sum_{i_1,\ldots,i_n\in I} P(a_{i_1},\ldots,a_{i_n})
\otimes Q(b_{i_1},\ldots,b_{i_n})$, if $r_\A$ has the form 
$\sum_{i\in I} a_i\otimes b_i$. 
Our next main result is (see Theorem \ref{thm:unique:sol}, 
Lemma \ref{lemma:3.3}, Corollary \ref{cor:malta} and 
Proposition \ref{prop:equivalence}) 

\begin{thm} \label{big:thm:QYBE}
There exists a family $(\varrho^\varpi_n)_{n\geq 1}$ in $\prod_{n\geq 1}
F_n$, such that $\varrho^\varpi_1 = x_1\otimes x_1$, and if $\A$ is  Lie 
algebra and $r_\A\in\A\otimes\A$ is a 
solution of the classical Yang-Baxter equation (CYBE), then  
$\cR^\varpi(\sum_{n\geq 1} \kappa_n^{(ab)}(r_\A)(\varrho^\varpi_n))$ is a 
solution of QYBE in $\Sh^\varpi(\A)^{\otimes 3}$. 
\end{thm}

We find the family $(\varrho^\varpi_n)_{n\geq 1}$ as the solution of a
system of equations in $\prod_{n\geq 1} F_n$, the 
{\it universal Lie QYB equations} (see Section 
\ref{sect:univ:Lie}). More precisely, we prove 
that the family $(\varrho^\varpi_n)_{n\geq 1}$
is the unique solution to these equations (see Theorem 
\ref{thm:unique:sol}). 

For any $\varpi$, we have 
$\varrho^\varpi_2 = {1\over 8}[x_1,x_2] \otimes [x_1,x_2]$, so that 
the expansion of $\rho^\varpi(r_\A) = \sum_{n\geq 1}
\kappa_n^{(ab)}(r_\A)(\varrho^\varpi_n)$ is 
$$
\rho^\varpi(r_\A) = \sum_i a_i\otimes b_i + {1\over 8} \sum_{i,j\in I} 
[a_i,a_j] \otimes [b_i,b_j] + \cdots. 
$$

The universal Lie QYB equations are equalities 
in an algebra $F^{(3)}$, which is a particular instance of a 
family of algebras $F^{(N)}$ 
(see Appendices \ref{proof:prop:key} and \ref{proof:prop:delta:4}). 
The algebras $F^{(N)}$ are universal algebras for  the pairs $(\A,r_\A)$ of
a Lie algebra $\A$ and a solution $r_\A\in\A\otimes\A$ of CYBE. 
When $A$ is any algebra, we connect the problem of
deforming a solution $r_A\in A\otimes A$ of CYBE into a solution of 
QYBE with the Lie coalgebra structure on $A$ defined by $r_A$, more
precisely, with the corresponding Lie coalgebra cohomology 
(Proposition \ref{prop:leeds:first}; the complete statement 
is in Theorem \ref{def:theory} of Appendix \ref{proof:problem:3}). 
Adapting this result to a family of universal shuffle algebras $\Sh^{(F)}_k$
(see Theorem \ref{thm:univ} and Appendix \ref{sect:proof:thm:univ}), we  
formulate the universal Lie QYB equations in terms of cohomology 
groups $H^2_n$ and $H^3_n$, constructed in terms of the family 
$(F^{(N)})_{N = 2,3,4}$. 
We then show vanishing statements for the groups $H^2_n$ and $H^3_n$
(Theorem \ref{results:cohom} and Appendix \ref{app:cohomologies}), 
using identities in free Lie algebras (Propositions \ref{reut}
and \ref{prereut}, see also \cite{Reut}) . This shows the existence and 
unicity of the solution $(\varrho^\varpi_n)_{n\geq 1}$ to the universal 
Lie QYB equations (Theorem \ref{thm:unique:sol}). 

We then apply this construction to the problem of universal quantization 
of Lie bialgebras (Section \ref{sect:quantization}). 
To a finite-dimensional Lie bialgebra $\G$ over $\KK$ 
are attached its
dual Lie bialgebra $\G^*$, its double $\D$, and the canonical 
$r$-matrix $r_\G\in\G\otimes\G^*$ (see Section \ref{sect:basics} and
\cite{QG}). $\G$ and $\G^*$ are Lie subalgebras of $\D$, so 
$r_\G\in\D\otimes\D$; $r_\G$ is then a solution of CYBE. 
Then 
$$
\cR^\varpi(\rho^\varpi(\hbar r_\G))\in \Sh^\varpi(\G)\otimes\Sh^\varpi(\G^*)[[\hbar]]
\subset \Sh^\varpi(\D)^{\otimes 2}[[\hbar]],
$$ 
and as an element of the  latter algebra, $\cR^\varpi(\rho(\hbar r_\G))$ 
is a solution of QYBE  
(here $\hbar$ is a formal parameter). As we explained above, 
$\cR^\varpi(\rho^\varpi(\hbar r_\G))$ induces a Hopf algebra morphism $\ell$ 
from a Hopf algebra $T_{\hbar}^{\varpi}(\G)$ dual to 
$\Sh^\varpi(\G^*)$ to 
$\Sh^\varpi(\G)[[\hbar]]$. We then define $U^\varpi_\hbar\G$ as the image 
$\Imm(\ell)$ of $\ell$ and show that 
is a quantization of the Lie bialgebra $\G$ (Theorem 
\ref{quantization:QLBA}). We also prove that $\Imm(\ell)$
is divisible in $\Sh^\varpi(\G)[[\hbar]]$ 
(as a $\KK[[\hbar]]$-module). The fact 
that $\cR^\varpi( \rho^\varpi(\hbar r_\G))$ is a solution of QYBE plays an essential role in the
proof of both results; the idea from \cite{Enr} that the subalgebras
of shuffle algebras are necessarily torsion-free is also at the 
basis of the construction. 
We show that the map $(\G,[\ ,\ ]_\G,\delta_\G) \mapsto U^\varpi_\hbar\G$ 
defines a functor from the category of Lie bialgebras to that 
of Hopf algebras (Proposition \ref{prop:functor}). We comment on 
the extension of this functor to infinite-dimensional Lie 
bialgebras in Remark \ref{rem:inf:dim}. 

We also characterize the quantized formal series Hopf
(QFSH) algebra $\cO^\varpi_\hbar(G^*)$ corresponding to 
$U^\varpi_\hbar\G$ as the intersection of 
$\Imm(\ell)$ with a subalgebra $\Sh^\varpi_\hbar(\G)$ of 
$\Sh^\varpi(\G)[[\hbar]]$
(Proposition \ref{prop:QFSH}).
According to a general result
(\cite{QG,Gav}), $\cO^\varpi_\hbar(G^*)$ is a quantization   
of the Hopf algebra of functions on the formal Poisson-Lie
group corresponding to $\G^*$; we show that $\cO^\varpi_\hbar(G^*)$ coincides
with the dual $(U^\varpi_\hbar\G^*)^*$ of $U^\varpi_\hbar\G^*$.

According to Theorem \ref{quantization:QLBA} and 
Proposition \ref{prop:functor}, there is a unique map 
$\al_\KK : \SSh(\KK) \to \{$universal quantization functors
of Lie bialgebras over $\KK\}$ 
such that $\al_\KK(\varpi) = (\G\mapsto 
U_\hbar^\varpi(\G))$. On the other hand, we define 
a map $\beta_\KK : \{$universal quantization functors
of Lie bialgebras over $\KK\} \to \SSh(\KK)$ as follows. 
When 
$V$ is a vector space, denote by $F(V)$ the free Lie algebra of 
$V$. If $\G$
is a Lie algebra, then $F(\G^*)$ is naturally equipped with a 
structure of Lie bialgebra. One shows that the restriction 
of a quantization functor $Q$ to bialgebras of the form $F(\G^*)$
yields an element $\beta_\KK(Q)$ of $\SSh(\KK)$.  

If $\G\mapsto Q(\G)$ is a universal quantization functor 
for Lie bialgebras, let us say that it is compatible with the
operations of dual and double if the following holds: if
$\G$ is any Lie bialgebra, then there is are canonical isomorphisms
between $Q(\G^*)$ and $(Q(\G)^*)^\vee$ ($O^\vee$ is the quantized 
universal enveloping algebra
attached to any quantized formal series Hopf algebra $O$, see 
\cite{QG,Gav}) and between $Q(D(\G))$ and the quantum double of $Q(\G)$
($D(\G)$ is the double of a Lie algebra $\G$). In that case, we will 
say that $Q$ is a compatible functor.  
 
Let us denote by $\cG_0$ the group of all functorial assignments $\A\mapsto 
\rho_\A$, where $\A$ is a Lie bialgebra and 
$\rho_\A\in \End(\A)[[\hbar]]$, such that 
$\rho_\A$ is an iterate of tensor products of the Lie bracket and 
cobracket, and satisfies the identities $\rho_\A([x,y]) = [\rho_\A(x),y]$, 
$\rho_{\A^*} = \rho_\A^t$ and $\rho_\A = \id_\A + o(\hbar)$. 
Then $\cG_0$ may be viewed as a group of universal transformations
of the solutions of CYBE arising from Lie bialgebras (see Lemma 
\ref{lemma:G0} and Remark \ref{rem:G0}), and as a ``linear'' subgroup 
of the group of functorial assignments $(\A\mapsto\rho_\A)$, 
such that if $[\ ,\ ]_\A$ and $\delta_\A$ are the Lie bracket and
cobracket of $\A$, then $(\A,[\ ,\ ]_\A,(\rho_\A\otimes\rho_\A)\circ
\delta_\A\circ\rho_\A^{-1})$ is a Lie bialgebra. 

In Section \ref{sect:last}, 
we prove 
\begin{thm} \label{thm:last}
The map $\al_\KK$ is an injection and $\beta_\KK$ is a retraction of 
$\al_\KK$, that is $\beta_\KK\circ\al_\KK = \id_{\SSh(\KK)}$. 
These maps depend functorially on $\KK$ and are
 equivariant with respect to 
the natural actions of $\on{Aut}(\KK[[\hbar]])$
on $\SSh(\KK)$ and 
$\{$universal quantization functors
of Lie bialgebras over $\KK\}$. 

We have $\al_\KK(\SSh(\KK))\subset \{$compatible  
quantization functors for Lie bialgebras$\}$. There is a  
bijection $\ve_\KK$ between $\SSh(\KK)\times\cG_0$
and  $\{$compatible 
quantization functors for Lie bialgebras$\}$, such that 
$\beta_\KK\circ\ve_\KK$ is the projection on the first factor and the
restriction of $\ve_\KK$ on $\SSh(\KK)\times \{$neutral element$\}$ 
coincides with $\al_\KK$. 
\end{thm}

In particular, if $\cG_0$ reduces to multiplication by scalars, 
then $\al_\KK$ sets up a bijection between
$\SSh(\KK)$ and $\{$compatible quantization functors for Lie 
bialgebras$\}$.  

The proof of this Theorem is in Section \ref{sect:last}. It 
relies on the following idea. If $\G$ is a Lie algebra, then 
the universal enveloping algebra of $F(\G)$ is 
the Hopf-co-Poisson algebra $(T(\G),\delta_\G)$. The first part 
of the proof of 
Theorem \ref{thm:last} involves the fact that the quantizations of 
$(T(\G),\delta_\G)$ given by $\Sh^\varpi(\G)^*$ and 
$U_{\hbar}^{\varpi}(F(\G^*))$
are naturally isomorphic (Section \ref{natural:isom}). The second 
part of the proof relies on the fact that any Lie bialgebra $\G$
may be viewed as the image of a Lie bialgebra morphism from 
$F(\G)$ to $F(\G^*)^*$. 

\medskip 

The map $\Assoc(\KK)\to \SSh(\KK)$ may help 
to study how the quantization functors of Etingof and Kazhdan 
depend on associators. More precisely, the image of the 
map $EK$ is contained in $\{$compatible quantization functors
for Lie bialgebras$\}$, so if $\cG_0$ is trivial, then the map 
$EK(\Assoc(\KK))\to \SSh(\KK)$ is an injection. 

For some classes of Lie algebras, the 
Etingof-Kazhdan quantizations
are known to be independent of the associator (for example in the case
of Kac-Moody algebras (see \cite{EK:KM}) and 
nondegenerate triangular Lie bialgebras (\cite{Moreno})). 

\medskip 
\subsection*{Acknowledgements}
 
I would like to express my hearty thanks to N.\ Andruskie-witsch;  
this project was started in collaboration with him during 
my visit at C\'ordoba (Argentina) in July 1998.  
I would also like to thank P.\ Etingof for discussions about
the relation of the present work with his and D.\ Kazhdan's
work, and M.\ Jimbo, who explained 
to me how he obtained the quantum Kac-Moody algebras in \cite{Jimbo}. 
I thank Y.\ Kosmann-Schwarzbach, J.-H.\ Lu and D.\ Manchon
for discussions about this work, and 
P.\ Bressler and B.\ Tsygan who pointed out an important
mistake in an earlier version of it.  
Finally, I am grateful to J.\ Stasheff for several ``red-ink pen'' 
comments on the first drafts of the paper.

\medskip\noindent 

\newpage 

\section{$\SSh(\KK)$ and shuffle algebras of Lie algebras}
\label{sect:sh:g}

\subsection{Definition of $\SSh(\KK)$} \label{sect:ssh}

Let $R$ be a commutative ring. Let us define $FL_n(R)$ as the 
part of the free Lie algebra over $R$ with $n$ generators,  
homogeneous of degree $1$ in each generator. If $\la$ belongs to $R$, 
let us define 
$\cB(R)_\la$ as the set of all families $(B_{pq})_{p,q\geq 0}$, where 
each $B_{pq}$ belongs to $FL_{p+q}(R)$, such that $(B_{pq})_{p,q\geq 0}$ 
satisfies the equations
$$
B_{10}(x) = B_{01}(x) = x, \quad B_{p0} = B_{0p} = 0
\ \on{if}\ p\neq 1, 
\quad B_{11}(x_1,x_2) = \la [x_1,x_2], 
$$ 
and 
\begin{align} \label{cond:Bpq}
& \sum_{\al>0} 
\sum_{(p_\beta)_{\beta = 1,\ldots,\al}\in\Part_{\al}(p),
(q_\beta)_{\beta = 1,\ldots,\al}\in\Part_{\al}(q)} 
B_{\al r} \big( B_{p_1q_1}(x_1,\ldots,x_{p_1}|y_1,\ldots,y_{q_1}) 
\cdots \\ & \nonumber 
\cdots B_{p_\al q_\al}(
x_{\sum_{\beta = 1}^{\al - 1}p_\beta + 1},
\ldots, x_p|y_{\sum_{\beta = 1}^{\al - 1}q_\beta + 1},
\ldots, y_q)
|z_1,\ldots,z_r \big)
\\ & \nonumber 
= \sum_{\al>0} 
\sum_{(q_\beta)_{\beta = 1,\ldots,\al}\in\Part_{\al}(q),
(r_\beta)_{\beta = 1,\ldots,\al}\in\Part_{\al}(r)} 
B_{p \al} \big( x_1,\ldots, x_p| 
\\ & \nonumber 
B_{q_1r_1}(y_1,\ldots,y_{q_1}|z_1,\ldots,z_{r_1}) 
\cdots 
B_{q_\al r_\al}(
y_{\sum_{\beta = 1}^{\al - 1}q_\beta + 1},
\ldots, y_q|z_{\sum_{\beta = 1}^{\al - 1}r_\beta + 1},
\ldots, z_r)
\big)
\end{align}
for any integers $p,q,r>0$. Here $\Part_\al(n)$ is the set of 
$\al$-partitions of $n$, that is the set of $\al$uples
$(n_1,\ldots,n_\al)$ of nonnegative integers such that 
$\sum_{\beta = 1}^{\al} n_\beta = n$. We also set $\cB(R) = 
\coprod_{\la\in R}\cB(R)_\la$. 

Let us define $\cG(R)$ as the subset of the product 
$\prod_{n\geq 1} FL_n(R)$ of elements $(P_n)_{n\geq 1}$ such that 
$P_1(x) = x$. 
If $P = (P_n)_{n\geq 1}$ and $Q = (Q_n)_{n\geq 1}$ belong to $\cG(R)$, 
let us define $(P*Q)_n$ as the element of $FL_n(R)$ equal to
\begin{align} \label{group:G}
& (P*Q)(x_1,\ldots,x_n) = \sum_{\al>0} \sum_{(n_1,\ldots,n_\al)
\in\Part_\al(n)} P_\al(Q_{n_1}(x_1,\ldots,x_{n_1}),\ldots,
\\ & \nonumber Q_{n_\al}(x_{\sum_{\al' = 1}^\al n_{\al'} + 1},\ldots,x_n)) . 
\end{align}
Then $*$ defines a group structure on $\cG(R)$. 

If $P = (P_n)_{n\geq 1}$ belongs to $\cG(R)$ and 
$B = (B_{pq})_{p,q\geq 0}$ belongs to $\cB(R)$, let us define 
$(P*B)_{pq}$ as the element of $FL_{p+q}$ defined by 
\begin{align*}
& (P*B)_{pq}(x_1,\ldots,x_p|y_1,\ldots,y_q) = 
\sum_{\al,\beta>0} \sum_{(p_1,\ldots,p_\al)\in \Part_\al(p), 
(q_1,\ldots,q_\beta)\in \Part_\beta(q)}
\\ & 
B_{\al\beta}\big(P_{p_1}(x_1,\ldots,x_{p_1}),\ldots,P_{p_\al}
(x_{\sum_{\al' = 1}^{\al - 1} p_{\al'} + 1},\ldots,x_p)|
\\ & 
P_{q_1}(y_1,\ldots,y_{q_1}),\ldots,P_{q_\al}(y_{\sum_{\beta' = 1}^{\beta - 1}
q_{\al'} + 1},\ldots,y_q)\big) . 
\end{align*}
Then $*$ defines an action of $\cG(R)$ on each $\cB(R)_\la$ and on 
$\cR(R)$.   

We define  $\wt{\SSh}(R)$ and $\wt{\SSh}(R)_\la$ as the quotient sets $\cB(R) / \cG(R)$
and $\cB(R)_\la / \cG(R)$.  
If $\KK$ is a field of characteristic zero, we set 
$\SSh(\KK) = \wt\SSh(\KK[[\hbar]])$, where $\hbar$ is a formal 
variable. We also set $\SSh(\KK)_\la = \wt{\SSh}(\KK[[\hbar]])_\la$, 
if $\la\in\KK[[\hbar]]$. 

If $R$ is a ring, then there is a unique involution $\omega\mapsto 
\omega^\vee$ of $\cB(R)$, such that if $\omega = (B_{pq})_{p,q\geq 0}$, 
then $\omega^\vee = (B^\vee_{pq})_{p,q\geq 0}$, where 
$B^\vee_{pq}(x_1,\ldots,x_p|y_1,\ldots,y_q) = 
B_{pq}(x_p,\ldots,x_1|y_q,\ldots,y_1)$. The map $\omega\mapsto 
\omega^\vee$ induces an involution $\varpi\mapsto\varpi^\vee$ of
$\wt\SSh(R)$ and of $\SSh(\KK)$, if $\KK$ is any field. 

The group $R^\times$ also acts on $\cB(R)$ by the rule 
$r\cdot (B_{pq})_{p,q} = (r^{p+q - 1} B_{pq})_{p,q}$. This action
induces an action of $R^\times$ on $\wt{\SSh}(R)$.   
The group $\on{Aut}(R)$ also acts on $\wt{\SSh}(R)$ in the obvious way. 

\begin{remark} 
Define $\cB_n(R)$ as the set of families 
$(B_{pq})_{p,q|p,q\geq 0,p+q\leq n}$ 
satisfying the relations (\ref{cond:Bpq}), for any 
$p,q,r$ such that $p+q+r\leq n$. Define $\cG_n(R)$
as the set of families $(P_k)_{k|1\leq k\leq n}$. 
Then the rule (\ref{group:G}) equips  $\cG_n(R)$
with a group structure. Moreover, there are natural 
projection maps $\cB_{n+1}(R)\to \cB_n(R)$ and 
$\cG_{n+1}(R)\to\cG_n(R)$, which induce maps 
$\wt\SSh_{n+1}(R) \to \wt\SSh_n(R)$. Then $\wt\SSh(R)$
is the projective limit $\limm_{\leftarrow n}\wt\SSh_n(R)$. 
\end{remark}

\subsection{Hopf-Poisson structures on shuffle algebras}

Let $\KK$ be a field of characteristic zero and let  
$V$ be a vector space over $\KK$. The shuffle algebra
$\Sh(V)$ is a commutative Hopf algebra; it is defined as follows.   
As a vector space, $\Sh(V)$ is isomorphic with the tensor 
algebra $T(V)$. If $v_1,\ldots,v_n$ are elements of $V$, we denote by 
$(v_1\ldots v_n)$ the element of $\Sh(V)$ corresponding to 
$v_1\otimes\cdots\otimes v_n$ by this isomorphism.  

If $p$ and $q$ are integers, let us denote by $\SG_{p,q}$ 
the subset of $\SG_{p+q}$ of all permutations $\sigma$ such that 
if $(i,j)$ is a pair of integers such that $1\leq i <j\leq p$
or $p+1\leq i < j\leq p+q$, then $\sigma(i) < \sigma(j)$
(shuffle permutations). The space $\Sh(V)$ is equipped with a 
multiplication $m_0$ defined by 
$$
m_0((v_1\ldots v_n)\otimes (v_{n+1}\ldots v_{n+m})) = 
\sum_{\sigma\in\SG_{n,m}} (v_{\sigma(1)}\ldots v_{\sigma(n+m)}) 
$$
and a comutiplication $\Delta_{\Sh(V)}$ defined by 
$$
\Delta_{\Sh(V)}((v_1\ldots v_n)) = \sum_{k = 0}^n 
(v_1\ldots v_k)\otimes (v_{k+1}\ldots v_n).  
$$
If we define the unit of $\Sh(V)$ as $1\in T^0(V)$, its counit 
$\varepsilon_{\Sh(V)}$ as the projection map on $T^0(V) = \KK$ parallel to the sum 
$\oplus_{n>0} T^n(V)$, and the antipode by the formula 
$S_{\Sh(V)}((v_1\ldots v_p)) = (-1)^p (v_p\ldots v_1)$, then  
$(\Sh(V),m_0,\Delta_{\Sh(V)},\varepsilon_{\Sh(V)},S_{\Sh(V)})$
is a commutative Hopf algebra. 

Assume that $V$ is the underlying vector space of a Lie algebra $\G$. 
Then there is a unique linear map $m_1 : \Sh(\G)^{\otimes 2}\to\Sh(\G)$, 
such that 
\begin{align*}
& m_1((x_1\ldots x_n)\otimes (x_{n+1}\ldots x_{n+m})) 
\\ & 
= \sum_{i = 1}^n \sum_{j = 1}^m
\sum_{\sigma'\in\SG_{i-1,j-1}, \sigma''\in\SG_{n-i,m-j} }
( y^{(i,j)}_{\sigma'(1)} \ldots y^{(i,j)}_{\sigma'(i+j-2)}
[x_i,x_{n+j}]
z^{(i,j)}_{\sigma''(1)} \ldots z^{(i,j)}_{\sigma''(n+m - i - j)}) 
\end{align*}
for any $x_1,\ldots,x_{n+m}$ in $\G$, where we set 
$$
(y^{(i,j)}_1,\ldots,y^{(i,j)}_{i+j-2}) = (x_1,\ldots,x_{i-1},
x_{n+1},\ldots,x_{n+j-1}),
$$  
$(z^{(i,j)}_1,\ldots,z^{(i,j)}_{i+j-2}) = (x_{i+1},\ldots,x_{n},
x_{n+j+1},\ldots,x_{n+m})$ and $\SG_{i,j}$ is the set of shuffle 
permutations of $\SG_{i+j}$ associated to the partition 
$\{1,\ldots,i+j\} = \{1,\ldots,i\}\cup \{i+1,\ldots,i+j\}$.  

The map $m_1$ defines a Poisson structure on the commutative 
algebra $\Sh(\G)$. 
Moreover, this Poisson structure is compatible with $\Delta_{\Sh(\G)}$, 
so $(\Sh(\G),m_0,m_1,\Delta_{\Sh(\G)},$ $\varepsilon_{\Sh(\G)},S_{\Sh(\G)})$
is a commutative Hopf-Poisson algebra; moreover, the map $\G\mapsto 
(\Sh(\G),m_0,m_1,\Delta_{\Sh(\G)})$ defines in a natural way a
functor from the category of Lie algebras to that of Hopf-Poisson 
algebras. 

In the same way, the tensor algebra $T(V)$ is equipped with a unique
Hopf structure $(T(V),m_{T(V)},\Delta_{T(V)},\ve_{T(V),S_{T(V)}})$ 
such that the elements of $T^1(V)$ are 
primitive. If $V$ is the underlying space of a Lie coalgebra
$(\HH,\delta_\HH)$, then there exists a unique co-Poisson map 
$\delta_{T(\HH)} : T(\HH) \to T(\HH)^{\otimes 2}$, whose
restriction to $T^1(\HH)$ coincides with $\delta_\HH$. 
$(T(\HH),m_{T(\HH)},\Delta_{T(\HH)},\delta_{T(\HH)})$ is then 
a Hopf-co-Poisson algebra, and the map $\HH\mapsto 
(T(\HH),m_{T(\HH)},\Delta_{T(\HH)},\delta_{T(\HH)})$ is 
a functor from the category of Lie coalgebras to that of 
Hopf co-Poisson algebras.

\subsection{Quantizations of shuffle algebras}

Let  $\omega = (B_{pq})_{p,q\geq 0}$ belong to $\cB(\KK[[\hbar]])$. If 
$\G$ is a Lie algebra over $\KK$, define $\Sh^\omega(\G)$
as follows. As a vector space, $\Sh^\omega(\G)$ is isomorphic to 
$T(\G)[[\hbar]]$. If $x_1,\ldots,x_n$ belong to $\G$, let us denote by 
$(x_1 \ldots x_n)$ the element corresponding to $x_1\otimes\cdots\otimes x_n$. 
Define a map $m_{\Sh^\omega(\G)} : \Sh^\omega(\G)^{\wh\otimes 2}\to \Sh^\omega(\G)$ 
by the rule
\begin{align} \label{main}
& m_{\Sh^\omega(\G)} \big((x_1\ldots x_n), (y_1\ldots y_m) \big)
= \sum_{k\geq 0} \sum_{(p_1 \ldots p_k) \on{\ partition \ of\ }n, 
(q_1 \ldots q_k) \on{\ partition \ of\ }m}
\\ & \nonumber 
\big( B_{p_1q_1}(x_1\ldots x_{p_1}|y_1\ldots y_{q_1}) 
\ldots 
B_{p_kq_k}(x_{p_1 + \cdots + p_{k-1}+1}\ldots x_{p}
|y_{q_1 + \cdots + q_{k-1}+1}\ldots y_{q}) \big) 
\end{align}
for any $x_1,\ldots,y_m$ in $\G$ ($\wh\otimes$ is the $\hbar$-adically completed 
tensor product)
and $\Delta_{\Sh^\omega(\G)}$ as the unique linear map from  $\Sh^\omega(\G)$
to $\Sh^\omega(\G)\wh\otimes \Sh^\omega(\G)$ such that 
$$
\Delta_{\Sh^\omega(\G)}((x_1 \ldots x_n))
= \sum_{i = 0}^n (x_1\ldots x_i) \otimes (x_{i+1}\ldots x_n),  
$$
for any $x_1,\ldots,x_n$ in $\G$. 
Let us define $\ve_{\Sh^\omega(\G)}$ as the unique linear map from
 $\Sh^\omega(\G)$ to $\KK$, such that the restriction of $\ve_{\Sh^\omega(\G)}$ to 
$\G^{\otimes 0} = \KK$ is the identity map, and the restriction of 
$\ve_{\Sh(\G)}$ to $\oplus_{i>0}\G^{\otimes i}$ is zero.  
 
Finally, let us define inductively $S_{\Sh^\omega(\G)}$ as the unique 
endomorphism  of $\Sh^\omega(\G)$ such that 
$$
S_{\Sh^\omega(\G)}(1) = 1, \quad 
S_{\Sh^\omega(\G)}((x_1\ldots x_n)) = - \sum_{i=0}^{n-1} 
S_{\Sh^\omega(\G)}((x_1\ldots x_i))(x_{i+1}\ldots
x_n)
$$ 
for any $x_1,\ldots,x_n$ in $\G$. One checks that 
$$
S_{\Sh^\omega(\G)}(x_1\ldots x_n) = 
\sum_{k\geq 0} (-1)^k \sum_{(p_1,\ldots,p_k)
\on{\ partition\ of\ }n} (x_1\ldots x_{p_1}) \cdots 
(x_{p_1 + \cdots + p_{k-1} + 1}\ldots x_{p_1 + \cdots + p_k}). 
$$

\begin{prop}
$(\Sh^\omega(\G),m_{\Sh^\omega(\G)},\Delta_{\Sh^\omega(\G)},\ve_{\Sh^\omega(\G)},
S_{\Sh^\omega(\G)})$
is a Hopf algebra. If $\omega$ belongs to $\cB(\KK)_{1/2}$, then 
the subspace of $\Sh^\omega(\G)$ of all symmetric tensors 
is a Hopf subalgebra of $\Sh^\omega(\G)$, canonically isomorphic with 
$U\G$. Moreover, the mapping $\G\mapsto\Sh^\omega(\G)$ defines a functor from the 
category of Lie algebras to that of Hopf algebras. 

If $\omega$ and $\omega'$ belong to the same orbit of the action of 
$\cG(\KK[[\hbar]])$ on $\cB(\KK[[\hbar]])$, then  
the Hopf algebras $\Sh_\omega(\G)$ and $\Sh_{\omega'}(\G)$
are canonically isomorphic (that is, the isomorphisms are 
depend functorially on $\G$). If $\varpi$ belongs to $\SSh(\KK)$, 
we will denote by $\Sh^\varpi(\G)$ any of the Hopf algebras 
$\Sh^\omega(\G)$, if $\varpi$ is the coset of $\omega$.  
\end{prop}

{\em Proof.} $\Delta_{\Sh^\omega(\G)}$ is obviously coassociative. Moreover, 
it is also clear that $m_{\Sh^\omega(\G)}$ is a coalgebra map. The associativity 
of $m_{\Sh^\omega(\G)}$ follows from the identities (\ref{cond:Bpq}).
The fact that if $x$ and $y$ belong to $\Sh(\G)$, then 
$S_{\Sh(\G)}(xy) = S_{\Sh(\G)}(x)S_{\Sh(\G)}(y)$
is proved by induction on the degree of $x$ any $y$. 
The other Hopf algebra axioms are checked directly. 
This proves the first statement.

The form of the product implies that the symmetric tensors form a subalgebra
of $\Sh^\omega(\G)$. Let $\iota$ be the algebra morphism from $T(\G)$
to $\Sh^\omega(\G)$ such that for any $x$ in $\G$, $\iota(x) = (x)$. Then 
$\iota(x \otimes y - y \otimes y - [x,y])$ is equal to zero, so 
$\iota$ induces an algebra morphism $\bar\iota$ from $U\G$ to the subalgebra of
summetric elements of $\Sh^\omega(\G)$. The associated graded of this morphism 
coincides with the identity map of $S(\G)$, therefore $\bar\iota$ is an 
isomorphism. 

Assume that $\omega = P * \omega'$, where $P = (P_i)_{i\geq 1}\in \cG(\KK[[\hbar]])$, 
then the canonical isomorphism $i_{\omega,\omega'}$ 
from $\Sh^\omega(\G)$ to $\Sh^{\omega'}(\G)$
is given by by the rule
$$
i_{\omega,\omega'}(x_1 \ldots \ldots x_n) = \sum_{\al>0}
\sum_{(n_1,\ldots,n_\al)\in \Part_\al(n)}
\big( P_{n_1}(x_1,\ldots,x_{n_1})\cdots 
P_{n_\al}(x_{\sum_{\al' = 1}^{\al - 1} n_{\al'} +1}\ldots x_n)\big). 
$$
\hfill \qed\medskip 

Let us explain why we consider $\Sh^\omega(\G)$ as a quantization of the 
Hopf-Poisson algebra $(\Sh(\G),m_0,m_1,\Delta_{\Sh(\G)},\ve_{\Sh(\G)},S_{\Sh(\G)})$. 
Let us denote by $\Sh^{\omega,\leq i}(\G)$ the subspace of $\Sh^\omega(\G)$
equal to $\oplus_{i'\leq i} T^{i'}(\G)[[\hbar]]$. Then the sequence of 
inclusions $\Sh^{\omega,\leq 0}(\G) $ $\subset \Sh^{\omega,\leq 1}(\G) \subset
\cdots$ is a Hopf algebra filtration of $\Sh^\omega(\G)$. The associated 
graded algebra is commutative and inherits therefore a Poisson bracket. 
If $\omega$ belongs to $\cB(\KK[[\hbar]])_{\la}$, where $\la\in {1/2} + 
o(\hbar)$, the reduction modulo $\hbar$ of this Poisson algebra is   
$(\Sh(\G),m_0,m_1)$.

\subsection{$\SSh(\KK)$ and quantizations of tensor algebras}

Let $(\A,\delta_\A)$ be a Lie coalgebra. Let us consider the 
following quantization problem:  to construct a topologically 
free $\KK[[\hbar]]$-Hopf algebra, 
$(T,m_T,\Delta,\ve_T,S_T)$ quantizing the Hopf-co-Poisson structure 
$(T(\A),m_{T(\A)},\Delta^0_{T(\A)},\delta_{T(\A)},\ve_{T(\A)},S_{T(\A)})$
(we call $(T,m_T,\Delta_T,\ve_T,S_T)$ a quantized tensor algebra). 

The corresponding functorial problem 
is to construct all functors from the category of 
Lie coalgebras to that of quantized tensor algebras with the suitable
classical limit. Let us say that such a functor is universal if 
the structure constants of the quantization of $T(\A)$ depends polynomially 
in those of $\A$. 
 The purpose of this Section is to construct a bijection between 
$\SSh(\KK)$ and $\{$universal quantization functors of the tensor algebras 
$T(\G)\}$.    

\begin{prop} \label{prop:lezard}
Let $(\A,\delta_\A)$ be a Lie coalgebra, and let 
$(T,m_T,\Delta_T,\ve_T,S_T)$ be a quantization of 
$(T(\A),m_{T(\A)},\Delta^0_{T(\A)},\delta_{T(\A)},
\ve_{T(\A)},S_{T(\A)})$. 

Then there exists an algebra isomorphism
$\theta : T\to T(\A)[[\hbar]]$, such that if 
$\Delta: T(\A)[[\hbar]] \to T(\A)^{\otimes 2}[[\hbar]]$ is the map 
$(\theta\otimes\theta)\circ\Delta_T \circ\theta^{-1}$, then   
$\Delta(T^i(\A)) \subset \wh\bigoplus_{p,q|p+q\geq i}
\hbar^{p+q - i} T^p(\A)\otimes T^q(\A)$
($\wh\oplus$ is the completed direct sum).  

Moreover, any isomorphism with this 
property is of the form $\kappa\circ\theta$, where 
$\kappa$ is the automorphism of $T(\A)[[\hbar]]$
induced by a map $\A\to\wh\oplus_{n\geq 1}
\hbar^{n-1} T^n(\A)[[\hbar]]$ whose reduction modulo $\hbar$
is the identity. 
\end{prop}

{\em Proof.} Let $\phi : \A\to T$ be any section of the projection map 
$T\to T/\hbar T = T(\A)\to \A$ (the second map is the projection 
on $\A$ parallel to $\oplus_{i\neq 1} T^i(\A)$). This map 
extends to a unique algebra morphism $\iota : T(\A)[[\hbar]] \to T$. Since $T$
is a flat deformation of $T(\A)$, $\iota$ is 
a linear isomorphism. Let us set $\Delta_T^\iota = 
(\iota\otimes\iota)\circ\Delta_T\circ\iota^{-1}$. 
Then $\Delta_T^\iota$ is a coproduct map on 
$T(\A)[[\hbar]]$. We will now construct an algebra automorphism
$\phi$ of $T(\A)[[\hbar]]$, such that the reduction of $\phi$  
mod $\hbar$ is the identity and $(\phi\otimes\phi)
\circ \Delta_T^\iota\circ\phi^{-1}$ has the required property. 

Assume that we constructed an automorphism $\phi_{(n)}$ 
of $T(\A)[[\hbar]]$, whose reduction modulo $\hbar$ is 
the identity, and such that $\Delta_{T,n} = (\phi_{(n)}
\otimes\phi_{(n)})\circ\Delta_T^\iota
\circ \phi_{(n)}^{-1}$ has the property
\begin{equation} \label{lezard}
\Delta_{T,n}(\A) \subset \bigoplus_{p,q|p+q\geq 1}
\hbar^{p+q - 1} T^p(\A) \otimes T^q(\A) +
 \hbar^n T(\A)\otimes T(\A)[[\hbar]] . 
\end{equation}
Let $\Delta_{T,n,k}$ be the linear endomorphisms of $T(\A)$
such that $\Delta_{T,n} = \sum_{k\geq 0}\hbar^k \Delta_{T,n,k}$. 
If we denote by $\Delta_{T(\A)}$ the usual (undeformed) 
coproduct of $T(\A)$, the coassociativity of $\Delta_{T,n}$
implies that 
$(\Delta_{T(\A)}\otimes\id - \id\otimes \Delta_{T(\A)})\circ 
\Delta_{T,n,n}$ maps $T(\A)$ to $\oplus_{p,q,r|p+q+r\leq n+1}
T^p(\A) \otimes T^q(\A)\otimes T^r(\A)$. 

Let us denote by $\Delta_{>n+1}$ the composition of
$(\Delta_{T,n,n})_{|\A}$ with the projection of $T(\A)^{\otimes 2}$
on $\oplus_{p,q|p+q>n+1} T^p(\A)\otimes T^q(\A)$ 
parallel to $\oplus_{p,q|p+q\leq n+1} T^p(\A)\otimes T^q(\A)$. 
Then $(\Delta_{T(\A)}\otimes \id  - \id\otimes \Delta_{T(\A)})
\circ \Delta_{>n+1} = 0$.  
Let $\psi_n$  be the composition $ - \ve_{T(\A)} \circ \Delta_{>n+1}$; 
$\phi'_n$ is a linear map from $\A$ to $T(\A)$. 
Then the relation
$$
\id_{T(\A)\otimes T(\A)} = 
\Delta_{T(\A)} \circ (\ve_{T(\A)}) 
+ (\ve\otimes \id_{T(\A)}^{\otimes 2}) \circ (\Delta_{(T(\A))}
\otimes \id_{T(\A)} - \id_{T(\A)}\otimes \Delta_{(T(\A))})
$$ 
implies that $\Delta_{>n+1} = \Delta_{T(\A)}\circ \psi_n$. 
Let us define 
$\phi'_n$ as the automorphism of $T(\A)[[\hbar]]$ induced 
by the map $1 - \hbar^n \psi_n$, and set 
$\phi_{(n+1)} = \phi'_n\circ\phi_{(n)}$, then $(\phi_{(n+1)}\otimes\phi_{(n+1)})
\circ\Delta_{T,n} \circ (\phi_{(n+1)})^{-1}$ satisfies (\ref{lezard})
with $n$ replaced by $n+1$. 

The sequence $(\phi_{(n)})_{n\geq 0}$ has a unique $\hbar$-adic limit, 
wich we denote by $\phi$. Then if we set $\Delta_{T,\infty}
= (\phi_{(\infty)}\otimes \phi_{(\infty)})\circ 
\Delta_T^\iota \circ\phi_{(\infty)}^{-1}$, then we have 
$$
\Delta_{T,\infty}(\A) \subset \bigoplus_{p,q|p+q\geq 1}
\hbar^{p+q - 1} T^p(\A)\otimes T^q(\A) . 
$$
Since $\Delta_{T,\infty}$
is an algebra morphism, this relation implies 
$$
\Delta_{T,\infty}(T^k(\A)) \subset \bigoplus_{p,q|p+q\geq k}
\hbar^{p+q - k} T^p(\A) \otimes T^q(\A)  
$$
for any $k$. So we set $\theta = \phi\circ\iota$. 
\hfill \qed\medskip 

Assume that $(\A,\delta_\A)$ is a finite-dimensional (or
nonnegatively graded with finite dimensional components) Lie 
coalgebra. Then $\A^*$ (or the restricted dual of $\A$) is a 
Lie algebra; denote it by $\B$.  In the situation of 
Proposition \ref{prop:lezard}, 
define $\delta_{pq}^\B$ as the  
composition of the restriction $(\theta\otimes\theta)
\circ\Delta\circ\theta^{-1}$ with the projection of
$T(\A)\otimes T(\A)[[\hbar]]$ on $T^p(\A)\otimes T^q(\A)[[\hbar]]$
parallel to all other $T^{p'}(\A)\otimes T^{q'}(\A)[[\hbar]]$. 

Define $B^\B_{pq}$ as the linear
map from $\B^{\otimes p+q}$ to $\B[[\hbar]]$ dual to 
$\delta_{pq}^\B$. Equip $\Sh(\B)[[\hbar]]$ with its
usual coproduct and the product $m_{T,\theta}$ given by formula 
(\ref{main}), with the $B_{pq}$ replaced by $B_{pq}^{\B}$. 
  
\begin{lemma} \label{lemma:lez}
In the situation of Proposition \ref{prop:lezard}, 
there is a unique Hopf algebra structure 
$(\Sh(\B)[[\hbar]],m_{T,\theta},\Delta_{\Sh(\B)},\ve_{\Sh(\B)},
S_{T,\theta})$ on $\Sh(\B)[[\hbar]]$ with product $m_{T,\theta}$
and coproduct $\Delta_{\Sh(\B)}$. The pairing $T(\A)[[\hbar]]
\times \Sh(\B)[[\hbar]] \to \KK[[\hbar]][\hbar^{-1}]$ defined by 
$\langle a_1\otimes\cdots\otimes a_p, (b_1\ldots b_q)\rangle 
= {{\delta_{pq}}\over{\hbar^p}} \prod_{i = 1}^p \langle a_i,b_i
\rangle$ induces 
a Hopf pairing between this Hopf algebra and the Hopf algebra structure on 
$T(\A)[[\hbar]]$ defined in Proposition \ref{prop:lezard}.     
\end{lemma}

We now construct a map $\gamma_\KK : \SSh(\KK) \to\{$universal 
quantization
functors of the tensor algebras $T(\A)\}$. If $n$ is an 
integer and $P$ belongs to $FL_n(\KK)$, and if $(\A,\delta_\A)$ is a 
Lie coalgebra, let us define the map $\delta_\A^{(P)} : \A\to
\A^{\otimes n}$ as follows. Let $FA_n(\KK)$ be the 
part of the free algebra with $n$ generators, homogeneous in each 
generator. There is a unique element $(P_\sigma)_{\sigma\in\SG_n}$
in $\KK^{\SG_n}$ such that $P = \sum_{\sigma_in\SG_n}
P_\sigma x_{\sigma(1)} \cdots x_{\sigma(n)}$. Let us set 
$$
\delta^{(P)}_\A(a) = 
{1\over{n}}
\big( (\id_{\A}^{\otimes n}\otimes\delta_\A)
\circ \cdots\circ (\id_\A\otimes\delta)\circ\delta(a) 
\big)^{(\sigma(1)\ldots
\sigma(n))},
$$ 
for any $a$ in $\A$. Then Proposition \ref{reut} implies that 
if $\B$ is any Lie algebra dual to $\A$, then 
$$
\langle \delta_\A^{(P)}(a), b_1\otimes\cdots \otimes b_n
 \rangle_{\A^{\otimes n}\times\B^{\otimes n}}
= 
\langle a, P(b_1,\ldots,b_n)
\rangle_{\A^{\otimes n}\times\B^{\otimes n}} 
$$
(here $\langle\ ,\  \rangle_{\A\times\B}$ is the 
pairing between $\A$ and $\B$ and 
$\langle\ ,\  \rangle_{\A^{\otimes n}\times\B^{\otimes n}}$ is 
its $n$th tensor power). If $P$ belongs to $\KK[[\hbar]]$, the
map $\delta^{(P)} : \A\to \A^{\otimes n}[[\hbar]]$ is defined in the 
same way. 

If $\omega$ belongs to $\cB(\KK[[\hbar]])$, and $\A$ is any 
Lie coalgebra, let us define  
$T_\hbar^\omega(\A)$ 
as follows. As an algebra, $T_\hbar^\omega(\A)$  is isomorphic 
to $T(\A)[[\hbar]]$. Define $\mu_{pq}$ as the  
map from $\A^{\otimes p+q}$ to $T(\A)\otimes T(\A)$
obtained by the composition of the isomorphism $\A^{\otimes (p+q)}
\to \A^{\otimes p}\otimes \A^{\otimes q}$ with the injections 
of $\A^{\otimes p}$ and $\A^{\otimes q}$ in each factor of 
$T(\A)\otimes T(\A)$. Then there exists a unique coalgebra
map $\Delta_{T_\hbar^\omega(\A)}$ on $T(\A)[[\hbar]]$, such that for any $a$ in $\A$,
$$
\Delta_{T_\hbar^\omega(\A)}(a) = \sum_{p,q\geq 0} \hbar^{p+q - 1}
\mu_{p,q}(\delta^{(B_{pq})}(a)). 
$$ 
The usual augmentation map $\ve_{T_\hbar^\omega(\A)}$ is a counit for this coalgebra 
structure, and there exists a unique corresponding antipode
$S_{T_\hbar^\omega(\A)}$. $(T_\hbar^\omega(\A),m_{T_\hbar^\omega(\A)},
\Delta_{T_\hbar^\omega(\A)},\ve_{T_\hbar^\omega(\A)},S_{T_\hbar^\omega(\A)})$ 
is then a quantized tensor algebra. It is easy to see
that it is a quantization of $(T(\A),m_{T(\A)}\Delta_{T(\A)},
\delta_{T(\A)})$. Moreover, for any fixed $\omega$, the map 
$\A\mapsto T_\hbar^\omega(\A)$ is a quantization functor of the 
tensor algebras $T(\A)$, and if $\omega$ and $\omega'$ are in the same
$\cG(\KK[[\hbar]])$-orbit, then there are functorial (in $\A$)
Hopf algebra isomorphisms between $T_\hbar^{\omega}(\A)$ and 
$T_\hbar^{\omega'}(\A)$. If $\varpi$ belongs to $\SSh(\KK)$, let us 
define $(\A\mapsto T_\hbar^\varpi(\A))$ as any of the quantization 
functors $(\A\mapsto T_\hbar^\omega(\A))$, where $\varpi$ is the 
coset of $\omega$.   Summarizing, we  have 

\begin{prop} \label{prop:Ta}
The assignment $\varpi\mapsto (\A\mapsto T_\hbar^\varpi(\A))$ defines
a map $\gamma_\KK$ from $\SSh(\KK)$ to $\{$universal quantization functors
of the tensor algebras $T(\A)\}$. 
\end{prop}

We will now show

\begin{prop} \label{prop:quant:functors} 
The restriction to $\SSh(\KK)_{1/2}$  of the 
map $\gamma_\KK$ is a bijection. It is functorial in $\KK$
and equivariant with respect to the natural actions of 
$\on{Aut}(\KK[[\hbar]])$. 
\end{prop}

{\em Proof.} Let us construct the inverse map to $\gamma_\KK$. 

Let $\cF_n$ be the free Lie algebra with $n$ generators 
$x_1,\ldots,x_n$. 
If $\Psi$ is a quantization functor of the tensor algebras $T(\A)$, 
we may apply  Proposition \ref{prop:lezard} and Lemma \ref{lemma:lez}
to $\Psi(\cF_n^*)$. These statements have analogues in the situation
of $\NN^n$-graded Lie algebras. 
Functoriality  with respect to the Lie algebra automorphisms
of $\cF_n$ given by $x_i\mapsto \la_i x_i$ allows us to apply these refined
statements to $\cF_n$. We obtain maps $B^{\cF_n}_{pq} : (\cF_n)^{\otimes p+q}
\to \cF_n[[\hbar]]$, graded and well defined up to graded automorphisms. 
The maps $B^{\cF_n}_{pq}$ are graded, so 
if $p+q\leq n$, then $B^{\cF_n}_{pq}(x_1,\ldots,x_{p+q})$ belongs to 
$FL_{p+q}(\KK[[\hbar]])$. The associativity of the product of 
$\Psi(\cF_n)$ implies that the $(B^{\cF_n}_{pq})_{p,q|p+q\leq n}$
satisfy identities (\ref{cond:Bpq}), for $p+q+r\leq n$. The 
family $(B^{\cF_n}_{pq})_{p,q|p+q\leq n}$ therefore 
defines an element $\wt\Ts_n(\Psi)$ of $\SSh_n(\KK)$. 
Functoriality with respect to the map $\cF_{n}\to\cF_{n+1}$ 
sending each $x_i$ to $x_i$ implies that the family 
$(\wt\Ts_n(\Psi))_{n\geq 0}$ defines an element $\wt\Ts(\Psi)$
of $\SSh_n(\KK)$. Then there is a unique $\la(\Psi) \in 1 + o(\hbar)$ 
such that $\wt\Ts(\Psi)$ belongs to $\SSh_n(\KK)_{\la(\Psi)/2}$. Let 
$\Ts(\Psi)$ be the result of the action of $\la(\Psi)^{-1}\in\KK[[\hbar]]^\times$
on $\wt\Ts(\Psi)$. 
We have obviously $\Ts\circ \gamma_\KK = \id_{\SSh(\KK)_{1/2}}$. 

Let us prove that 
$\gamma_\KK  \circ \Ts= \id_{\SSh(\KK)_{1/2}}$. Let $\Psi$ be 
a universal quantization functor of the tensor algebras $T(\A)$. Apply 
Proposition \ref{prop:lezard} and Lemma \ref{lemma:lez} to each 
Lie coalgebra $\A$. Then for 
each Lie algebra $\G$ over $\KK$, we obtain  
linear maps $B_{pq}^{\G} : \G^{\otimes p+q}\to\G[[\hbar]]$
satisfying the relations (\ref{cond:Bpq}), depending
polynomially in the structure constants of $\G$, and such that 
the dual to $\Psi(\A)$ is the quantized shuffle algebra of $\A^*$
corresponding to the maps $B_{pq}^\G$. One the other hand, 
such linear maps $\G^{\otimes p+q}\to \G$ all correspond to 
substitution of elements of $\G$ in free Lie polynomials.
The family $(B_{pq}^\G)_{p,q\geq 0}$ therefore corresponds to an element of 
$\SSh(\KK)$. 
\hfill \qed\medskip

\subsection{Relation with the CBH formula}
\label{sect:cort}

Let $\G$ be a Lie algebra over $\KK$ and let 
$\Sym : S\G\to U\G$ be the symmetrization map, such that 
for any $x$ in $\G$ and any integer $n$, $\Sym(x^n) = x^n$. 
The pull-back of the product on the universal enveloping algebra 
$U\G$ by the symmetrization map 
$\Sym$ may be viewed as a star-product over 
$S\G$. This pull-back involves the 
Campbell-Baker-Hausdorff (CBH) series and 
is called PBW quantization of $S\G$.  

Let us recall how the CBH series is defined. 
It is an infinite series $B^{CBH}(x,y) = \sum_{p,q\geq 0} B^{CBH}_{pq}(x,y)$, 
where $B^{CBH}_{pq}(x,y)$ is an element of the free Lie algebra in two variables 
$x,y$, homogeneous of degrees $p$ in $x$ and $q$ in $y$. For $x,y$ 
formal variables in the neighborhood of zero in $\G$, the identity 
$$
e^x e^y = e^{B^{CBH}(x,y)}
$$ 
holds in the formal Lie group associated with $\G$.  

Then if $x$ and $y$ belong to $\G$,  we have the identity 
\begin{equation} \label{pouillard}
x^p y^q = \sum_{k\geq 0}
{{p!q!}\over {k!}} \sum_{(p_1,\ldots,p_k)\in \Part_k(p),  
(q_1,\ldots,q_k)\in \Part_k(q) }
\Sym ( B^{CBH}_{p_1q_1}(x,y)  \cdots  B^{CBH}_{p_kq_k}(x,y)) 
\end{equation}
in $U\G$. 

Assume that $\omega = (B_{pq})_{p,q\geq 0}$ belongs to $\cB(\KK[[\hbar]])_{1/2}$. 
Since the image of $x^p$ by the canonical map $U\G\to 
\Sh_\hbar^\omega(\G)$ is $p!(x\ldots x)$ ($x$ appears $p$ times), 
and identifying the terms of lower degree of the image of 
(\ref{pouillard}) by $U\G\to \Sh_\hbar^\omega(\G)$, we get 
$B_{pq}(x,\ldots,x|y,\ldots,y) = B^{CBH}(x,y)$. This identity
holds for any Lie algebra, in particular, if $\G$ is the free Lie 
algebra with two generators and $x,y$ are identified with these
generators. We have therefore 

\begin{prop}
If $\omega = (B_{pq})_{p,q\geq 0}$ belongs to 
$\cB(\KK[[\hbar]])_{1/2}$, then the identities 
$$
B_{pq}(x,\ldots,x|y,\ldots,y) = B_{pq}^{CBH}(x,y)
$$
hold in the free Lie algebra with two generators $x,y$;  
\end{prop}

\subsection{Explicit formulas} \label{sect:explicit}

One may show that if $\omega$ is any element 
in $\cB(\KK[[\hbar]])_{1/2}$, there exists $\omega' = 
(B_{pq})_{p,q\geq 0}$ in the orbit $\cG(\KK[[\hbar]]) * \omega$
of $\omega$ such that 
$$
B_{12}(x,x'|y) = {1\over {24}} ([x,[x',y]] + [x',[x,y]]), \quad 
B_{21}(x|y,y') = {1\over {24}}([y,[y',x]] + [y',[y,x]]). 
$$
In the algebra $\Sh^{\omega'}(\G)$, we have the relations
$$
(x)(y) = (xy) + (yx) + {1\over 2}([x,y]), 
$$
$$
(xx')(y) = (xx'y) + (xyx') + (yxx') + {1\over 2}(x[x',y])  
+ {1\over 2}([x,y]x') + {1\over{24}} ([x,[x',y]] + [x',[x,y]]),  
$$
$$
(x)(yy') = (xyy') + (yxy') + (yy'x) + {1\over 2}([x,y]y')  
+ {1\over 2}(y[x,y']) + {1\over{24}} ([y,[y',x]] + [y',[y,x]]).   
$$

\newpage

\section{$R$-matrices in shuffle algebras}
\label{sect:R:mat}

Let us fix an element $\omega$ of $\cB(\KK[[\hbar]])$. 
Let $\A$ and $\B$ be Lie algebras and let $r_{\A\B} = \sum_{i\in I}
a_i\otimes b_i$ be an element of 
$\A\otimes \B$. In this Section, we associate to $r_{\A\B}$
a family of elements $(\cR^\omega_n(r_{\A\B}))_{n\geq 0}$ of 
$\Sh^\omega(\A) \otimes\Sh^\omega(\B)$,
satisfying analogues of the quasitriangular identities $(\Delta\otimes
\id)(\cR) = \cR^{(13)} \cR^{(23)}$    and  $(\id \otimes\Delta)(\cR) =
\cR^{(13)} \cR^{(12)}$. 

If $n_1,\ldots,n_k$ are positive integers, denote by 
$x_1^{(i)},\ldots,x^{(i)}_{n_i}$ the generators of the 
$i$th factor of the tensor product $
(\otimes_{i = 1}^k FL_{n_i})
\otimes FL_{n_1 + \cdots + n_k}$ and by 
$y_1,\ldots,y_{n_1 + \cdots + n_k}$
the generators of the last tensor factor. Then 
the product of
symmetric groups $\prod_{i = 1}^k \SG_{n_i}$ acts on this tensor 
product, is such a way that $\SG_i$ permutes simultaneously 
the generators $(x^{(i)}_a)_{a = 1,\ldots,k_i}$ and 
$(y_{n_1 + \cdots + n_{i-1} + a})_{a = 1,\ldots,k_i}$. 
Let us denote by $F_{n_1,\ldots,n_k}$ the space 
of coinvariants of this action, so 
$$
F_{n_1,\ldots,n_k} = \big( 
(\bigotimes_{i = 1}^k FL_{n_i})
\otimes FL_{n_1 + \cdots + n_k}
\big)_{\prod_{i = 1}^k \SG_{n_i}} . 
$$
Let us define inductively $\la_{n_1,\ldots,n_k} \in F_{n_1,\ldots,n_k}$
as follows. 
Recall that if $p$ and $r$ are integers, $\Part_r(p)$ 
is the set of $r$-partitions of $p$, 
i.e.\ the set of families of integers $(p_1,\ldots,p_r)$
such that $p_1 + \cdots + p_r = p$. 
We set $\la_1 = x^{(1)}_1\otimes y_1$, and if $\sum_{l = 1}^k n_l >1$
\begin{align*}
& \la_{n_1,\ldots,n_k} = {1\over{n_1 + \cdots + n_k - 1}}
\sum_{\al_1,\ldots,\al_k\in\NN_{>0}}
\sum_{l = 1}^{k-1}
\\ & 
\sum_{(n_{ij})_j\in \Part_{\al_i}(n_i)}
\la''_{(n_{1i})_i} \otimes \cdots \otimes \la''_{(n_{ki})_i} 
\otimes
B^{\omega^\vee}_{\al_1 + \cdots + \al_l, \al_{l+1} + \cdots + \al_k}
\big(\la_{(n_{1j})_j}^{\prime (1)}, \ldots ,
\la_{(n_{kj})_j}^{\prime (\al_k)}
\big) , 
\end{align*}
where we set 
$$
\la_{(n_{ij})_{j = 1,\ldots,\al_i}} = 
\la^{\prime (1)}_{(n_{ij})_{j = 1,\ldots,\al_i}} 
\otimes\cdots \otimes 
\la^{\prime (\al_i)}_{(n_{ij})_{j = 1,\ldots,\al_i}} 
\otimes \la''_{(n_{ij})_{j = 1,\ldots,\al_i}} , 
$$
and the variables 
$(x^{(i)}_{a})_{a = 1,\ldots,n_i}$ (resp., 
$y_{ \sum_{\al <i} n_\al + \sum_{\beta <j} n_{i\beta} + 1}, 
\cdots, y_{\sum_{\al <i} n_\al + \sum_{\beta \leq j} n_{i\beta}}$)
are substituted in $\la''_{(n_{ij})_j}$ (resp., in 
$\la^{\prime (j)}_{(n_{ij})_j}$). 

Let us denote by $\conc_{\Sh^\omega(\B)}$ the concatenation product on the tensor 
algebra $T(\B)[[\hbar]]$, which we identify linearly with 
$\Sh^\omega(\B)$. We have 
$\conc_{\Sh^\omega(\B)}((x_1\ldots x_k)\otimes (x_{k+1}\ldots x_n)) = 
(x_1\ldots x_n)$. 

There is a unique map $\kappa(r_{\A\B})$ from $F_{n_1,\ldots,n_k}$
to $\Sh^\omega(\A)\otimes \B$, such that for any $P_1\in FL_{n_1}$, 
$\ldots, P_k\in FL_{n_k}$, $Q\in FL_{n_1 + \cdots + n_k}$, we have 
\begin{align*}
&
\kappa(r_{\A\B})(\otimes_{i = 1}^k P_i(x^{(i)}_a) \otimes 
Q(y_1,\ldots,y_{n_1 + \cdots + n_a})) 
= \sum_{i_1,\ldots i_{n_1 + \cdots + n_k\in I}}
\\ & 
\big( P_1(a_{i_1},\ldots,a_{i_{n_1}}) \ldots 
P_k(a_{i_{n_1 + \cdots + n_{k-1} + 1}},\ldots, a_{
i_{n_1 + \cdots + n_k}}) \big)\otimes 
Q(b_{i_1},\ldots,b_{i_{n_1 + \cdots + n_k}}) . 
\end{align*}
We also denote by $\kappa(r_{\A\B})$ the unique 
extension of this map to a linear map from $\oplus_{n_1,\ldots,n_k}
F_{n_1,\ldots,n_k}$ to $\Sh^\omega(\A)\otimes \B$. 

Let us denote by $\Psi_{\Sh(\B)}$ the linear automorphism of 
$\Sh(\B)$ such that  for any $x_1,\ldots,x_n$ in $\B$, 
$\Psi_{\Sh(\B)}((x_1\ldots x_n)) = (x_n\ldots x_1)$. 

\begin{lemma} \label{Psi}
If $\Delta'_{\Sh(\B)}$ is the composition of $\Delta_{\Sh(\B)}$ with the   
permutation of factors, then $\Psi_{\Sh(\B)}$ is a Hopf algebra 
isomorphism from $(\Sh^{\omega^\vee}(\B),m_{\Sh^{\omega^\vee}(\B)},
\Delta_{\Sh^{\omega^\vee}(\B)})$ 
to $(\Sh^\omega(\B),m_{\Sh^\omega(\B)},\Delta'_{\Sh^\omega(\B)})$, so 
$m_{\Sh^\omega(\B)} \circ (\Psi_{\Sh(\B)} \otimes \Psi_{\Sh(\B)}) = 
\Psi_{\Sh(\B)} \circ m_{\Sh^{\omega^\vee}(\B)}$ and 
$\Delta'_{\Sh^\omega(\B)}\circ \Psi_{\Sh(\B)} = 
(\Psi_{\Sh(\B)} \otimes \Psi_{\Sh(\B)}) \circ \Delta_{\Sh^{\omega^\vee}(\B)}$.  
\end{lemma}

If $m$ is an associative
operation, we will denote by  $m^{(k)}$ the $k$th fold product associated to  
$m$.  If $n$ is a positive integer, we set $\la_n = \sum_{k>0}
\sum_{(n_1,\ldots,n_k)\in \Part_k(n)} \la_{n_1,\ldots,n_k}$.

\begin{prop} \label{prop:QT}
Let us define the family $(\cR'_n(r_{\A\B}))_{n\geq 0}$ of elements
of $\Sh^\omega(\A)\otimes\Sh^{\omega^\vee}(\B)$ as follows. 
$\cR'_0(r_{\A\B}) = 1$, 
$\cR'_1(r_{\A\B}) = r_{\A\B}$ and 
\begin{align*}
& \cR'_n(r_{\A\B}) 
\\ & = \sum_{k \geq 1} \sum_{(n_1,\ldots,n_k)\in \Part_k(n)}
\big( m_{\Sh(\A)}^{(k)} \otimes\conc_{\Sh(\B)}^{(k)} \big) 
\big( \kappa(r_{\A\B})(\la_{n_1})^{(1,k+1)} \cdots  
\kappa(r_{\A\B})(\la_{n_k})^{(k,2k)} \big) . 
\end{align*}
If we set 
$$
\cR^\omega_n(r_{\A\B}) = (\id\otimes \Psi_{\Sh(\B)})(\cR'_n(r_{\A\B})), 
$$
then the family $(\cR_n(r_{\A\B}))_{n\geq 0}$ satisfies 
\begin{align} \label{ids:Delta}
& (\Delta_{\Sh^\omega(\A)}\otimes \id)(\cR^\omega_n(r_{\A\B})) = 
\sum_{k = 0}^n (\cR^\omega_k(r_{\A\B}))^{(13)} 
(\cR_{n - k}^\omega(r_{\A\B}))^{(23)}, 
\\ & \nonumber 
(\id \otimes \Delta_{\Sh^{\omega}(\B)})(\cR_n(r_{\A\B})) = 
\sum_{k = 0}^n (\cR^\omega_k(r_{\A\B}))^{(13)} 
(\cR_{n - k}^\omega)(r_{\A\B})^{(12)}.  
\end{align}
Moreover, we have 
\begin{equation} \label{behavior:S}
(S_{\Sh^\omega(\A)} \otimes \id_{\Sh^\omega(\B)}) (\cR^\omega_n(r_{\A\B})) 
= (\id_{\Sh^\omega(\A)}\otimes S^{-1}_{\Sh^\omega(\B)}) (\cR^\omega_n(r_{\A\B})) . 
\end{equation}
\end{prop}

{\em Proof.}
It is immediate that $(\cR'_n(r_{\A\B}))_{n\geq 0}$ satisfies 
\begin{equation} \label{k:lo}
(\id \otimes \Delta_{\Sh^{\omega^\vee}(\B)})(\cR'_n(r_{\A\B})) = 
\sum_{k = 0}^n \cR_k(r_{\A\B})^{\prime (12)} 
\cR_{n - k}(r_{\A\B})^{\prime (13)},
\end{equation} 
so 
$(\cR'_n(r_{\A\B}))_{n\geq 0}$ satisfies  the second part of 
(\ref{ids:Delta}). Let us prove by induction on $m$ that the relation 
\begin{equation} \label{joas}
(\Delta_{\Sh^\omega(\A)}\otimes \id)(\cR'_m(r_{\A\B})) = 
\sum_{k = 0}^m \cR'_k(r_{\A\B})^{(13)} 
\cR'_{m - k}(r_{\A\B})^{(23)} 
\end{equation}
is also satisfied. The relation clearly holds for $m = 0,1$. 
Assume that it is satisfied up to order $m = n-1$. Let us set 
$$
Z = 
(\Delta_{\Sh^\omega(\A)}\otimes \id)(\cR'_n(r_{\A\B}))  -  
\sum_{k = 0}^n \cR'_k(r_{\A\B})^{(13)} 
\cR'_{n - k}(r_{\A\B})^{(23)}  
$$
and let us apply $(\id\otimes \Delta_{\Sh^{\omega^\vee}(\B)})$ to $Z$. We find 
\begin{align*}
&
(\id\otimes \Delta_{\Sh^{\omega^\vee}(\B)})(Z) = 
(\Delta_{\Sh^\omega(\A)}\otimes \Delta_{\Sh^{\omega^\vee}(\B)})
(\cR'_n(r_{\A\B}))  
\\ & 
-   \sum_{k = 0}^n \sum_{k',k''}
 \cR'_{k'}(r_{\A\B})^{(13)} \cR'_{k - k'}(r_{\A\B})^{(14)} 
\cR'_{k''}(r_{\A\B})^{(23)}  \cR'_{n - k - k''}(r_{\A\B})^{(24)}  , 
\end{align*}
and equation (\ref{k:lo}) and the induction hypothesis imply 
that this is $Z^{(123)} + Z^{(124)}$. Since the space of primitive
elements of $\Sh^{\omega^\vee}(\B)$ coincides with $\B$, 
$Z$ belongs to $\Sh^\omega(\A)^{\otimes 2}\otimes \B$. 

Let us show that $Z$ is zero. 
Let us set $\wt\Delta_{\Sh^\omega(\A)}(x) = \Delta_{\Sh^\omega(\A)}(x) - 
x\otimes 1 - 1\otimes x + \varepsilon_{\Sh^\omega(\A)}(x) (1\otimes 1)$. 
The restriction of 
$\wt\Delta_{\Sh^\omega(\A)}$ to $\Ker(\ve_{\Sh^\omega(\A)})$ is a 
linear map from 
$\Ker(\ve_{\Sh^\omega(\A)})$ to $\Ker(\ve_{\Sh^\omega(\A)})^{\otimes 2}$. 
Let us define $\wt\conc_{\Sh_\A}$ as the linear map from 
$\Ker(\ve_{\Sh^\omega(\A)})^{\otimes 2}$ to $\Ker(\ve_{\Sh^\omega(\A)})$
such that if $y = \sum_{k\geq 2} y_k$, with $y_k\in \oplus_{k',k''|
k' + k'' = k} \A^{\otimes k'}
\otimes \A^{\otimes k''}$, then $\wt\conc_{\Sh^\omega(\A)}(y)
= \sum_{k \geq 2} {1\over {k-1}} \conc(y_k)$. Let us also denote by 
$\wt\conc^{(2)}$ the map from $\Ker(\ve_{\Sh^\omega(\A)})^{\otimes 3}$ 
to $\Ker(\ve_{\Sh^\omega(\A)})^{\otimes 2}$ whose restriction to 
$\oplus_{k,k',k''|k + k' + k'' = n} \A^{\otimes k} 
\otimes \A^{\otimes k'} \otimes \A^{\otimes k''}$ is equal to 
${1\over {n-1}}(\conc_{\Sh^\omega(\A)}\otimes \id - \id \otimes\conc_{\Sh^{\omega^\vee}(\B)})$. 
Then we have the following homotopy formula
\begin{equation} \label{homotopy}
\id_{\Ker(\ve_{\Sh^\omega(\A)})^{\otimes 2}} = 
\wt\Delta_{\Sh^\omega(\A)} \circ  \wt\conc_{\Sh(\B)} 
+ \wt\conc^{(2)} \circ (\wt\Delta_{\Sh^\omega(\A)}  \otimes \id  
- \id \otimes \wt\Delta_{\Sh(\A)} ) . 
\end{equation}
$(\ve_{\Sh^\omega(\A)}\otimes \id\otimes \id)(Z) = (\id \otimes \ve_{\Sh^\omega(\A)}\otimes \id)(Z)
= 0$, so $Z$ belongs to $\Ker(\ve_{\Sh^\omega(\A)})^{\otimes 2}\otimes \B$. 
Let us apply identity (\ref{homotopy}) to the first two tensor factors of $Z$. 
Since $(\wt\conc_{\Sh^\omega(\A)}\otimes \id_{\B})(Z)$ belongs to 
$\Sh^\omega(\A)^{\otimes 2}\otimes \B$, this term coincides with
$(\wt\conc_{\Sh^\omega(\A)}\otimes \pr_{\B})(Z)$, therefore it
is equal to the sum of
$(n-1)\kappa(r_{\A\B})(\la_n)$ and contributions of the 
$\cR'_k(r_{\A\B})$, where $k<n$. By the construction of $\la_n$, 
this sum is zero. 
 
Let us compute now 
\begin{equation} \label{imanu'h}
\big( (\wt\Delta_{\Sh^\omega(\A)}  \otimes \id_{\Sh^\omega(\A)}  
- \id_{\Sh^\omega(\A)} \otimes \wt\Delta_{\Sh^\omega(\A)} ) \otimes \id_{\Sh(\B)} \big)(Z).
\end{equation} 
Recall that $Z$ is equal to 
$(\wt \Delta_{\Sh^\omega(\A)} \otimes \id_{\Sh(\B)})
(\cR'_n(r_{\A\B})) - 
\sum_{k = 1}^{n - 1} \cR'_k(r_{\A\B})^{(13)}
\cR'_{n-k}(r_{\A\B})^{(23)}$; since $(\wt\Delta_{\Sh^{\omega^\vee}(\B)}\otimes 
\id_{\Sh(\B)} - \id_{\Sh(\B)}\otimes
\wt\Delta_{\Sh^{\omega^\vee}(\B)} )\circ \wt\Delta_{\Sh^{\omega^\vee}(\B)} = 0$, (\ref{imanu'h}) is equal to 
$$
- \big( (\wt\Delta_{\Sh^\omega(\A)}  \otimes \id_{\Sh(\A)}  
- \id_{\Sh(\A)} \otimes \wt\Delta_{\Sh^\omega(\A)} ) \otimes 
\id_{\Sh^\omega(\B)} \big)
\big( \sum_{k = 1}^{n - 1} \cR'_k(r_{\A\B})^{(13)}
\cR'_{n-k}(r_{\A\B})^{(23)} \big) . 
$$
The induction hypothesis implies that 
$$
\big( \wt\Delta_{\Sh(^\omega\A)}  \otimes \id_{\Sh(\A)}  
\otimes \id_{\Sh(\B)} \big)
\big( \sum_{k = 1}^{n - 1} \cR'_k(r_{\A\B})^{(13)}
\cR'_{n-k}(r_{\A\B})^{(23)} \big)
$$
and 
$\big( \id_{\Sh(\A)} \otimes \wt\Delta_{\Sh^\omega(\A)}  
\otimes \id_{\Sh(\B)} \big)
\big( \sum_{k = 1}^{n - 1} \cR'_k(r_{\A\B})^{(13)}
\cR'_{n-k}(r_{\A\B})^{(23)} \big)$ are both equal to 
$\sum_{k,k',k''|k, k', k'' >0, k + k' + k''= n}
\cR'_{k}(r_{\A\B})^{(14)} \cR'_{k'}(r_{\A\B})^{(24)}
\cR'_{k''}(r_{\A\B})^{(34)}$. Therefore (\ref{imanu'h})
vanishes. The homotopy formula (\ref{homotopy}) then implies 
that $Z$ is zero. This proves the induction step. 

Equations (\ref{ids:Delta}) imply that if $t$ is a formal parameter and
$\cR^\omega(r_{\A\B}) = \sum_{n\geq 0} t^n \cR^\omega_n(r_{\A\B})$, we
have  $(S_{\Sh^\omega(\A)}\otimes \id)(\cR^\omega(r_{\A\B})) = \cR^\omega(r_{\A\B})^{-1}$; 
as $S_{\Sh^{\omega}(\B)}^{-1}$ is the antipode corresponding to the
coproduct $\Delta'_{\Sh^{\omega}(\B)}$, we have also  
and $(\id\otimes S_{\Sh^\omega(\B)}^{-1})(\cR^\omega(r_{\A\B})) 
= \cR^\omega(r_{\A\B})^{-1}$. These equalities imply 
(\ref{behavior:S}). 
\hfill \qed\medskip

\begin{remark} \label{rem:dual:R}
The sequence $(\cR^\omega_n(r_{\A\B}))_{n\geq 0}$ is uniquely 
determined by the conditions that it satisfies (\ref{ids:Delta}), 
$\cR^\omega_0(r_{\A\B}) = 1$, $\cR^\omega_1(r_{\A\B}) = r_{\A\B}$, 
and $(\pr_\A\otimes\pr_\B)(\cR^\omega(r_{\A\B})) = 0$ for $n\geq 2$. 
It follows that $\cR'_n(r_{\A\B})$ is equal to 
$\cR'_n(r_{\A\B}^{(21)})^{(21)}$. 
\end{remark}

\begin{remark} \label{rem:QT} 
If the sum $\cR^\omega(r_{\A\B}) = \sum_{i\geq 0} \cR^\omega_i(r_{\A\B})$ makes sense
(for example if we work with $\hbar$-adic completions and $r_{\A\B}$
is in $O(\hbar)$), it satisfies the $R$-matrix identities
$$
(\Delta_{\Sh^\omega(\A)}\otimes \id_{\Sh(\B)})(\cR^\omega(r_{\A\B}))  =
(\cR^\omega)^{(13)}(r_{\A\B})(\cR^\omega)^{(23)}(r_{\A\B}),  
$$
$$ 
(\id_{\Sh(\A)} \otimes \Delta_{\Sh^\omega(\B)})(\cR^\omega(r_{\A\B}))  
=(\cR^\omega)^{(13)}(r_{\A\B})(\cR^\omega)^{(12)}(r_{\A\B}).  
$$ 
\end{remark}

\begin{remark}
If $\omega$ is $\omega'$ of Section \ref{sect:explicit}, 
then the first $\cR^\omega_i(r_{\A\B})$ are 
\begin{align*}
& \cR^\omega_2(r_{\A\B})  = \sum_{i,j} (a_i a_j) \otimes (b_i b_j) + (a_i a_j)\otimes (b_j b_i)
- {1\over 2}([a_i,a_j])\otimes (b_i b_j) 
+  {1\over 2}(a_i a_j)\otimes ([b_i,b_j])   
\\ & 
= \sum_{i,j} {1\over 2}([a_i,a_j]) \otimes (b_j b_i)  + (a_i a_j) \otimes (b_i)(b_j)  
= \sum_{i,j} {1\over 2}(a_ia_j) \otimes ([b_i, b_j])  + (a_i)(a_j) \otimes (b_j b_i) 
\end{align*}
\begin{align*}
& \cR^\omega_3(r_{\A\B}) = 
\sum_{i,j,k} (a_ia_ja_k) \otimes \big( (b_ib_jb_k) + \on{all\ permutations\ in\ } i,j,k  \big)
\\ & 
+ {1\over 2}(a_ia_ja_k) \otimes \big( [b_i,b_j]b_k  + b_k[b_i,b_j]  + [b_i,b_k]b_j  
+ b_j[b_i,b_k] + [b_j,b_k]b_i  + b_i[b_j,b_k] \big)  
\\ & 
+ {1\over 2}\big( [a_i,a_j]a_k  + a_k[a_i,a_j]  + [a_i,a_k]a_j  + a_j[a_i,a_k] 
+ [a_j,a_k]a_i  + a_i[a_j,a_k] \big)  
\otimes  (b_kb_jb_i) 
\\ & 
+ (a_ia_ja_k) \otimes (L_3(b_i,b_j,b_k)) + (L_3(a_i,a_j,a_k))\otimes (b_kb_jb_i) 
\\ & 
+ {1\over 4}([a_i,a_j]a_k)\otimes \big( [b_j,b_k]b_i + b_j [b_i,b_k]\big)   
+ {1\over 4}(a_i[a_j,a_k])\otimes \big( b_k[b_i,b_j] + [b_i,b_k] b_j\big)  
\\ & 
= \sum_{i,j,k}  (a_i)(a_j)(a_k) \otimes  (b_kb_jb_i) 
+  {1\over 2}(a_i a_j)(a_k)  \otimes   (b_k[b_i,b_j]) 
+ {1\over 2}(a_i)(a_ja_k) \otimes  ([b_j,b_k]b_i) 
\\ & + (a_i a_j a_k) \otimes  ( L_3(b_i,b_j,b_k) ) , 
\end{align*}
where 
$$
L_3(x,y,z) = {1\over 6} \big( [x,[y,z]] + [[x,y],z] \big).  
$$
\end{remark}

\newpage 

\section{Construction of solutions of QYBE} \label{sect:QYBE}

In this Section, we again fix an element $\omega\in\cB(\KK[[\hbar]])$. 
Let us assume that the Lie algebras $\A$ and $\B$ are the same, 
and that that we are given an element $r_\A = \sum_{i\in I} 
a_i\otimes b_i$ of $\A\otimes\A$, 
satisfying the classical Yang-Baxter equation (CYBE), that is 
$$
[r_\A^{(12)},r_\A^{(13)}] + [r_\A^{(12)},r_\A^{(23)}] + [r_\A^{(13)},
r_\A^{(23)}] = 0. 
$$ 

Our aim in this Section is to construct a family 
$(\rho^\omega_n(r_\A))_{n\geq 1}$ of elements  
of $\A\otimes\A$, expressed in terms of Lie polynomials of degree $n$
in the $a_i$ and $b_i$, such that  if we set $\rho^\omega(r_\A) = 
\sum_{n\geq 0}
\rho^\omega_n(r_\A)$, 
and if we denote by 
$\cR^\omega_i(\rho^\omega(r_\A))[p]$ the homogeneous component of degree $p$
in $r_\A$ of $\cR^\omega_i(\rho^\omega(r_\A))$, then for each integer $n\geq 0$, 
$\rho^\omega(r_\A)$ satisfies the equation
\begin{align*}
& \sum_{p,p',p''|p + p' + p'' = n}
\sum_{i,j,k\geq 0} \cR^\omega_i(\rho^\omega(r_\A))[p]^{(12)} 
\cR^\omega_j(\rho^\omega(r_\A))[p']^{(13)}  
\cR^\omega_k(\rho^\omega(r_\A))[p'']^{(23)}  
\\ & = \sum_{p,p',p''|p + p' + p'' = n} 
\sum_{i,j,k\geq 0}\cR^\omega_k(\rho^\omega(r_\A))[p'']^{(23)} 
 \cR^\omega_j(\rho^\omega(r_\A))[p']^{(13)}  
\cR^\omega_i(\rho^\omega(r_\A))[p]^{(12)} .     
\end{align*}

In particular, if the sum $\sum_i \cR^\omega_i(\rho^\omega(r_\A))$ converges
(for example if $r_\A$ has positive $\hbar$-adic valuation), 
then $\cR^\omega(\rho^\omega(r_\A)) = \sum_{i\geq 0}\cR^\omega_i(\rho^\omega(r_\A))$ satisfies 
the quantum Yang-Baxter equation (QYBE) 
\begin{align*}
& 
(\cR^\omega(\rho^\omega(r_\A)))^{(12)} (\cR^\omega(\rho^\omega(r_\A)))^{(13)}  
(\cR^\omega(\rho^\omega(r_\A)))^{(23)}  
\\ & 
= (\cR^\omega(\rho^\omega(r_\A)))^{(23)}  (\cR^\omega(\rho^\omega(r_\A)))^{(13)}  
(\cR^\omega(\rho^\omega(r_\A)))^{(12)} .    
\end{align*}
After we reduce this problem to a Lie algebraic problem
(Section \ref{sect:reduction}), we look for universal formulas
expressing $\rho(r)$ and translate the Lie algebraic problem in 
terms of these universal formulas (Section \ref{sect:universal}). 
We then show that the latter problem can be formulated in 
cohomological terms (Section \ref{sect:coh:formulation}) and 
compute the relevant cohomologies 
(Section \ref{sect:comp:cohomologies}). We then gather
our results to show the existence and unicity of universal 
formulas for $\rho^\omega(r_\A)$ (Section \ref{sect:thm}).

\subsection{Reduction to a Lie algebraic problem}
\label{sect:reduction}

Let $\hbar$ be a formal variable and let  
$\rho$ be an element of $\hbar(\A\otimes\A)[[\hbar]]$. Then 
$\cR^\omega(\rho)$ is an element of $\Sh(\A)^{\otimes 2}[[\hbar]]$.    
Replacing the ground field $\KK$ by the ring $\KK[[\hbar]]/(\hbar^n)$ in the
definition of $\cR^\omega$, we define a map $\cR^\omega : \hbar\A^{\otimes 2}[[\hbar]] / 
\hbar^n \A^{\otimes 2}[[\hbar]] \to \Sh(\A)^{\otimes 2}[[\hbar]] / 
(\hbar^n)$. The following Proposition shows that the condition that  
$\cR^\omega(\rho)$ be a solution of the associative QYBE is equivalent to 
$\rho$ being a solution of equation (\ref{Lie:QYBE}), 
a Lie algebraic equation, 
which we call the Lie version of QYBE.

\begin{prop} \label{prop:equivalence}
Let $\pr_\A$ denote the projection of $\Sh^\omega(\A)$ on $\A$ parallel 
to $\oplus_{i\neq 1}\A^{\otimes i}$. Let $n$ be an integer and $\rho_n$ belong 
to $\hbar\A^{\otimes 2}[[\hbar]] / \hbar^n \A^{\otimes 2}[[\hbar]]$. 
Then the following statements are equivalent

{\em i)}  $\cR^\omega(\rho_n)$ satisfies QYBE in $\Sh^\omega(\A)^{\otimes 3}[[\hbar]] / (\hbar^n)$ 

{\em ii)} the equation 
\begin{align} \label{eq:rho:n}
& \pr_\A^{\otimes 3}(
(\cR^\omega(\rho_n))^{(12)} (\cR^\omega(\rho_n))^{(13)}  (\cR^\omega(\rho_n))^{(23)} ) 
\\ & \nonumber 
= \pr_\A^{\otimes 3}((\cR^\omega(\rho_n))^{(23)}  
(\cR^\omega(\rho_n))^{(13)}  (\cR^\omega(\rho_n))^{(12)} )  
\end{align}
holds in $\A^{\otimes 3}[[\hbar]] / (\hbar^n)$. 

It follows that if $\rho$ belongs to $\hbar\A^{\otimes 2}[[\hbar]]$,
the statements 

{\em iii)}  $\cR^\omega(\rho)$ satisfies QYBE in $\Sh^\omega(\A)^{\otimes 3}[[\hbar]]$ 

{\em iv)} the equation 
\begin{equation} \label{eq:rho}
\pr_\A^{\otimes 3}(
(\cR^\omega(\rho))^{(12)} (\cR^\omega(\rho))^{(13)}  (\cR^\omega(\rho))^{(23)} ) 
= \pr_\A^{\otimes 3}((\cR^\omega(\rho))^{(23)}  
(\cR^\omega(\rho))^{(13)}  (\cR^\omega(\rho))^{(12)} )  
\end{equation}
holds in $\A^{\otimes 3}[[\hbar]]$. 
\end{prop}

{\em Proof.} Let us prove the equivalence between {\em i)} and {\em ii)} by induction 
over $n$. The case $n = 0$ is trivial. Assume that the equivalence holds at order $n$
and let us show it at order $n+1$.  The direct sense is obvious. 
Let us assume now that $\rho_{n+1}$ satisfies 
(\ref{eq:rho:n}) at order $n+1$ and let us show that 
$$
Z = \cR^\omega(\rho_{n+1})^{(12)} \cR^\omega(\rho_{n+1})^{(13)}  \cR^\omega(\rho_{n+1})^{(23)}  
- \cR^\omega(\rho_{n+1})^{(23)}  \cR^\omega(\rho_{n+1})^{(13)}  \cR^\omega(\rho_{n+1})^{(12)}    
$$
(an element of $\Sh^\omega(\A)^{\otimes 3}[[\hbar]] / (\hbar^{n+1})$) is zero. 
By the induction hypothesis, 
$Z$ belongs to $\hbar^n \Sh^\omega(\A)^{\otimes 3}[[\hbar]] / 
\hbar^{n+1}\Sh^\omega(\A)^{\otimes n+1}[[\hbar]]$. 
Set 
$$
L = \cR^\omega(\rho_{n+1})^{(12)} \cR^\omega(\rho_{n+1})^{(13)} , 
$$ 
$$
R = \cR^\omega(\rho_{n+1})^{(23)}   
(\cR^\omega(\rho_{n+1})^{(13)}  \cR^\omega(\rho_{n+1})^{(12)})  
(\cR^\omega(\rho_{n+1})^{(23)} )^{-1} . 
$$
Since   $L - R = Z(\cR^\omega(\rho_{n+1})^{(23)})^{-1}$ and 
$\cR^\omega(\rho_{n+1})$ belongs to 
$$
1 + \hbar \Sh^\omega(\A)^{\otimes 2}[[\hbar]] / 
\hbar^{n+1}\Sh^\omega(\A)^{\otimes 2}[[\hbar]],
$$
we get $L - R = Z$. On the other hand, 
both $L$ and $R$ satisfy  
\begin{equation} \label{51}
(\Delta_{\Sh^\omega(\A)}\otimes \id^{\otimes 2})(L) = L^{(134)} L^{(234)} \quad \on{and}\quad   
(\Delta_{\Sh^\omega(\A)}\otimes \id^{\otimes 2})(R) = R^{(134)} R^{(234)}.   
\end{equation}
Let us write $L = \sum_{i\geq 0}\hbar^i L_i$, $R = \sum_{i\geq 0}\hbar^i R_i$, 
with $L_i,R_i\in \Sh^\omega(\A)^{\otimes 3}$. Then the induction hypothesis implies
that $L_i = R_i$ for $i\leq n$ and  (\ref{51}) implies that 
$$
(\wt\Delta_{\Sh^\omega(\A)}\otimes \id^{\otimes 2})(L_{n+1}) = \sum_{i = 1}^n L_i^{(134)} L_{n+1-i}^{(234)} , 
\quad
(\wt\Delta^{\Sh^\omega(\A)}\otimes \id^{\otimes 2})(R_{n+1}) = \sum_{i = 1}^n R_i^{(134)} R_{n+1-i}^{(234)}, 
$$
where $\wt\Delta_{\Sh^\omega(\A)}(x) = \Delta(x) - x\otimes 1 - 1\otimes x$.  Therefore 
$(\wt\Delta_{\Sh^\omega(\A)}\otimes \id^{\otimes 2})(L_{n+1}) = 
(\wt\Delta_{\Sh^\omega(\A)}\otimes \id^{\otimes 2})(R_{n+1})$, so 
$L_{n+1} - R_{n+1}$ belongs to $\A\otimes \Sh^\omega(\A)^{\otimes 2}$. Let us define 
$Z_{n+1}$ as $\hbar^{-n-1} Z$ mod $\hbar$. Then $Z_{n+1}$, then $Z_{n+1}
 = L_{n+1} - R_{n+1}$, therefore $Z_{n+1}$ belongs to 
$\A\otimes \Sh^\omega(\A^{\otimes 2})$.  

Applying the same reasoning to the comparison of 
$$
\cR^\omega(\rho_{n+1})^{(13)}  \cR^\omega(\rho_{n+1})^{(23)} 
\ \on{and} \   
(\cR^\omega(\rho_{n+1})^{(12)})^{-1} \big( \cR^\omega(\rho_{n+1})^{(23)}  
\cR^\omega(\rho_{n+1})^{(13)}  \big) \cR^\omega(\rho_{n+1})^{(12)} ,    
$$
resp., of 
$$
\cR^\omega(\rho_{n+1})^{(12)} \cR^\omega(\rho_{n+1})^{(23)}  
\quad \on{and} \quad 
\sum_{\la}
B_\la^{(3)} \cR^\omega(\rho_{n+1})^{(23)}  \cR^\omega(\rho_{n+1})^{(13)}  \cR(\rho_{n+1})^{(12)}  
A_\la^{(1)} , 
$$
where $\sum_\la A_\la\otimes B_\la$ is the inverse of $\cR^\omega(\rho_{n+1})^{-1}$
in $\Sh^\omega(\A)\otimes\Sh^\omega(\A)^{opp}[[\hbar]] / (\hbar^{n+1})$, where $\Sh^\omega(\A)^{opp}$ is the opposite
algebra of $\Sh^\omega(\A)$,  we find that 
$Z_{n+1}$ belongs to $\Sh^\omega(\A^{\otimes 2}) \otimes \A$, resp., to $\Sh^\omega(\A)\otimes \A\otimes 
\Sh^\omega(\A)$. Therefore, $Z_{n+1}$ belongs to $\A^{\otimes 3}$, so it may be 
identified with its image by $\pr_\A^{\otimes 3}$. Since by assumption  
$\pr_\A^{\otimes 3}(Z)$ is zero, we have  $\pr_\A^{\otimes 3}(Z_{n+1}) = 0$, 
therefore $Z_{n+1}$ is equal to zero and $Z$ is zero as well. 
This proves the induction step. 

The equivalence between {\em iii)} and {\em iv)} follows easily from that of 
{\em i)} and {\em ii)}. 
\hfill \qed\medskip 

We will call equation (\ref{eq:rho}) the {\em Lie version of QYBE} (Lie
QYBE for short). This equation may be expressed as follows.

For $a_1,\ldots,a_n$ elements of a Lie  algebra $\G$, 
let us denote by $L^\G_n(a_1,\ldots,a_n)$ the projection of the first summand of 
$\Sh^\omega(\G) = \G \oplus (\oplus_{i\neq 1}\G^{\otimes i})$ of 
the product $(a_1)\cdots(a_n)$; so $L^\G_n(a_1,\ldots,a_n)$ belongs to $\G$.  
When $n = 0$, we set  $L^\G_n(x_1,\ldots,x_n) = 0$. 
There is a unique element of $FL_n$, which we denote by 
$L_n$, such that 
$L^\G_n(a_1,\ldots,a_n)$ is obtained from $L_n$
by substituting $a_i$ to the $i$th generator of $FL_n$, $i = 1,\ldots,n$. 
$L_n$ depends on $\omega$. 
If $\un = (n_1,\ldots,n_k)$, we set 
$\la_\un = \la^{\prime (1)}_\un \otimes\cdots\otimes 
\la^{\prime (k)}_\un \otimes \la''_\un$. 

Then we have 
\begin{align*}
& \pr_\A^{\otimes 3}\big( 
\cR^\omega(\rho)^{(12)}\cR^\omega(\rho)^{(13)}\cR^\omega(\rho)^{(23)}
\big) 
= \sum_{\xi',\xi'',\eta\geq 0} \sum_{\al_1,\ldots,\al_\eta\geq 0}
\sum_{n_1,\ldots,n_{\xi' + \xi''}\geq 0, m_{11},\ldots,m_{\eta\al_\eta}\geq 0}
\\ & 
\Big( 
L_{\al_1 + \cdots + \al_\eta + \xi'}(\la^{\prime (1)}_{\um_1},
\ldots,\la^{\prime (\al_\eta)}_{\um_\eta},\la^{\prime (1)}_\un,
\ldots,\la^{\prime (\xi')}_\un)
\otimes
\\ & 
\otimes 
L_{\eta+\xi''}(\la''_{\um_1},\ldots,\la''_{\um_\eta},
\la^{\prime(\xi'+1)}_\un,\ldots,\la^{\prime(\xi' + \xi'')}_\un)
\otimes \la''_\un\Big) (\rho) , 
\end{align*}
where we set $\un = (n_i)_{i = 1,\ldots,\xi'+\xi''}$, 
$\um_\beta = (m_{\beta i})_{i = 1,\ldots,\al_\beta}$, for
$\beta = 1,\ldots,\eta$. 
In the same way, 
\begin{align*}
& \pr_\A^{\otimes 3}\big( 
\cR^\omega(\rho)^{(23)}\cR^\omega(\rho)^{(13)}\cR^\omega(\rho)^{(12)}
\big) 
= \sum_{\xi',\xi'',\eta\geq 0} \sum_{\al_1,\ldots,\al_\eta\geq 0}
\sum_{n_1,\ldots,n_{\xi' + \xi''}\geq 0, m_{11},\ldots,m_{\eta\al_\eta}\geq 0}
\\ & 
\Big( 
L_{\al_1 + \cdots + \al_\eta + \xi''}(
\la^{\prime(\xi'+1)}_\un,\ldots,\la^{\prime(\xi' + \xi'')}_\un,
\la^{\prime (1)}_{\um_1},
\ldots,\la^{\prime (\al_\eta)}_{\um_\eta})
\otimes
\\ & 
\otimes 
L_{\eta+\xi'}( \la^{\prime (1)}_\un, \ldots,\la^{\prime (\xi')}_\un,
\la''_{\um_1},\ldots,\la''_{\um_\eta})
\otimes \la''_\un\Big) (\rho) .  
\end{align*}
So we have 
\begin{cor} \label{cor:malta}
The Lie QYB equation (\ref{eq:rho}) is equivalent to the  
following equation
\begin{align} \label{Lie:QYBE}
& \nonumber
\sum_{\xi',\xi'',\eta\geq 0} \sum_{\al_1,\ldots,\al_\eta\geq 0}
\sum_{n_1,\ldots,n_{\xi' + \xi''}\geq 0, m_{11},\ldots,m_{\eta\al_\eta}\geq 0}
\\ & \nonumber
\Big( 
L_{\al_1 + \cdots + \al_\eta + \xi'}(\la^{\prime (1)}_{\um_1},
\ldots,\la^{\prime (\al_\eta)}_{\um_\eta},\la^{\prime (1)}_\un,
\ldots,\la^{\prime (\xi')}_\un)
\otimes
\\ & \nonumber
\otimes 
L_{\eta+\xi''}(\la''_{\um_1},\ldots,\la''_{\um_\eta},
\la^{\prime(\xi'+1)}_\un,\ldots,\la^{\prime(\xi' + \xi'')}_\un)
\otimes \la''_\un\Big) (\rho) 
\\ & \nonumber
= \sum_{\xi',\xi'',\eta\geq 0} \sum_{\al_1,\ldots,\al_\eta\geq 0}
\sum_{n_1,\ldots,n_{\xi' + \xi''}\geq 0, m_{11},\ldots,m_{\eta\al_\eta}\geq 0}
\\ & \nonumber
\Big( 
L_{\al_1 + \cdots + \al_\eta + \xi''}(
\la^{\prime(\xi'+1)}_\un,\ldots,\la^{\prime(\xi' + \xi'')}_\un,
\la^{\prime (1)}_{\um_1},
\ldots,\la^{\prime (\al_\eta)}_{\um_\eta})
\otimes
\\ & 
\otimes 
L_{\eta+\xi'}( \la^{\prime (1)}_\un, \ldots,\la^{\prime (\xi')}_\un,
\la''_{\um_1},\ldots,\la''_{\um_\eta})
\otimes \la''_\un\Big) (\rho) ,   
\end{align}
where $\rho$ is an element of $\hbar (\A\otimes\A)[[\hbar]]$. 
\end{cor}

\subsection{Universal formulation of the problem}
\label{sect:universal}

Let us explain more precisely the nature of the function 
$r_\A\mapsto \rho_n(r_\A)$.  

If $n$ is an integer, let us denote by $Free_n$ the part of the 
free Lie algebra with $n$ generators, homogeneous of degree one in 
each generator.  The action of $\SG_n$ by permutation of the generators
induces a $\SG_n$-module structure on $Free_n$.  

If $\Gamma$ is a group acting on a vector space $M$, we denote
by $M_\Gamma$ the space of coinvariants of $M$ with respect to $\Gamma$; 
it is defined as $M_\Gamma = M / \Span\{\gamma m - m, \gamma\in \Gamma, m\in M\}$.

\begin{lemma} \label{RH}
Let us set $F_n = (Free_n \otimes Free_n)_{\SG_n}$. Let $(\A,r_\A)$
be the pair of a Lie algebra and an element $r_\A = \sum_{i\in I}
a_i\otimes b_i\in \A\otimes\A$ satisfying of CYBE.   
There is a unique map $\kappa^{(ab)}_{r_\A} : F_n\to\A\otimes \A$, such that 
if $P(x_1,\ldots,x_n)$ and $Q(x_1,\ldots,x_n)$ belong to 
$Free_n$, then the image by $\kappa^{(ab)}_{r_\A}$ of the class of 
$P\otimes Q$ is
$$
\sum_{i_1,\ldots,i_n\in I} P(a_{i_1},\ldots,a_{i_n})
\otimes Q(b_{i_1},\ldots,b_{i_n}). 
$$
\end{lemma}

{\em Proof.} This follows from the fact that for any $\sigma$ in 
$\SG_n$, we have 
$$
\sum_{i_1,\ldots,i_n\in I} P(a_{i_1},\ldots,a_{i_n})
\otimes Q(b_{i_1},\ldots,b_{i_n})
= 
\sum_{i_1,\ldots,i_n\in I} P(a_{i_{\sigma(1)}},\ldots,a_{i_{\sigma(n)}})
\otimes Q(b_{i_{\sigma(1)}},\ldots,b_{i_{\sigma(n)}}). 
$$
\hfill \qed\medskip 

We will introduce equations for a system of elements 
$(\varrho_n)_{n\geq 1}$, where $\varrho_n$ belongs to $F_n$, 
which we will call the universal Lie QYB equations. 

The universal Lie QYB equations have the following property: 
if $\A$ is a Lie algebra and $r_\A\in \A\otimes \A$ is a solution of CYBE, if 
$(\varrho_n)_{n\geq 1}\in \prod_{n\geq 1}F_n$ is a solution of the 
universal Lie QYB equations, and if  
we set $\rho_n(r_\A) = \kappa^{(ab)}_{r_\A}(\varrho_n)$, then the family 
$(\rho_n(r_\A))_{n\geq 0}$ satisfies the Lie QYB equations in 
$\A^{\otimes 3}$. 

After we establish this fact, we will show that the  
universal Lie QYB equations have a unique solution such that 
$\varrho_1 = x_1\otimes y_1$. 

\subsubsection{Insertion map}

Let $p,q,r$ be integers. Denote by $FL_{p,q,r}$ the space 
$$
FL_{p,q,r} = \big( 
FL_{q+r}\otimes FL_{p+r} \otimes FL_{p+q}
\big)_{\SG_p\times\SG_q\times\SG_r} ; 
$$
here $FL_{q+r}$ (resp., $FL_{p+r}$, $FL_{p+q}$) 
is generated by the variables 
$v_1,\ldots,v_q$, $w_1,\ldots,w_r$
(resp., by $u_1,\ldots,u_p$, $w'_1,\ldots,w'_r$ 
and by $u'_1,\ldots,u'_p$, $v'_1,\ldots,v'_q$), and
$\SG_p$, $\SG_q$ and $\SG_r$
acts by simultaneously permuting the variables
$x_i$ and $x'_i$ ($x = u,v,w$).

Then there is a unique linear map 
$$
\ins : FL_{p,q,r} \otimes \big( \prod_{n\geq 1} F_n \big) \to 
\prod_{p',q',r'\geq 0} FL_{p',q',r'}, 
$$
such that the $(p',q',r')$ component of 
$\ins\big((P\otimes Q\otimes R)\otimes (\sum_{n\geq 1}\sigma_n) \big)$
is 
\begin{align*}
& \sum_{(n_i),(m_i),(s_i)|\sum_{i = 1}^p n_i = p', 
\sum_{i = 1}^p m_i = q', \sum_{i = 1}^p s_i = r'} 
P(\wt\sigma'_{m_1},\ldots,\wt\sigma'_{m_q}, 
\wt\sigma'_{s_1},\ldots,\wt\sigma'_{s_r}) 
\\ & 
\otimes Q(\wt\sigma'_{n_1},\ldots,\wt\sigma'_{n_p},
\wt\sigma''_{s_1},\ldots,\wt\sigma''_{s_r}) \otimes 
R(\wt\sigma''_{n_1},\ldots,\wt\sigma''_{n_p},\wt\sigma''_{m_1},\ldots,\wt\sigma''_{m_q}) , 
\end{align*}
where we set $\sigma_n = \sigma'_n\otimes\sigma''_n$ and  
$\wt\sigma'_{m_i} = \sigma'_{m_i}(v_{m_1+\cdots + m_{i-1} + 1}
,\ldots,v_{m_1 + \cdots + m_i})$, 
$\wt\sigma''_{m_i} = \sigma''_{m_i}(v'_{m_1+\cdots + m_{i-1} + 1}
,\ldots,v'_{m_1 + \cdots + m_i})$, 
$\wt\sigma'_{n_i} = \sigma'_{n_i}(u_{n_1+\cdots + n_{i-1} + 1}
,\ldots,u_{n_1 + \cdots + n_i})$, 
$\wt\sigma''_{n_i} = \sigma''_{n_i}(u'_{n_1+\cdots + n_{i-1} + 1}
,\ldots,u'_{n_1 + \cdots + n_i})$, and 
$\wt\sigma'_{s_i} = \sigma'_{s_i}(w_{s_1+\cdots + s_{i-1} + 1}
,\ldots,w_{s_1 + \cdots + s_i})$, 
$\wt\sigma''_{s_i} = \sigma''_{s_i}(w'_{s_1+\cdots + s_{i-1} + 1}
,\ldots,w'_{s_1 + \cdots + w_i})$.

\subsubsection{Universal version of the Lie QYBE identity} 
\label{sect:univ:Lie}\label{vacances}

When $N$ is an integer $\geq 0$, let us denote by
$F_N^{(aab)}$ and $F_N^{(abb)}$ the direct sums 
$$
F_N^{(aab)} = \oplus_{p,q\geq 1|p+q = N} 
(FL_p \otimes FL_q \otimes 
 FL_N)_{\SG_p \times\SG_q}, 
$$
$$
F_N^{(abb)} = \oplus_{p,q\geq 1|p+q = N} 
(FL_N\otimes FL_p \otimes FL_q)_{\SG_p \times\SG_q}, 
$$
where the $Free_N$ (resp., $Free_p$ and $Free_q$) is endowed with 
the action of $\SG_p\times\SG_q$ provided by the 
inclusion map $\SG_p\times\SG_q\to\SG_N$ (resp., the projection map 
of $\SG_p\times\SG_q$ on $\SG_p$ and $\SG_q$).   So 
$F_N^{(aab)} = \oplus_{p,q|p+q = N} FL_{p,q,0}$ and 
$F_N^{(abb)} = \oplus_{p,q|p+q = N} FL_{0,p,q}$. 

In the same way as Lemma \ref{RH}, we have
\begin{lemma}
There are unique linear maps $\kappa_r^{(aab)} : F_N^{(aab)} \to \A^{\otimes 3}$ 
and $\kappa_r^{(abb)} : F_N^{(abb)}\to \A^{\otimes 3}$, such that for $P\in Free_N$
and $P',P''$ in $Free_p$ and $Free_q$, the image by $\kappa_N^{(abb)}$ of the class of 
$P\otimes P'\otimes P''$ is 
$$
\sum_{i_1,\ldots,i_N\in I} P(a_{i_1},\ldots, a_{i_N})\otimes
P'(b_{i_1},\ldots,b_{i_p})\otimes P''(b_{i_{p+1}},\ldots,b_{i_N})
$$
and the image by $\kappa_N^{(aab)}$ of the class of 
$P'\otimes P''\otimes P$ is 
$$
\sum_{i_1,\ldots,i_N\in I} P'(a_{i_1},\ldots,a_{i_p})
\otimes P''(a_{i_{p+1}},\ldots,a_{i_N}) 
\otimes P(b_{i_1},\ldots, b_{i_N}). 
$$
\end{lemma}

\begin{prop} \label{prop:key}
If $p,q,r$ are integers, there is a unique map 
$\mu_{Lie}^{p,q,r} : FL_{p,q,r} \to F^{(aab)}_{p+q+r} \oplus
F^{(abb)}_{p+q+r}$, such that if  $(\A,r_\A)$ is the pair of a Lie 
algebra and a solution $r_\A = \sum_{i\in I} a_i\otimes b_i
\in\A\otimes \A$ of CYBE, then 
\begin{align*}
& (\kappa_{p+q+r}^{(aab)} \oplus \kappa_{p+q+r}^{(abb)})
(\mu_{Lie}^{p,q,r}(P\otimes Q\otimes R))
= \sum_{i_1,\ldots,k_r\in I}
P(a_{j_1},\ldots,a_{j_q},a_{k_1},\ldots,a_{k_r})
\\ & \otimes 
Q(a_{i_1},\ldots,
a_{i_p},b_{k_1},\ldots,b_{k_r})\otimes R(b_{i_1},\ldots,b_{i_p},
b_{j_1},\ldots,b_{j_q})
\end{align*}
We define $\mu_{Lie}$ as the direct sum 
$\oplus_{p,q,r} \mu_{Lie}^{p,q,r}$.  
\end{prop}

{\em Proof.} See Appendix \ref{proof:prop:key}. \hfill \qed\medskip

Let $(\varrho_n)_{n\geq 1}$ be an element of $\prod_{n\geq 1} F_n$. 
We say that $(\varrho_n)_{n\geq 1}$ is a solution of the {\em universal 
Lie QYB equations} if for any integer $N$, the equality 
\begin{align} \label{univ:Lie:QYBE}
& \nonumber
\sum_{\xi',\xi'',\eta\geq 0} \sum_{\al_1,\ldots,\al_\eta\geq 0}
\sum_{n_1,\ldots,n_{\xi' + \xi''}\geq 0, m_{11},\ldots,m_{\eta\al_\eta}\geq 0}
\\ & \nonumber
\mu_{Lie} \circ \ins \Big(
\big( 
L_{\al_1 + \cdots + \al_\eta + \xi'}(\la^{\prime (1)}_{\um_1},
\ldots,\la^{\prime (\al_\eta)}_{\um_\eta},\la^{\prime (1)}_\un,
\ldots,\la^{\prime (\xi')}_\un)
\otimes
\\ & \nonumber
\otimes 
L_{\eta+\xi''}(\la''_{\um_1},\ldots,\la''_{\um_\eta},
\la^{\prime(\xi'+1)}_\un,\ldots,\la^{\prime(\xi' + \xi'')}_\un)
\otimes \la''_\un\big) \otimes (\sum_n \varrho_n) \Big)_N
\\ & \nonumber
= \sum_{\xi',\xi'',\eta\geq 0} \sum_{\al_1,\ldots,\al_\eta\geq 0}
\sum_{n_1,\ldots,n_{\xi' + \xi''}\geq 0, m_{11},\ldots,m_{\eta\al_\eta}\geq 0}
\\ & \nonumber
\mu_{Lie} \circ \ins \Big(
\big( 
L_{\al_1 + \cdots + \al_\eta + \xi''}(
\la^{\prime(\xi'+1)}_\un,\ldots,\la^{\prime(\xi' + \xi'')}_\un,
\la^{\prime (1)}_{\um_1},
\ldots,\la^{\prime (\al_\eta)}_{\um_\eta})
\otimes
\\ & 
\otimes 
L_{\eta+\xi'}( \la^{\prime (1)}_\un, \ldots,\la^{\prime (\xi')}_\un,
\la''_{\um_1},\ldots,\la''_{\um_\eta})
\otimes \la''_\un\big) \otimes (\sum_n \varrho_n) \Big)_N,   
\end{align}
holds, where the index $N$ means the homogeneous component in 
$F_N^{(aab)}\oplus F_N^{(abb)}$. We call the above equation the universal Lie QYB equation 
of degree $N$. We supplement this equation by the 
condition $\varrho_1  = x\otimes x$, where $x$ is the canonical 
generator of $Free_1$. 

Then we have 
\begin{lemma} \label{lemma:3.3}
Assume that  $(\varrho_n)_{n\geq 1}$ is a solution of the universal 
Lie QYB equations. Let $\A$ be any Lie algebra and let $r_\A\in \A\otimes \A$
be a  solution of CYBE. Set $\rho_n  = \kappa^{(ab)}_{r_\A}(\varrho_n)$.  
Then $\sum_{n\geq 1} \hbar^n \rho_n$ is a solution of the Lie QYB
equation (\ref{Lie:QYBE}) in $\A^{\otimes 3}[[\hbar]]$.
\end{lemma} 

{\em Proof.} This follows at once from Proposition 
\ref{prop:key}. \hfill \qed\medskip 

\subsection{Cohomological formulation}
\label{sect:coh:formulation}

Our aim is to solve equations (\ref{univ:Lie:QYBE}): 
we will show that these equations have a unique solution. 
For this, we will formulate equations (\ref{univ:Lie:QYBE})
in cohomological terms.  

For any integer $n$, there is a unique linear map 
$\delta_{3}^{(F)} : F_n \to F^{(aab)}_{n+1} + F^{(abb)}_{n+1}$, 
such that for any $P$ and $Q$ in $FL_n$,  
\begin{align*}
&
\delta_3^{(F)}(P\otimes Q) 
=  [w_1,P(v_1,\ldots,v_n)] \otimes w'_1\otimes Q(v'_1,\ldots,v'_n) 
\\ & + 
\mu_{Lie}^{1,0,n}(w_1\otimes [w'_1,P(u_1,\ldots,u_n)]\otimes Q(u'_1,\ldots,u'_n))
+ v_1\otimes P(u_1,\ldots,u_n)\otimes [v'_1,Q(u'_1,\ldots,u'_n)]
\\ & 
+  [P(w_1,\ldots,w_n),v_1] \otimes v'_1\otimes Q(w'_1,\ldots,w'_n) 
+  \mu_{Lie}^{1,0,n}(P(w_1,\ldots,w_n)\otimes [Q(w'_1,\ldots,w'_n),
u_1]\otimes u'_1
\\ & 
+ P(v_1,\ldots,v_n)\otimes u_1\otimes [Q(v'_1,\ldots,v'_n),u'_1] . 
\end{align*}

\begin{lemma}
Let $\A$ be any Lie algebra and $r_\A\in \A\otimes\A$ be a solution
of CYBE. Let us define $\delta_{3,r_\A} : \A^{\otimes 2} \to 
\A^{\otimes 3}$ by 
$\delta_{3,r_\A}(x) = [r_\A^{12},x^{13}] 
+ [r_\A^{12},x^{23}] +  [r_\A^{13},x^{23}] 
+ [x^{12},r_\A^{13}]  + [x^{12},r_\A^{23}] +  [x^{13},r_\A^{23}]$. Then  
the diagram 
$$
\begin{array}{ccc}
F_n & \stackrel{\delta_3^{(F)}}{\rightarrow} & F_{n+1}^{(aab)} \oplus F_{n+1}^{(abb)}
\\  
\downarrow{\kappa_{r_\A}} & \; & 
\downarrow{\kappa_{r_\A}^{(aab)} \oplus \kappa_r^{(abb)}}
\\
\A^{\otimes 2} & \stackrel{\delta_{3,r_\A}}{\rightarrow} & \A^{\otimes 3}
\end{array}
$$
is commutative.  
\end{lemma}

{\em Proof.} This follows from Proposition \ref{prop:key}. 
\hfill \qed\medskip

Then there is a unique map $\Phi_N : \prod_{i = 2}^{N-2} F_i
\to F_N^{(aab)} \oplus F_N^{(abb)}$, which is polynomial 
of degree $3$, such that  the left side of equation 
(\ref{univ:Lie:QYBE}) is equal to 
$\delta_3^{(F)}(\varrho_{N-1}) + \Phi_N(\varrho_2,\ldots,\varrho_{N-2})$ 
(recall that $\varrho_1$ is the canonical generator of $F_1$). 
Equation 
(\ref{univ:Lie:QYBE}) may then be rewritten as follows
\begin{equation} \label{induction:Lie:QYBE}
\delta_3^{(F)}(\varrho_{N-1}) + \Phi_N(\varrho_2,\ldots,\varrho_{N-2}) = 0.  
\end{equation}

If $\A$ is a Lie algebra and $r_\A\in\A\otimes\A$ is a solution of CYBE, 
let us define $\delta_{4,r_\A}$ as the linear map from $\A^{\otimes 3}$
to $\A^{\otimes 4}$ such that for any $x\in\A^{\otimes 3}$, 
\begin{align*}
\delta_{4,r_\A}(x) = & [r_\A^{12} + r_\A^{13} + r_\A^{14},x^{234}] 
- [-r_\A^{12} + r_\A^{23} + r_\A^{24},x^{134}] 
\\ & + [-r_\A^{13} - r_\A^{23} + r_\A^{34},x^{124}]  
- [-r_\A^{14} - r_\A^{24} - r_\A^{34},x^{123}]  . 
\end{align*}

In the same way as $\delta_3^{(F)}$ is universal version of 
$\delta_{3,r_\A}$, we associate to $\delta_{4,r_\A}$ its universal version 
$\delta_4^{(F)}$. For this, let us first define the spaces
$F^{(x_1,\ldots,x_n)}_N$. If $n$ is an integer $\geq 1$ and 
$(x_1,\ldots,x_n)$ is a sequence of $\{a,b\}^n$, we define 
$F^{(x_1,\ldots,x_n)}_N$ as follows. Let us set 
$K = \{k | x_k = a\}$ and $L = \{l|x_l = b\}$. 
$K$ and $L$ therefore form a partition of $\{1,\ldots,n\}$. 
If $(p_{kl})_{(k,l)\in K\times L}$ belongs to $\NN^{K\times L}$, 
let us set  $f((p_{kl}),i) = \sum_{l\in L} p_{il}$ if $i\in K$ 
and $f((p_{kl}),i) = \sum_{k\in K} p_{ki}$ if $i\in L$. We then set 
\begin{equation} \label{def:lie:space}
F^{(x_1,\ldots,x_n)}_{(p_{kl})_{(k,l)\in K\times L}} = 
\Big( \bigotimes_{i = 1}^n Free_{f((p_{kl}),i)}\Big)_{\prod_{k\in K,l\in L} \SG_{p_{kl}}} , 
\end{equation}
where each $\SG_{p_{kl}}$ acts by simultaneously permuting the 
generators of index $\sum_{j =1}^{l-1}p_{kj} + 1, \ldots, \sum_{j = 1}^{l} p_{kj}$ 
of the $k$th tensor factor, and the generators of 
index $\sum_{i = 1}^{k-1} p_{il} + 1,\ldots, \sum_{i = 1}^k p_{il}$ of the 
$l$th tensor factor. 

\begin{equation} \label{big:lie:space}
F^{(x_1,\ldots,x_n)} = \bigoplus_{(p_{kl})\in \NN^{K\times L}| 
p_{kl} = 0\on{\ if\ }k>l}
\Big( \bigotimes_{i = 1}^n Free_{f((p_{kl}),i)}\Big)_{\prod_{k\in K,l\in L} \SG_{p_{kl}}}
\end{equation}
and 
$$
F^{(x_1,\ldots,x_n)}_N = 
\bigoplus_{ (p_{kl})\in \NN^{K\times L}| p_{kl} = 0\on{\ if\ }k>l
\on{\ and\ }\sum_{k\in K,l\in L} p_{kl} = N }
\Big( \bigotimes_{i = 1}^n Free_{f((p_{kl}),i)}\Big)_{\prod_{k\in K,l\in L} \SG_{p_{kl}}}. 
$$
For example, 
$$
F_N^{(aaab)} = \bigoplus_{p,p',p''|p+p'+p'' = N}
\big( Free_{p}\otimes Free_{p'}\otimes Free_{p''} 
\otimes  Free_N 
\big)_{\SG_p \times \SG_{p'} \times \SG_{p''}} 
$$
and 
$$
F_N^{(aabb)} = \bigoplus_{p,p',q,q'|p+p' + q + q' = N}
\big( Free_{p+p'} \otimes Free_{q+q'}\otimes Free_{p+q}\otimes Free_{p'+q'} 
\big)_{\SG_p \times \SG_{p'} \times \SG_{q}\times\SG_{q'}} , 
$$
where in the last equality, $\SG_p$ acts by permutation of the 
$p$ first generators of $Free_{p+p'}$ and $Free_{p+q}$, 
$\SG_{p'}$ acts by permutation of the $p'$ last (resp., first)
generators of $Free_{p+p'}$ (resp., $Free_{p'+q'}$), $\SG_q$ acts
by permutation of the $q$ first (resp., last) generators of $Free_{q+q'}$
(resp., $Free_{p+q}$), and $\SG_{q'}$ acts by permutation of the 
$q'$ last generators of $Free_{q+q'}$ and $Free_{p'+q'}$. 

As before, we associate to any Lie algebra $\A$ and any element $r_\A\in \A\otimes\A$, the 
map 
$$ 
\kappa_{r_\A}^{(x_1,\ldots,x_n)} : F^{(x_1\ldots x_n)}\to \A^{\otimes n}
$$  
such that if $(p_{kl})\in \NN^{K\times L}$ and we set 
$N = \sum_{k\in K,l\in L} p_{kl}$, and if $P_\al\in Free_{f((p_{kl}),\al)}$, then 
$$
\kappa_{r_\A}^{(x_1,\ldots,x_n)}(\otimes_{i = 1}^n P_i) = 
\sum_{\al(1),\ldots,\al(N)\in I}\bigotimes_{i = 1}^n A_i,  
$$
where $A_i =  P_i(a_{\al(p_1+\cdots + p_{i-1})+1},\ldots,a_{\al(p_1+\cdots + p_{i})})$
if $i\in K$ and 
$$
A_i =  P_i(b_{\beta(q_1+\cdots + q_{i-1})+1},\ldots,b_{\beta(q_1+\cdots + q_{i})})
$$
if $i\in L$, and we set $p_k = \sum_{l\in L}p_{kl}$ if $k\in K$,  
$q_l = \sum_{k\in K}p_{kl}$ if $l\in L$, and 
$\beta(q_1 + \cdots + q_{i-1} + p_{1i} + \cdots + p_{k-1,j} + s) = \al(p_1 + 
\cdots + p_{i-1} + p_{i1} + \cdots + p_{i,j-1} + s)$ 
if $1\leq s \leq p_{kj}$. 

\begin{prop} \label{prop:delta:4}
There exists a map 
$$
\delta_4^{(F)} : \bigoplus_{x\in \{a,b\}} F_n^{(axb)} 
\to \bigoplus_{x,y\in \{a,b\}} F_{n+1}^{(axyb)}  
$$
such that for any pair $(\A,r_\A)$ of a Lie algebra and a 
solution of CYBE, the diagram 
$$
\begin{array}{ccc}
\bigoplus_{x\in \{a,b\}} F_n^{(axb)}  & \stackrel{\delta_4^{(F)}}{\rightarrow} & 
\bigoplus_{x,y\in \{a,b\}} F_{n+1}^{(axyb)} 
\\  
\downarrow{\oplus_{x\in \{a,b\}} \kappa_{r_\A}^{(axb)}} & \; & 
\downarrow{\oplus_{x,y\in \{a,b\}}\kappa_{r_\A}^{(axyb)}}
\\
\A^{\otimes 3} & \stackrel{\delta_{4,r_\A}}{\rightarrow} & \A^{\otimes 4}
\end{array}
$$
is commutative.  
\end{prop}

{\em Proof.} See Appendix \ref{proof:prop:delta:4}.
\hfill \qed\medskip

\begin{prop} \label{prop:delta:diff}
1) If $\A$ is any Lie algebra and $r_\A\in \A\otimes\A$ is any solution of CYBE, 
we have $\delta_{4,r_\A}\circ \delta_{3,r_\A} = 0$.  

2) We have also $\delta_4^{(F)} \circ \delta_3^{(F)} = 0$.  
\end{prop}

{\em Proof.} 1) is a direct computation. 
The proof of 2) is in Appendix \ref{proof:prop:delta:4}. 
\hfill \qed\medskip 

\begin{remark}
When $r_\A + r_\A^{(21)}$ is $\A$-invariant, $\delta_{3,r_\A}$ 
and $\delta_{4,r_\A}$ are differentials of the 
Lie coalgebra cohomology complex of $\A$, endowed with the 
Lie coalgebra structure given by $\delta(x) = [r_\A,x\otimes 1 + 1\otimes x]$.
This explains the relation  $\delta_{4,r_\A}\circ \delta_{3,r_\A} = 0$.  
On the other hand, $\delta_3^{(F)}$ and $\delta_4^{(F)}$ are also 
differentials of a complex whose degree $n$ part is 
$$
\oplus_{x_1,\ldots,x_{n-2}\in \{a,b\}} F^{(a x_1\ldots x_{n-2}b)} ,  
$$
which may be viewed as a universal version of the Lie coalgebra
cohomology complex of $\A$.  
\hfill \qed\medskip  
\end{remark}

\begin{thm} \label{training}
Let $\A$ be a Lie algebra and $\rho_1,\ldots,\rho_N$ be elements of  
$\A^{\otimes 2}$. Set $r = \rho_1$ and assume that  
$\rho_1,\ldots,\rho_{N-1}$ satisfy the equations 
\begin{equation} \label{rached}
\delta_{3,r}(\rho_{M-1}) + \Phi_M(\rho_2,\ldots,\rho_{M-2}) = 0 
\end{equation}
for $M = 1,\ldots,N$. (In particular, $r$ is a solution of CYBE.)
Then 
$$
\delta_{4,r}(\Phi_{N+1}(\rho_2,\ldots,\rho_{N-1})) = 0.
$$ 
\end{thm}

{\em Proof.} 
The proof of this Theorem relies on the following Proposition, 
which will be proved in Appendix \ref{proof:problem:3}.

\begin{prop} \label{prop:leeds:first} 
Let $A$ be an algebra and $r_A$ belong to $A\otimes A$.  
Let $N$ be an integer $\geq 3$ and assume that 
$\cR_1,\ldots,\cR_{N-2}$ in $A\otimes A$ satisfy $\cR_1 = r_A$ and 
\begin{equation} \label{QYBE:low:first}
\lbrack\!\lbrack r_A,\cR_i \rbrack\!\rbrack  
= - \sum_{p+q+r = i+1, p,q,r>0} \cR_p^{12}\cR_q^{13}\cR_r^{23} 
+  \sum_{p+q+r = i+1,p,q,r>0} \cR_r^{23} \cR_q^{13} \cR_p^{12} ,  
\end{equation}
for $i = 1,\ldots,N-2$ (in particular, $r_A$ is a solution of 
CYBE). Then 
\begin{equation} \label{test:term:first}
\delta(r_A|
\sum_{p+q+r = N,p,q,r>0} \cR_p^{12}\cR_q^{13}\cR_r^{23} 
-  \cR_r^{23} \cR_q^{13} \cR_p^{12})
\end{equation}
is zero. Here we set, if $a,b\in A$, 
$$
\lbrack\!\lbrack a,b \rbrack\!\rbrack   
= [a^{12},b^{13}] + [a^{12},b^{23}] + [a^{13},b^{23}]  
+ [b^{12},a^{13}] + [b^{12},a^{23}] + [b^{13},a^{23}],
$$ 
and if $T\in A^{\otimes 3}$, 
\begin{align*}
\delta(a|T) = &
[a^{12} + a^{13} + a^{14}, T^{234}] 
- [- a^{12} + a^{23} + a^{24}, T^{134}] 
\\ & 
+ [-a^{13} - a^{23} + a^{34}, T^{124}] 
+ [a^{14} + a^{24} + a^{34}, T^{123}].
\end{align*}  
\end{prop}

Let us now prove Theorem \ref{training}. 
Let us define $\rho^{(N)}$ as the element 
$\sum_{i = 1}^N \hbar^i\rho_i$ of $\hbar\A^{\otimes 2}[[\hbar]]
/ \hbar^{N+1}\A^{\otimes 2}[[\hbar]]$. Then $\cR^\omega(\rho^{(N)})$ 
is an element of $\Sh^\omega(\A)^{\wh\otimes 3} / (\hbar^{N+1})$, 
and by Proposition \ref{prop:equivalence}, it satisfies QYBE. 
Let us expand  $\cR^\omega(\rho^{(N)})$ as
$\sum_{i = 1}^{N} \hbar^i \cR_i$. We may apply Theorem \ref{training}
with $A = \Sh(\A)$ and $N - 2$  replaced by $N$. Then the conclusion of this Theorem 
says that $\delta(r| \sum_{i,j,k>0, i+j+k = N+1} 
\cR^{12}_i \cR^{13}_j \cR^{23}_k   - 
\cR^{23}_k   \cR^{13}_j  \cR^{12}_i)$ is zero. 
To apply $\pr_\A^{\otimes 4}$ to this 
identity, we use the following Lemma. Let us denote 
by $\iota$ the natural embedding of $\A$ in $\Sh(\A)$, 
sending $x$ to $(x)$. 

\begin{lemma} \label{lemma:proj} 
If $a\in \A$ and $T\in \Sh(\A)$, then $\pr_\A([\iota(a),T]) =  
[a,\pr_\A(T)]$. 
\end{lemma}

{\em Proof.} We may assume that $T$ is homogeneous. 
The statement is obvious when $T$ has degree $\leq 1$. 
Assume that  $T = (x_1\ldots x_n)$ with $n\geq 2$, then 
$\pr_\A([\iota(a),T]) =   B_{1,n}(a|x_1\ldots x_n) - 
B_{n,1}(x_1\ldots x_n|a)$. Set $\beta_n(a|x_1\ldots x_n) = 
B_{1,n}(a|x_1,\ldots,x_n) - B_{n,1}(x_1,\ldots,x_n|a)$. 
Then $\beta_n$ is an element of $FL_{n+1}$ and we 
should prove that $\beta_n(a|x_1\ldots x_n)$ is 
identically zero. Let us prove this 
by induction on $n$. Assume that we have shown that 
when $m<n$, $\beta_m(a|x_1\ldots x_m) = 0$. 
Then we have 
$$
[(a),(x_1\ldots x_n)] = \sum_{i = 1}^n (x_1,\ldots [a,x_i]
\ldots x_n) + (\beta_n(a|x_1\ldots x_n)). 
$$
Then the Jacobi identity implies that 
\begin{align}
& \big( [a',\beta_n(a|x_1,\ldots,x_n)] + \sum_{i = 1}^n 
\beta_n(a'|x_1,\ldots,[a,x_i],\ldots,x_n) \big) 
\\ & 
- \big( [a,\beta_n(a'|x_1,\ldots,x_n)] + \sum_{i = 1}^n 
\beta_n(a|x_1,\ldots,[a',x_i],\ldots,x_n) \big) 
\\ & = \beta_n([a',a]|x_1,\ldots,x_n) ;   
\end{align}
since $\beta_n$ is linear in each argument, this 
implies that $\beta_n([a',a]|x_1,\ldots,x_n) = 0$. 
This is a universal identity, valid in $FL_{n+2}$. 
Therefore the element $\beta_n$ of $FL_{n+1}$ is zero. 
\hfill \qed\medskip 

We have already seen that $\pr_\A^{\otimes 4} \big( 
\delta(r| \sum_{i,j,k>0, i+j+k = N+1} 
\cR^{12}_i \cR^{13}_j \cR^{23}_k   - 
\cR^{23}_k   \cR^{13}_j  \cR^{12}_i) \big)$  is zero. 
It follows from Lemma \ref{lemma:proj} that this is 
the image by $\delta_{4,r}$ of 
$$
\pr_{\A}^{\otimes 3} \big( \sum_{i,j,k>0, i+j+k = N+1} 
\cR^{12}_i \cR^{13}_j \cR^{23}_k   - 
\cR^{23}_k   \cR^{13}_j  \cR^{12}_i \big) , 
$$ 
which is equal to 
$\Phi_{N+1}(\rho_2,\ldots,\rho_{N-1})$. Therefore the image of 
$\Phi_{N+1}(\rho_2,\ldots,\rho_{N-1})$ by $\delta_{4,r}$
is zero. This proves Theorem \ref{training}. 
\hfill \qed\medskip

Theorem \ref{training} has the following ``universal'' counterpart.

\begin{thm} \label{thm:univ}
Assume that $p$ is an integer and  $\varrho_1,\ldots,\varrho_p$ 
belong to $F_1,\ldots,F_p$
and satisfy the universal Lie QYB equations of order $\leq p$
(this is the system of equations (\ref{univ:Lie:QYBE}), 
where $N$ takes the values $1,\ldots,p$). Then the image by 
$\delta_4^{(F)}$ of the right side of equation 
(\ref{induction:Lie:QYBE}), in which 
$N$ takes the value $p+1$, is zero. 
\end{thm}

{\em Proof.} See Appendix \ref{sect:proof:thm:univ}. \hfill \qed\medskip

\subsection{Cohomology groups $H^2_n$ and $H^3_n$}
\label{sect:comp:cohomologies}

\subsubsection{Definition of $H^2_n$ and $H^3_n$}

When $n$ and $N$ are integers, let us set 
$$
F^{Lie,(n)}_N = \bigoplus_{x_1,\ldots,x_{n-2}\in\{a,b\}}
F^{(ax_1\ldots x_{n-2}b)}_N. 
$$
So we have $F^{Lie,(2)}_N = F_N$ and $\oplus_{N\in \NN}
F^{Lie,(n)}_N = \bigoplus_{x_1,\ldots,x_{n-2}\in\{a,b\}}
F^{(ax_1\ldots x_{n-2}b)}$. Then 
$\delta^{(3)}_F$ maps $F^{Lie,(2)}_N$ to $F^{Lie,(3)}_{N+1}$
and  $\delta^{(4)}_F$ maps $F^{Lie,(3)}_N$ to $F^{Lie,(4)}_{N+1}$. 

Let us set 
$$
H^2_N = \Ker(\delta^{(F)}_{3|F^{Lie,(2)}_N}) 
$$
and 
$$
H^3_N = \Ker(\delta^{(F)}_{4|F^{Lie,(2)}_N})  
/ \delta^{(F)}_3(F^{Lie,(2)}_{N-1}). 
$$
Then if we set 
$$
H^2 = \Ker(\delta^{Lie,(F)}_3) \quad\on{and}\quad 
H^3 = \Ker(\delta^{Lie,(F)}_4) / \Imm(\delta^{Lie,(F)}_3),  
$$
we have 
$H^2 = \oplus_{N\geq 0} H^2_N$ and $H^3 = \oplus_{N\geq 0} H^3_N$.

\subsubsection{Results}

\begin{thm} \label{results:cohom}

1) $H^2_N$ is zero if $N\neq 1$, and $H^2_1$ is one-dimensional, spanned
by the class of $r = x_1\otimes y_1$. 

2) $H^3_N$ is zero if $N\neq 2$, and $H^3_2$ is two-dimensional, 
spanned by the classes of $[r^{(13)},r^{(23)}]
\in F_2^{(aab)}$ and $[r^{(12)},r^{(13)}] \in F_2^{(abb)}$.  
\end{thm}

{\em Proof.} See Appendix \ref{app:cohomologies}. 
\hfill \qed\medskip

\subsection{Solution of the universal Lie QYB equations}
\label{sect:thm}

\begin{thm} \label{thm:unique:sol}
There exists a unique solution 
$(\varrho_n)_{n\geq 1}\in \prod_{n\geq 1}F_n$ 
to the universal Lie QYB 
equations (\ref{univ:Lie:QYBE}), such that 
$\varrho_1 = x_1\otimes y_1$. 
\end{thm}

{\em Proof.} $\varrho_2$ should satisfy $\varrho_2\in F_2$ and 
\begin{equation} \label{eq:rho:2}
\delta_3^{(F)}(\varrho_2) = - 2\mu_{Lie}
([w_1,v_1]\otimes [w'_1,u_1]\otimes[v'_1,u'_1]) .  
\end{equation}
Denote by $\varsigma$ the right side of (\ref{eq:rho:2}). Then  
$\varsigma$ belongs to $F^{Lie,(3)}_3$, and it follows from 
Theorem \ref{thm:univ}
that $\delta^{(F)}_4(\varsigma) = 0$. Then the second part of 
Theorem \ref{results:cohom} 
implies the existence of $\varrho_2$ satisfying (\ref{eq:rho:2}), and 
the first part of this Theorem implies the unicity of $\varrho_2$.  
The sequence $(\varrho_n)_{n\geq 3}$ is then constructed inductively 
in the same way. 
\hfill \qed\medskip 

\newpage

\section{Quantization of Lie bialgebras} \label{sect:quantization}

In this Section, we again fix an element $\omega$ of 
$\cB(\KK[[\hbar]])$.  

\subsection{Lie bialgebras} \label{sect:basics}

Let $\G$ be a finite-dimensional Lie bialgebra (we may also 
assume that $\G$ is a positively graded Lie bialgebra, whose homogeneous
components are finite-dimensional). We denote by $[\ ,\ ]_\G
: \G\otimes\G\to \G$ and by $\delta_\G : \G\otimes\G \to \G$
the bracket and cobracket of $\G$. We also 
denote by $[\ ,\ ]_{\G^*}$ and $\delta_{\G^*}$ the bracket 
and cobracket of $\G^*$, so $[\ ,\ ]_{\G^*} = (\delta_\G)^*$ 
and $\delta_{\G^*} = ([\ ,\ ]_\G)^*$. We denote by $\D$
the double Lie bialgebra of $\G$, by $[\ , \, ]_\D$ and 
$\delta_{\D}$ its bracket and cobracket. We have $\D = \G\oplus\G^*$,
$([\ , \, ]_\D)_{|\G\times\G} = [\ , \, ]_\G$,  
$([\ , \, ]_\D)_{|\G^*\times\G^*} = [\ , \, ]_{\G^*}$
$(\delta_\D)_{|\G} = \delta_\G$ and $(\delta_\D)_{|\G^*} = 
- \delta_{\G^*}$. 

Moreover, the symmetric nondegenerate bilinear form 
$\langle\ ,\ \rangle_\D$ defined on $\D$ by 
$\langle (x,\xi), (y,\eta)\rangle_\D = \eta(x) + \xi(y)$ is nondegenerate
and invariant, and if $r_\G\in \G\otimes \G^*$ is the canonical element
of $\G\otimes\G^*$ corresponding to the pairing between $\G$ and
$\G^*$, we have $\delta_\D(x) = [x\otimes 1 + 1\otimes x,r_\G]_\D$.

\subsection{Subalgebras of shuffle algebras}

Recall that $\Sh^\omega(\D)$ is a topological $\KK[[\hbar]]$-Hopf algebra, 
linearly isomorphic to $T(\D)[[\hbar]] = 
\oplus_{k\geq 0}\D^{\otimes k}[[\hbar]]$. Then its subspaces $T(\G)$ 
and $T(\G^*)$ are Hopf subalgebras, isomorphic to $\Sh^\omega(\G)$
and $\Sh^\omega(\G^*)$. 

Let us define $\Sh^\omega_\hbar(\G)$ (resp., $\Sh^\omega_\hbar(\G^*)$)
as the subspace $\wh{\bigoplus}_{k\geq 0} \hbar^k\G^{\otimes k}[[\hbar]]$ 
(resp., $\wh{\bigoplus}_{k\geq 0} \hbar^k(\G^*)^{\otimes k}[[\hbar]]$)
of $\Sh^\omega(\D)$, 
where $\wh{\oplus}$ means the complete direct sum. 
Then $\Sh^\omega_\hbar(\G)$ (resp., $\Sh^\omega_\hbar(\G^*)$) is a topological Hopf 
subalgebra of $\Sh^\omega(\G)$ (resp., of $\Sh^\omega(\G^*)$). 

\begin{remark} If we emphasize the dependence of  the shuffle algebra 
$\Sh^\omega(\A)$ of a Lie algebra $(\A,[\ ,\ ]_\A)$ in the Lie bracket of $\A$ 
by denoting it $\Sh(\A,[\ ,\ ]_\A)$,  then $\Sh^\omega_\hbar(\G)$ may be viewed as a
completion of $\Sh^\omega(\G[[\hbar]],\hbar[\ ,\ ]_\G)$. 
\end{remark}

\subsection{Tensor algebra of $(\G,\delta_\G)$ and Hopf pairing}

Recall that the Hopf algebra $T^\omega_\hbar(\G)$ 
is the vector space $T(\G)[[\hbar]]$, 
equipped with the undeformed multiplication  
$m_{T(\G)}$ the tensor algebra
$T(\G)$ and the comultiplication defined by 
``reversing the arrows'' in the definition of the shuffle algebra 
$\Sh^\omega(\A)$ (Proposition \ref{prop:Ta}). 

There is a unique bilinear map 
$$
\langle\  ,\ \rangle_{
\Sh^\omega(\G^*)\times T^\omega_\hbar\G} : 
\Sh^\omega(\G^*) \times T^\omega_\hbar\G
\to \KK[[\hbar]][\hbar^{-1}] , 
$$
such that if $\xi_1,\ldots,\xi_m$
belong to $\G^*$ and $x_1,\ldots,x_n$ belong to $\G$, 
then $\langle (\xi_1\cdots \xi_m) 
, x_1\otimes \cdots \otimes x_n 
\rangle_{\Sh^\omega(\G^*) \times T^\omega_\hbar\G} = 
\hbar^{-n}\delta_{n,m} \prod_{i = 1}^n \langle \xi_i, x_i\rangle_{
\G^*\times\G}$, 
and $\langle \ ,\ \rangle_{\G^*\times\G}$ is the canonical 
pairing between $\G^*$ and $\G$. 

Then $\langle\  ,\ \rangle_{\Sh^\omega(\G^*) \times T^\omega_\hbar\G}$ 
is a Hopf pairing, which means that we have 
\begin{equation} \label{duality:1}
\langle \xi\eta , x 
\rangle_{\Sh^\omega(\G^*) \times T^\omega_\hbar\G}
=  \sum_i 
\langle \xi, x^{(1)}  \rangle_{\Sh^\omega(\G^*) \times T^\omega_\hbar\G} 
\langle \eta, x^{(2)} \rangle_{\Sh^\omega(\G^*) \times T^\omega_\hbar\G} 
\end{equation}
and 
\begin{equation} \label{duality:2}
\langle \xi, xy  \rangle_{\Sh^\omega(\G^*) \times T^\omega_\hbar\G} 
=  \sum_i 
\langle \xi^{(1)},x \rangle_{\Sh^\omega(\G^*) \times T^\omega_\hbar\G} 
\langle \xi^{(2)},x \rangle_{\Sh^\omega(\G^*) \times T^\omega_\hbar\G} 
\end{equation}
for any $\xi,\eta$ in 
$\Sh^\omega(\G^*)$ and any $x,y$ in $T^\omega_\hbar(\G)$, 
where we set $\Delta_{T^\omega_\hbar(\G)}(x) = \sum 
x^{(1)} \otimes x^{(2)}$ and $\Delta_{\Sh^\omega(\G^*)}(\xi) = 
\sum \xi^{(1)}\otimes \xi^{(2)}$.

\subsection{$R$-matrix and associated Hopf algebra morphism}

Recall that $r_\G$ is an element of $\G\otimes\G^*\subset \D\otimes\D$ 
and let us set 
$$
\cR = \cR^\omega(\sum_{n\geq 0} \rho^\omega_n(\hbar r_\G)) .  
$$
Then $\cR$ belongs to 
$\Sh^\omega(\G)\wh\otimes\Sh^\omega(\G^*)$, which is a subalgebra of 
$\Sh^\omega(\D)\wh\otimes\Sh^\omega(\D)$. It follows from Proposition
\ref{prop:QT} and Remark \ref{rem:QT} that $\cR$ satisfies the 
quasitriangularity identities
\begin{equation} \label{QT:part}
(\Delta_{\Sh^\omega(\G)}\otimes \id_{\Sh^\omega(\G^*)})(\cR) = \cR^{(13)}\cR^{(23)},  
\quad 
(\id_{\Sh^\omega(\G)} \otimes \Delta_{\Sh^\omega(\G^*)})(\cR) = \cR^{(13)}\cR^{(12)} 
\end{equation}
and
\begin{equation} \label{behavior:S:new}
(S_{\Sh^\omega(\G)} \otimes \id_{\Sh^\omega(\G^*)}) (\cR) 
= (\id_{\Sh^\omega(\G)}\otimes S^{-1}_{\Sh^\omega(\G^*)})(\cR) ,  
\end{equation}
and from Theorem \ref{big:thm:QYBE} that it satisfies the QYBE  
\begin{equation} \label{QYBE:part}
\cR^{(12)}\cR^{(13)}\cR^{(23)} = \cR^{(23)}\cR^{(13)}\cR^{(12)}. 
\end{equation}

\begin{lemma}
The rule 
$$
\ell(x) = \langle \cR , \id\otimes x
\rangle_{\Sh^\omega(\G^*)\times T^\omega_\hbar(\G)} 
$$
for $x\in T^\omega_\hbar(\G)$ defines a linear map $\ell$ from $T^\omega_\hbar(\G)$
to $\Sh^\omega(\G)$. 
\end{lemma}

{\em Proof.} Let us write $\cR = \sum_{n\geq 0} \hbar^n \cR_n$,
then $\cR_n$ has the following form 
$$
\cR_n\in \sum_{i_1,\ldots,i_n\in I,\sigma\in\SG_n} 
(a_{i_1}\cdots a_{i_n})\otimes (b_{i_{\sigma(1)}}\cdots b_{i_{\sigma(n)}}) + 
\bigoplus_{(k,k')|k\leq n,k'\leq n, (k,k')\neq (n,n)}
\G^{\otimes k}\otimes (\G^*)^{\otimes k},  
$$
where $(a_i)_{i\in I}$ is a basis of $\G$ and $(b_i)_{i\in I}$ is the dual basis 
of $\G^*$. Then if $x_1,\ldots,x_n$ are in $\G$ and $x = x_1\otimes\cdots \otimes x_k$, 
$\langle \cR_k, \id\otimes x
\rangle_{\Sh^\omega(\G^*)\times T^\omega_\hbar(\G)}
= 0$ if $k<n$, so $\ell(x)$ has nonnegative $\hbar$-adic valuation. 
More precisely, we have
$$
\cR_n\in \sum_{i_1,\ldots,i_n\in I} (a_{i_1})\cdots (a_{i_n}) \otimes
(b_{i_1}\cdots b_{i_n}) + (\bigoplus_{k|k\leq n}\G^{\otimes k})
\otimes (\bigoplus_{k'|k'<n} (\G^*)^{\otimes k'}) , 
$$
so $\langle \cR_n, \id\otimes x
\rangle_{\Sh^\omega(\G^*)\times T^\omega_\hbar(\G)}
= \hbar^{-n}(x_1)\cdots(x_n)$. Therefore 
\begin{equation} \label{reduction}
\ell(x_1\otimes\cdots \otimes x_k) = (x_1)\cdots(x_k) + O(\hbar).  
\end{equation}
\hfill \qed\medskip 

Then it follows from equations (\ref{QT:part}), 
(\ref{behavior:S:new}), 
(\ref{duality:1}) and (\ref{duality:2}) that $\ell$ is a morphism 
of Hopf algebras from $T^\omega_\hbar(\G)^{opp}$ to $\Sh^\omega(\G)$
($T^\omega_\hbar(\G)^{opp}$ is the opposite algebra of $T^\omega_\hbar(\G)$, 
endowed with the same coproduct as $T^\omega_\hbar(\G)$). 

\subsection{Construction of $U^\omega_\hbar\G$}

Let us set $U^\omega_\hbar\G = \Imm(\ell)$. Since $\ell$ is a 
morphism of Hopf algebras, $U^\omega_\hbar\G$ is a 
Hopf subalgebra of $\Sh^\omega(\G)$. We will denote by 
$m_{U^\omega_\hbar\G}$ and $\Delta_{U^\omega_\hbar\G}$ the product and
coproduct of $U^\omega_\hbar\G$, and by $\varepsilon_{U^\omega_\hbar\G}$
and $S_{U^\omega_\hbar\G}$ its counit and antipode.    

\begin{thm} \label{quantization:QLBA}
$(U^\omega_\hbar\G,m_{U^\omega_\hbar\G},\Delta_{U^\omega_\hbar\G},\varepsilon_{U^\omega_\hbar\G},
S_{U_\hbar\G})$ is a quantization of $(\G,[\ ,\ ]_\G,\delta_\G)$. 

Moreover,  $\Imm(\ell)$ is a divisible submodule of $\Sh^\omega(\G)$, 
i.e.\ $\Imm(\ell)\cap \hbar\Sh^\omega(\G) = \hbar\Imm(\ell)$. 
\end{thm}

{\em Proof.}
Let us show that $U^\omega_\hbar\G / \hbar U^\omega_\hbar\G$ is isomorphic to 
$U\G$. $U^\omega_\hbar\G$ is isomorphic, as a Hopf algebra, 
to $T^\omega_\hbar(\G)^{opp} / \Ker(\ell)$, 
therefore $U^\omega_\hbar\G / \hbar U^\omega_\hbar\G$ is isomorphic
to the Hopf algebra $T^\omega_\hbar(\G)^{opp} / 
\big( \Ker(\ell) + \hbar T^\omega_\hbar(\G)^{opp} \big)$.  To identify 
the latter Hopf algebra with $U\G$, let us study the kernel 
$\Ker(\ell)$.  

\begin{lemma}
For any integers $p$ and $q$, there exists unique $(q+1)$-linear maps 
$\beta_{pq}$ and $\gamma_{qp} : \G^{q+1}\to \G^{\otimes p}[[\hbar]]$, such 
that for any $\xi_1,\ldots,\xi_p$
in $\G^*$ and any $x,x_1,\ldots,x_q$ in $\G$, we have the equalities 
$$
\langle B_{pq}(\xi_1,\ldots,\xi_p|x_1,\ldots,x_q), x\rangle_\D
= \langle \xi_1\otimes\cdots\otimes \xi_p, \beta_{pq}(x_1,\ldots,x_q,x) 
\rangle_{\D^{\otimes p}}
$$
and 
$$
\langle B_{qp}(x_1,\ldots,x_q|\xi_1,\ldots,\xi_p), x\rangle_\D
= \langle \xi_1\otimes\cdots\otimes \xi_p, \gamma_{qp}(x_1,\ldots,x_q,x) 
\rangle_{\D^{\otimes p}}
$$
(we view $\G$ and $\G^*$ as subspaces of $\D$ via the maps $x\mapsto 
(x,0)$ and $\xi\mapsto (0,\xi)$, and $\langle\ ,\ \rangle_{\D^{\otimes p}}$ 
is the $p$th tensor power of $\langle\ ,\ \rangle_{\D}$). 
\end{lemma}

Extend the pairing 
$\langle\ ,\ \rangle_{\Sh^\omega(\G^*)\times T^\omega_\hbar(\G)}$
to a bilinear map 
$$
\langle\ ,\ \rangle_{\Sh^\omega(\D)\times T^\omega_\hbar(\G)}:  
\Sh^\omega(\D)\times T^\omega_\hbar(\G) \to \KK[[\hbar]][\hbar^{-1}], 
$$
by the rule that if $x_1,\ldots,x_n$ are elements of $\D$, 
one of which belongs to $\G$, and $y$ is any element of 
$T^\omega_\hbar(\G)$, then $\langle (x_1\cdots x_n), y 
\rangle_{\Sh^\omega(\D)[[\hbar]]\times T^\omega_\hbar(\G)} = 0$. 

\begin{lemma} \label{dualnost}
There exist unique bilinear maps $\phi$ and $\psi : 
\Sh^\omega(\G)\times T^\omega_\hbar(\G) \to T^\omega_\hbar(\G)$ 
such that 
for any $(x,x',y)\in \Sh^\omega(\G) \times \Sh^\omega(\G^*)
\times T^\omega_\hbar(\G)$, we have 
$$
\langle x' x,  y\rangle_{\Sh^\omega(\D)\times
T^\omega_\hbar(\G)} =
\langle x' , \phi(x,y) \rangle_{\Sh^\omega(\D)\times
T^\omega_\hbar(\G)}
$$
and
$$ 
\langle  xx', y\rangle_{\Sh^\omega(\D)\times
T_\hbar^\omega(\G)} = 
\langle x' , \psi(x,y) \rangle_{\Sh^\omega(\D)\times
T_\hbar^\omega(\G)}. 
$$
Identify $\G$ with a subspace of $T^\omega_\hbar(\G)$ of  
tensors of degree $1$. Then if $x$ and $y$ belong to $\G$, 
then 
$$
\phi(x,y)\in {1\over 2}[x,y]_\G + \hbar T^\omega_\hbar(\G), \quad 
\psi(x,y)\in - {1\over 2}[x,y]_\G + \hbar T^\omega_\hbar(\G).  
$$
\end{lemma}

{\em Proof.} If $n$ and $m$ are integers $\geq 0$, if $x_1,\ldots,x_n$
and $y_1,\ldots,y_m$ are elements of $\G$, 
and if $x = (x_1\cdots x_n)$ and $y = y_1\otimes \cdots\otimes y_m$, 
we set 
\begin{align*}
& \phi(x,y) = \sum_{(\la_1,\ldots,\la_m)\on{\ partition\ of\ }n}
\sum_{k_1,\ldots,k_m >0} 
\\ & \hbar^{k_1 + \cdots + k_m - m}
\beta_{k_1\la_1}(x_1,\ldots,x_{\la_1},y_1)\otimes\cdots 
\otimes\beta_{k_m\la_m}(x_{\la_1 + \cdots + \la_{m-1}+1},\ldots,
x_{\la_1 + \cdots + \la_m},y_m)
\end{align*}
and 
\begin{align*}
& \psi(x,y) = \sum_{(\la_1,\ldots,\la_m)\on{\ partition\ of\ }n}
\sum_{k_1,\ldots,k_m >0} 
\\ & \hbar^{k_1 + \cdots + k_m - m}
\gamma_{\la_1 k_1}(x_1,\ldots,x_{\la_1},y_1)\otimes\cdots 
\otimes\gamma_{\la_m k_m}(x_{\la_1 + \cdots + \la_{m-1}+1},\ldots,
x_{\la_1 + \cdots + \la_m},y_m) . 
\end{align*}
In both right hand sides, only nonnegative powers of $\hbar$ occur, 
and the coefficient of each power of $\hbar$ is a finite sum, so 
that both right sides belong to $T^\omega_\hbar(\G)$. 
It is easy to check that these are the unique maps satisfying the
above requirements. The last statements follow from the equalities 
$\beta_{11}(x,y) = {1\over 2}[x,y]_\G$ and $\gamma_{11}(x,y) = 
- {1\over 2}[x,y]_\G$. Both equalities follow from the 
invariance of $\langle\ ,\ \rangle_\D$. 
\hfill \qed\medskip

\begin{lemma} \label{lemma:phi}
If $x$ and $y$ belong to $T^\omega_\hbar(\G)$, and we write 
$\Delta_{T^\omega_\hbar(\G)}(y) = \sum y^{(1)}\otimes y^{(2)}$, then   
$$
\sum y^{(1)} \phi(\ell(y^{(2)}),x) - 
\sum \psi(\ell(y^{(1)}), x) y^{(2)} 
$$
belongs to $\Ker(\ell)$. 
\end{lemma}

{\em Proof.} We have 
\begin{align*}
& \langle \cR^{(12)}\cR^{(13)}\cR^{(23)}, \id\otimes x\otimes y
\rangle_{\Sh^\omega(\D)\times T^\omega_\hbar(\G)} 
\\ & =
\sum  \langle \cR^{(12)}\cR^{(13)}\cR^{(24)}, \id\otimes x\otimes 
y^{(1)} \otimes y^{(2)}
\rangle_{\Sh^\omega(\D)\times T^\omega_\hbar(\G)}
\\ & 
= \sum \langle \cR^{(12)} (1\otimes \ell(y^{(2)})), \id\otimes x
\rangle_{\Sh^\omega(\D)\times T^\omega_\hbar(\G)} \ell(y^{(1)}). 
\end{align*}
Lemma \ref{dualnost} implies the identity 
$$
\langle \cR(1\otimes x), \id \otimes y\rangle_{\Sh^\omega(\D)\times
T^\omega_\hbar(\G)} = \ell(\phi(x,y)) , 
$$
and since $\ell$ is an algebra antihomomorphism from $T^\omega_\hbar(\G)$ to 
$\Sh^\omega(\G)$, we get 
$$
\langle \cR^{(12)}\cR^{(13)}\cR^{(23)}, \id\otimes x\otimes y
\rangle_{\Sh^\omega(\D)\times T^\omega_\hbar(\G)} 
= \ell(\sum y^{(1)} \phi(\ell(y^{(2)}), x) ) .  
$$
In the same way, 
$\langle \cR^{(23)}\cR^{(13)}\cR^{(12)}, \id\otimes x\otimes y
\rangle_{\Sh^\omega(\D)\times T^\omega_\hbar(\G)}
= \ell( \sum \psi(\ell(y^{(1)}), x) y^{(2)})$. 
Since $\cR$ satisfies QYBE, we have 
$\ell(\sum y^{(1)} \phi(\ell(y^{(2)}), x) )  = 
\ell( \sum \psi(\ell(y^{(1)}), x) y^{(2)})$. 
\hfill \qed\medskip 

{\em End of proof of Theorem \ref{quantization:QLBA}.} 
If  $x\in \G\subset T^\omega_\hbar(\G)$, then $\Delta_{T^\omega_\hbar(\G)}(x)
= x\otimes 1 + 1\otimes x + O(\hbar)$. 
Then applying Lemma \ref{lemma:phi} to the case when $x$ and
$y$ belong to $\G$ and using the end of 
Lemma \ref{dualnost}, we construct 
bilinear maps $m_n : \G\times\G\to T(\G)$, where
$n\geq 1$, such that if $x$ and $y$ belong to $\G$, 
\begin{equation} \label{relations}
y\otimes x - x \otimes y - [x,y]_\G - \sum_{n|n\geq 1} \hbar^n m_n(x,y)
\end{equation}
belongs to $\Ker(\ell)$.  

Let us denote by $I_0$ the complete 
two-sided ideal 
of $T_\hbar^\omega(\G)^{opp}$ generated by the elements (\ref{relations}). 
Then $I_0\subset \Ker(\ell)$. We have therefore a surjective 
morphism of $\KK[[\hbar]]$-algebras 
$T^\omega_\hbar(\G)^{opp} / I_0 \to T^\omega_\hbar(\G)^{opp} / \Ker(\ell)$. The reduction 
modulo $\hbar$ of this map is also surjective; it is a map 
$T_\hbar^\omega(\G)^{opp} / (I_0 + \hbar T^\omega_\hbar(\G)^{opp}) \to T^\omega_\hbar(\G)^{opp} / (\Ker(\ell) + 
\hbar T^\omega_\hbar(\G)^{opp})$. Due to the form of $I_0$, 
$T^\omega_\hbar(\G)^{opp} / (I_0 + \hbar T^\omega_\hbar(\G)^{opp})$ is isomorphic to $U\G$. 
As a result, we obtain that there is a unique surjective morphism 
$s: U\G\to T^\omega_\hbar(\G)^{opp} / (\Ker(\ell) + \hbar T^\omega_\hbar(\G)^{opp})$, such that for 
any $x\in\G$, $s(x) =$ the class of $x$. 

On the other hand, recall that $U\G$ is also the subalgebra of $\Sh^\omega(\G)$
generated by the elements of degree $1$. $\Imm(\ell)$ is a 
subalgebra of $\Sh^\omega(\G)$. It follows from equation 
(\ref{reduction}) that the image of $\Imm(\ell)$ by the 
morphism $\Sh^\omega(\G) \to \Sh^\omega(\G)$ given by the reduction 
modulo $\hbar$ is exactly $U\G$. This means that  we have an 
isomorphism of algebras between $U\G$ and $\Imm(\ell) / 
(\hbar\Sh^\omega(\G)\cap \Imm(\ell))$. The latter 
algebra is a quotient algebra of $\Imm(\ell) / \hbar\Imm(\ell)$. 
We obtain that there exists a surjective algebra morphism 
$s': \Imm(\ell) / \hbar \Imm(\ell) \to U\G$, such that for any 
$x\in\G$, $s'($the class of $\ell(x)) = x$. 

Moreover, $\ell$ induces an isomorphism between $T^\omega_\hbar(\G)^{opp}/\Ker(\ell)$ and
$\Imm(\ell)$, so its reduction modulo $\hbar$ induces 
an isomorphism between $T^\omega_\hbar(\G)^{opp} / (\Ker(\ell) + \hbar T^\omega_\hbar(\G)^{opp})$ and 
$\Imm(\ell) / \hbar \Imm(\ell)$. Denote by $\bar\ell$ this reduction, 
then we have $\bar\ell = s\circ s'$. Since $s$ and $s'$ is
surjective, $s$ is an 
isomorphism between $U\G$ and $\Imm(\ell) / \hbar \Imm(\ell)$.  In other
words, we have shown that $U^\omega_\hbar\G / \hbar U^\omega_\hbar\G$ is isomorphic to 
$U\G$. 

Moreover, $U^\omega_\hbar\G$ is $\hbar$-adically complete, and 
it is torsion-free, because it is a $\KK[[\hbar]]$-submodule 
of $\Sh^\omega(\G)$. $U^\omega_\hbar\G$ is therefore a topologically
free $\KK[[\hbar]]$-algebra, such that $U^\omega_\hbar\G / U^\omega_\hbar\G
= U\G$. 

Let us study now the co-Poisson structure on $U\G$ induced by 
this isomorphism. Consider $x\in \G$ as an element of $T^\omega_\hbar(\G)^{opp}$. 
Then 
$$
{1\over{\hbar}}( \Delta_{T^\omega_\hbar(\G)}(x) - \Delta'_{T^\omega_\hbar(\G)}(x))
\in \delta_\G(x) + \hbar T^\omega_\hbar(\G)\wh\otimes T^\omega_\hbar(\G),
$$
where $\hat\otimes$ is the $\hbar$-adically completed tensor product.  
Taking the image by $\ell$ of this identity, we find 
$$
{1\over{\hbar}}( \Delta_{U^\omega_\hbar \G}(\ell(x)) - \Delta'_{U^\omega_\hbar(\G)}(\ell(x)))
\in \ell^{\otimes 2}(\delta_\G(x)) + \hbar U^\omega_\hbar(\G)\wh\otimes U^\omega_\hbar(\G),
$$
which means that the co-Poisson structure on $U\G$ corresponding to 
$U^\omega_\hbar\G$ is given by $\delta_\G$. 

Let us now prove that 
$\Imm(\ell)$ is a divisible submodule of $\Sh^\omega(\G)$. 
We have shown that the map $s'$ is an isomorphism, which means 
that the surjective morphism $\Imm(\ell) / \hbar\Imm(\ell)
\to \Imm(\ell) / (\Imm(\ell) \cap \hbar\Sh^\omega(\G))$
is an isomorphism. This implies that 
$\Imm(\ell)\cap \hbar\Sh^\omega(\G) = \hbar\Imm(\ell)$.  
\hfill\qed\medskip 

\subsection{Functoriality}

\begin{prop} \label{prop:functor}
The map $(\G,[\ ,\ ]_\G,\delta_\G)\mapsto U^\omega_\hbar\G$ defines a universal 
quantization functor from the
category of finite-dimensional Lie bialgebras. 
\end{prop}  

{\em Proof.} Let $\phi$ be a morphism of Lie bialgebras 
from $(\G,[\ ,\ ]_\G,\delta_\G)$ to $(\HH,[\ ,\ ]_\HH,\delta_\HH)$. 
Then $\phi$ induces Lie algebra morphisms $\phi_{\G\HH} : 
\G\to \HH$ and $\phi_{\HH^*\G^*} : \HH^*\to \G^*$. The first morphism
induces a Hopf algebra morphism $\Sh(\phi_{\G\HH}) : \Sh^\omega(\G)\to\Sh^\omega(\HH)$, 
and the dual to $\Sh^\omega(\phi_{\HH^*\G^*})$ induces a Hopf algebra morphism 
$T(\phi_{\G\HH}) : T^\omega_\hbar(\G)\to T^\omega_\hbar(\HH)$. Let us denote by 
$\cR_\G$ and $\cR_\HH$ the analogues of $\cR$ for $\G$ and $\HH$; then 
$\cR_\G$ belongs to $\Sh^\omega(\G)\wh\otimes \Sh^\omega(\G^*)$ and 
$\cR_\HH$ belongs to $\Sh^\omega(\HH)\wh\otimes \Sh^\omega(\HH^*)$. 
Moreover, if $r_\G$ and $r_\HH$ are the canonical $r$-matrices of 
$\G$ and $\HH$, we have  
$(\id\otimes \phi_{\HH^*\G^*})(r_\G) = (\phi_{\G\HH}\otimes \id)(r_\HH)$, 
therefore  
$$
(\id\otimes \Sh(\phi_{\HH^*\G^*}))(\cR_\HH)= 
(\Sh(\phi_{\G\HH})\otimes \id)(\cR_\G).
$$ 
Let us denote by $\ell_\G$ and $\ell_\HH$ the analogues of $\ell$ corresponding to 
$\G$ and $\HH$. If $x\in T^\omega_\hbar\G$, we have 
\begin{align*}
 \Sh(\phi_{\G\HH})(\ell_\G(x)) & = \langle \id\otimes x, 
(\Sh(\phi_{\G\HH})\otimes \id)(\cR_\G) \rangle_{T^\omega_\hbar(\G)\times\Sh^\omega(\G^*)}
\\ & 
= \langle \id\otimes x, 
(\id\otimes \Sh(\phi_{\HH^*\G^*}))(\cR_\HH) \rangle_{T^\omega_\hbar(\G)\times\Sh^\omega(\G^*)}
\\ & 
= \langle \id\otimes T(\phi_{\G\HH})(x), 
\cR_\HH \rangle_{T^\omega_\hbar(\HH)\times\Sh^\omega(\HH^*)}
= \ell_\HH( T(\phi_{\G\HH})(x) ) , 
\end{align*}
so 
$\Sh(\phi_{\G\HH}) \circ \ell_\G = \ell_\HH \circ T(\phi_{\G\HH})$. 
Therefore the restriction of $\Sh(\phi_{\G\HH})$ to $U^\omega_\hbar\G$ 
induces a Hopf algebra morphism from $U^\omega_\hbar\G$ to $U^\omega_\hbar\HH$. 
Let us denote by $\phi^U$ this morphism; it is then clear that the reduction 
mod $\hbar$ of $\phi^U$ coincides with the morphism from $U\G$ to 
$U\HH$ induced by $\phi$. Moreover, if $\psi:\HH\to\K$ is a morphism of
Lie balgebras, we have $\Sh(\psi\circ \phi) = \Sh(\psi)\circ\Sh(\psi)$; 
the restriction of this identity to $U^\omega_\hbar\G$ yields $(\psi\circ\phi)^U 
= \psi^U\circ \phi^U$. 

Finally, the form taken by the relations (\ref{relations}) shows that 
the quantization functor 
$\G\mapsto U_\hbar^\omega\G$ is universal. 
\hfill \qed\medskip 

\begin{remark}
The condition $p_{kl} = 0$ if $k>l$ in the definition of the spaces 
$F^{(x_1\ldots x_n)}$ (see (\ref{big:lie:space})) 
seems to imply that the ideal 
$I_0$ generated by elements (\ref{relations}) is defined in 
terms of acyclic tensor calculus, in the sense of 
\cite{unsolved}. \hfill \qed\medskip 
\end{remark}

\begin{remark} \label{rem:inf:dim}
When $\G$ is infinite-dimensional, 
the maps $m_P : \G\to\G^{\otimes n}$ defined by 
$m_P(x) = \sum_{i_1,\ldots,i_n\in I}
\langle x, P(b_{i_1},\ldots,b_{i_n})\rangle_{\G\times \G^*} 
a_{i_1}\otimes \cdots\otimes a_{i_n}$, where $P$ is a Lie polynomial,  
do not make sense any more. However, if $\G$ is finite-dimensional, 
these maps are linear combinations of the maps 
$T_\sigma\circ 
(\delta\otimes \id^{\otimes (n-1)})\circ\cdots (\delta\otimes \id)\circ\delta$, 
where $\sigma\in \SG_n$ and $\sigma\mapsto T_\sigma$ is the action of the group 
$\SG_n$ by permutation of the factors of $\G^{\otimes n}$. It is 
easy to see that the map $\ell$ only involves linear combinations of 
compositions of the maps $m_P$ and of the Lie bracket of $\G$. It has 
therefore a natural analogue when $\G$ is an infinite-dimensional 
Lie bialgebra. It is natural to expect that the corresponding 
analogue of $\G\mapsto U^\omega_\hbar\G$  defines a quantization functor for the category of 
(possibly infinite-dimensional) Lie bialgebras. 
\hfill \qed\medskip 
\end{remark}

\subsection{The QFSH algebra of $U^\omega_\hbar\G$}

To any quantized universal enveloping algebra $U_\hbar\A$, one associates
in a canonical way a quantized formal series Hopf (QFSH) algebra $\cO_{\hbar}(A^*)$
(see \cite{QG,Gav}).  If $(\A,[\ ,\ ]_\A,\delta_\A)$ is the Lie bialgebra
corresponding to the semiclassical limit of $U_\hbar\A$, then $\cO_\hbar(A^*)$
is a quantization of the formal series Hopf algebra of functions on the 
formal group $A^*$ corresponding to the Lie algebra $(\A^*,\delta_\A^*)$, 
endowed with the Poisson-Lie structure corresponding to the cobracket $[\ ,\ ]_\A^*$. 

(For example, if $\delta_\A = 0$, then $A^*$ is the additive group 
$\A^*$, endowed with the Kostant-Kirillov Poisson bracket. If moreover
$U_\hbar\A$ is $U\A[[\hbar]]$, then $\cO_\hbar(A^*)$ is a deformation 
quantization of the formal series ring $S[[\A]]$ corresponding to 
the Kostant-Kirillov bracket, endowed with the cocommutative coproduct
such that the elements of $\A$ are primitive.)

The purpose of this Section is to express the QFSH algebra of $U^\omega_\hbar\G$
in terms of the above construction. 

Recall that we defined $\Sh^\omega_\hbar(\G)$ as the subalgebra 
${\wh{\bigoplus}}_{k\geq 0} \hbar^k\G^{\otimes k}[[\hbar]]$ of $\Sh^\omega(\G)$. 

\begin{prop} \label{prop:QFSH}
Let us define $\cO^\omega_\hbar(G^*)$ as the intersection 
$\Imm(\ell)\cap \Sh^\omega_\hbar(\G)$. 
Then $\cO^\omega_\hbar(G^*)$ is the QFSH algebra of $U^\omega_\hbar\G$. 
\end{prop}

{\em Proof.}
Let us first recall how the QFSH algebra of $U^\omega_\hbar\G$ is defined. According to 
\cite{QG,Gav}, one defines a functor $H\mapsto H'$ in the category of 
topologically free Hopf algebras over $\KK[[\hbar]]$. 
If $H$ is such an algebra, let us denote by $\Delta_H$ its coproduct, 
by $\Delta_H^{(n)}$ its $n$th fold coproduct  
and by $\varepsilon_H$ its counit. Let us set $\delta_n^H 
= (\id_H - \varepsilon_H)^{\otimes n}\circ \Delta_H^{(n)}$. Then $H'$ is defined as 
$\{h\in H|\forall n\geq 0, \delta_n^H(h)\in \hbar^n H^{\wh\otimes n} \}$.
One shows (\cite{Gav}) that $H'$ is then a QFSH algebra.  
We first show

\begin{lemma} \label{QFSH:Sh}
$(\Sh^\omega(\G))' = \Sh^\omega_\hbar(\G)$. 
\end{lemma}

{\em Proof.} Let $k$ and $n$ be integers. Recall that $\Sh^\omega(\G)$ is 
identified, as a vector space, with $T(\G)[[\hbar]]$ and we denote by $\conc$
the concatenation product on $T(\G)$. Denote by $\conc^{(n)}$ the $n$fold 
concatenation product; $\conc^{(n)}$ is a linear map from $T(\G)^{\otimes n}$
to $T(\G)$. Then if $x\in \G^{\otimes k}\subset \Sh^\omega(\G)$, 
we have 
$$
\conc^{(n)} \circ \delta^{\Sh^\omega(\G)}_n(x) = s(n,k)x, 
$$
where $s(n,k)$ is the number of ordered surjections from $\{1,\ldots,k\}$
to $\{1,\ldots,n\}$; so $s(n,k) = 0$ if $n<k$ and $s(n,k)>0$ else. 

Let $x$ be an element of $\Sh^\omega(\G)$. 
Set $x = \sum_{k\geq 0} x_k$, where $x_k\in \G^{\otimes k}[[\hbar]]$. 
Then if $x\in (\Sh^\omega(\G))'$, and $\nu$ is any integer, 
$\delta^{\Sh^\omega(\G)}_{\nu}(x)\in \hbar^{\nu}\Sh^\omega(\G)^{\wh\otimes \nu}$. 
Applying $\conc^{(\nu)}$ to this inclusion, we find $\sum_{n\geq 0}
s(n,k) x_k \in \hbar^{\nu}\Sh^\omega(\G)$, therefore 
$x_k\in \hbar^\nu \G^{\otimes k}[[\hbar]]$  if $k\geq \nu$. So for each $k$, 
$x_k\in \hbar^k \G^{\otimes k}[[\hbar]]$, which means that $x\in \Sh^\omega_\hbar(\G)$. 

Let us now show that $\Sh^\omega_\hbar(\G) \subset (\Sh^\omega(\G))'$. Let 
$x$ belong to $\hbar^k \G^{\otimes k}[[\hbar]]$, then we have, if $n\leq k$ 
$\delta^{\Sh^\omega(\G)}_n(x) \subset \hbar^k \Sh^\omega(\G)^{\wh\otimes n}
\subset \hbar^n \Sh^\omega(\G)^{\wh\otimes n}$, and if $n>k$, 
$\delta^{\Sh^\omega(\G)}_n(x) = 0$ so that $\delta^{\Sh^\omega(\G)}_n(x)$ 
is again contained in  $\hbar^n \Sh^\omega(\G)^{\wh\otimes n}$. So 
$x\in (\Sh^\omega(\G))'$. 
\hfill \qed\medskip 

{\em End of proof of Proposition \ref{prop:QFSH}.}
We should prove that $(\Imm(\ell))' = \Imm(\ell)\cap\Sh^\omega_\hbar(\G)$. 
By definition, $(\Imm(\ell))' = \{x\in \Imm(\ell) | \forall n\geq 0, 
\delta^{\Imm(\ell)}_n(x)\in \hbar^n\Imm(\ell)^{\wh\otimes n}\}$. 

It follows from the fact that 
$\Imm(\ell)$ is a divisible submodule of $\Sh^\omega(\G)$
(see Theorem \ref{quantization:QLBA}) that 
$\Imm(\ell)^{\wh\otimes n}\cap \hbar^n \Sh^\omega(\G)^{\wh\otimes n}
= \hbar^n \Imm(\ell)^{\wh\otimes n}$. Therefore 
$(\Imm(\ell))'$ is the same as 
$\{x\in \Imm(\ell) | \forall n\geq 0, 
\delta^{\Sh^\omega(\G)}_n(x)\in \Imm(\ell)^{\wh\otimes n}
\cap \hbar^n \Sh^\omega(\G)^{\wh\otimes n}\}$. 
Since $\delta^{\Sh^\omega(\G)}_n(\Imm(\ell))
\subset \Imm(\ell)^{\wh\otimes n}$, this set is the same as 
$\{x\in \Imm(\ell) | \forall n\geq 0, 
\delta^{\Sh^\omega(\G)}_n(x)\in 
\hbar^n \Sh^\omega(\G)^{\wh\otimes n}\}$, 
which is $\Imm(\ell)\cap (\Sh^\omega(\G))'$. It then 
follows from Lemma \ref{QFSH:Sh} that this is 
$\Imm(\ell)\cap \Sh^\omega_\hbar(\G)$. 
\hfill \qed\medskip

The dual $(U^\omega_\hbar\G)^* = \on{Hom}_{\KK[[\hbar]]}(U^\omega_\hbar\G,\KK[[\hbar]])$
is also a QFSH algebra associated to the formal group $G$. If we denote by 
$\ell_{\G^*}$ the analogue of $\ell$ for the Lie bialgebra $\G^*$, we get therefore  
two quantization functors from the category of Lie bialgebra to that of QFSH algebras,
namely $\G\mapsto (U^\omega_\hbar\G)^*$ and $\G\mapsto \cO^\omega_\hbar(G) = 
\Imm(\ell_{\G^*})\cap \Sh^\omega_\hbar(\G)$. 
Let us denote by $A\to A^\vee$ the functor associating to each 
Hopf algebra, its enveloping QUE (quantized universal enveloping)
algebra.

\begin{prop} \label{rem:dual}
The QFSH algebras $(U^\omega_\hbar\G)^*$ and $\cO^\omega_\hbar(G) = 
\Imm(\ell_{\G^*})\cap \Sh^\omega_\hbar(\G)$ 
are canonically isomorphic. The QUE algebras
$U^\omega_\hbar\G^*$ and $((U^\omega_\hbar\G)^*)^\vee$ are also canonically 
isomorphic. 
\end{prop}

{\em Proof.} We must 
construct a Hopf pairing between $U^\omega_\hbar\G$ and $\cO^\omega_\hbar(G)$. 
The  Hopf pairing $T^\omega_\hbar(\G)\times \Sh^\omega(\G^*) \to \KK[[\hbar]][\hbar^{-1}]$ 
induces a pairing 
$T^\omega_\hbar(\G)\times \Sh^\omega_\hbar(\G^*) \to \KK[[\hbar]]$. This pairing  
restricts to a pairing   
 $T^\omega_\hbar(\G)\times \big( \Imm(\ell_{\G^*})\cap \Sh^\omega_\hbar(\G^*) \big) \to \KK[[\hbar]]$; 
now for any $x\in T^\omega_\hbar(\G)$ and $y\in T^\omega_\hbar(\G^*)$, 
$$
\langle \ell_\G(x),y\rangle_{\Sh^\omega(\G)\times T^\omega_\hbar(\G^*)} 
= \langle y\otimes x, \cR\rangle_{(T^\omega_\hbar(\G^*)\times \Sh^\omega(\G))
\otimes (T^\omega_\hbar(\G)\times \Sh^\omega(\G^*))}
= \langle \ell_{\G^*}(y),x\rangle_{\Sh^\omega(\G^*) \times T^\omega_\hbar(\G)},
$$ 
so the latter pairing descends to a Hopf pairing   
$$
\big( T^\omega_\hbar(\G) / \Ker(\ell_\G) \big) 
\times \big( \Imm(\ell_{\G^*})\cap \Sh^\omega_\hbar(\G^*) \big) \to \KK[[\hbar]],
$$ 
which is the desired pairing $U^\omega_\hbar\G\times \cO^\omega_\hbar(G) \to\KK[[\hbar]]$. 
The second part of the Proposition follows from the results of 
 \cite{Gav}. 
\hfill \qed \medskip

\subsection{Behaviour of $\G\mapsto U_\hbar^\omega(\G)$ for the 
double operation}

If $\G$ is a Lie bialgebra, let us denote by $D(\G)$ is 
double bialgebra. If $U$ is a QUE algebra, its quantum 
double $D(U)$ is the $\KK[[\hbar]]$-module $U\otimes (U^*)^\vee$. 
It is a quasitriangular QUE algebra.

\begin{prop} \label{prop:double}
If $\G$ is a Lie bialgebra, then there is an isomorphism
$\iota_\G : U_\hbar^\omega(D(\G)) \to D(U_\hbar^\omega(\G))$. 
Let $\cR_{\G,can}$ be the canonical $R$-matrix of  $D(U_\hbar^\omega(\G))$.
Then $(\iota_\G\otimes \iota_\G)^{-1}(\cR_{\G,can})$ belongs to 
$(i_{\G,D(\G)})^U(U^\omega_\hbar\G)
\widehat\otimes (i_{\G^*,D(\G)})^U(U_\hbar^\omega(\G^*))$
($i_{\G,D(\G)}$  and $i_{\G^*,D(\G)}$ 
are the inclusions of $\G$ and $\G^*$ in $D(\G)$).    
We have $(\iota_{\G^*}\otimes\iota_{\G^*})^{-1}(\cR_{\G^*,can}) 
= (\iota_\G\otimes\iota_\G)^{-1}(\cR_{\G,can})^{(21)}$, and 
if $\phi : \G\to\HH$ is a Lie algebra morphism, then we have 
$(\phi^U\otimes \id)(\iota_\G\otimes\iota_\G)^{-1}(\cR_{\G,can}) 
= (\id\otimes (\phi^*)^U)(\iota_\HH\otimes\iota_\HH)^{-1}(\cR_{\HH,can})$. 
\end{prop}

{\em Proof.} The maps $(i_{\G,D(\G)})^U$ and $(i_{\G^*,D(\G)})^U$
are flat deformations of the inclusions $U\G\to U(D(\G))$   and 
$U\G^*\to U(D(\G))$, so the composition   of their tensor 
product with the multiplication map defines a linear 
isomorphism from $U_\hbar^\omega\G\widehat\otimes U_\hbar^\omega\G^*$
to $U_\hbar^\omega D(\G)$. Moreover, 
$(i_{\G,D(\G)})^U$ and $(i_{\G^*,D(\G)})^U$ are also 
Hopf algebra morphisms. 

Recall that we have defined a solution $\cR_\G$ in 
$$
\Sh^\omega(\G)\widehat\otimes\Sh^\omega(\G^*)\subset \Sh^\omega(D(\G))
\widehat\otimes \Sh^\omega(D(\G))
$$ 
of QYBE; $U_\hbar^\omega\G$ and
$U_\hbar^\omega\G^*$ are the Hopf subalgebras of $\Sh^\omega(D(\G))$
defined by $\cR_\G$. These are Hopf subalgebras of 
$U_\hbar^\omega(D(\G))$, therefore $\cR_\G$ belongs to 
$(i_{\G,D(\G)})^U(U_\hbar^\omega\G)\widehat\otimes
(i_{\G^*,D(\G)})^U(U_\hbar^\omega\G^*)$.

If $y$ belongs to $T_\hbar^\omega\G$, and $x = \ell_{\G}(y)$, 
then $x = \langle \cR_\G, 
\id\otimes y\rangle_{T_\hbar^\omega\G\times\Sh^\omega(\G)}$, and 
$\Delta_{U_\hbar^\omega\G}(x) = 
\langle \cR_\G^{(13)} \cR_\G^{(23)} , 
\id \otimes \id\otimes y\rangle_{T_\hbar^\omega\G\times\Sh^\omega(\G)}$. 
Since $\cR_\G$ satisfies QYBE, and $U_\hbar^\omega\G$ is the image of 
$\ell_\G$, we have 
$$
\cR_\G((i_{\G,D(\G)})^U)^{\otimes 2}
(\Delta_{U_\hbar^\omega(\G)}(x)) = 
((i_{\G,D(\G)})^U)^{\otimes 2}(\Delta'_{U_\hbar^\omega(\G)}(x)) 
\cR_\G,
$$ 
for any $x$ in 
$U_\hbar^\omega\G$. In the same way, one shows that 
$\cR_\G((i_{\G^*,D(\G)})^U)^{\otimes 2}
(\Delta_{U_\hbar^\omega(\G^*)}(x)) = 
((i_{\G^*,D(\G)})^U)^{\otimes 2}(\Delta'_{U_\hbar^\omega(\G^*)}(x)) 
\cR_\G$, for any $x$ in 
$U_\hbar^\omega\G^*$. This proves that 
$\cR_\G\Delta_{U_\hbar^\omega(D(\G))} 
= \Delta'_{U_\hbar^\omega(D(\G))}\cR_\G$. 
So $(U_\hbar^\omega D(\G),\cR_\G)$ satisfies the
axioms of the double QUE algebra of $U_\hbar^\omega\G$ and
may therefore be identified with the double of  $U_\hbar^\omega\G$. 

The relation between $\cR_{\G,can}$ and $\cR_{\G^*,can}$ 
follows from Remark \ref{rem:dual:R} 
and the functoriality of 
$(\iota_\G\otimes\iota_\G^{-1})(\cR_{\G,can})$
follows from that of $\cR_\G$. 
\hfill \qed\medskip 

\begin{cor} \label{cor:compatible}
For any $\omega\in\SSh(\KK)$, the functor $\G\mapsto U_\hbar^\omega\G$
is a compatible quantization functor (see the Introduction).  
\end{cor}

{\em Proof.} This follows from Propositions \ref{rem:dual} 
and \ref{prop:double}. \hfill \qed\medskip

\newpage
\section{Proof of Theorem \ref{thm:last}}
\label{sect:last}

\subsection{Identification of two quantizations of $F(\G^*)$}
\label{natural:isom}

Let $\G$ be a finite-dimensional Lie algebra. Then $\G^*$
is a Lie coalgebra. Let $F(\G^*)$ be the free Lie algebra
of $\G^*$ (here $\G^*$ is viewed as a vector space). 
Then the map $\delta_{\G^*} : \G^*\to\G^*\otimes\G^*$ 
dual to the bracket of 
$\G$  extends to a unique cocycle map $\delta_{F(\G^*)} : 
F(\G^*)\to F(\G^*)^{\otimes 2}$, which endows  
$F(\G^*)$ with a Lie bialgebra structure. The corresponding 
Hopf-co-Poisson algebra is $(T(\G^*),\delta_{T(\G^*)})$.  
An element $\omega$ of $\cB(\KK[[\hbar]])$ being fixed, we 
now have two quantizations of this Hopf-co-Poisson 
algebra, namely $U_\hbar^{\omega}(F(\G^*))$ and 
$T_\hbar^\omega(\G^*)$. 

\begin{prop} \label{can:isom}
$U_\hbar^{\omega}(F(\G^*))$ and 
$T_\hbar^\omega(\G^*)$ are
canonically isomorphic. 
\end{prop}

{\em Proof.} The inclusion of $\G^*$ in $F(\G^*)$ as its
part of degree $1$ is a morphism of Lie coalgebras. 
This morphism  induces a morphism of  Hopf algebras 
$i : T^\omega_\hbar(\G^*)\to T^\omega_\hbar(F(\G^*))$. 
Let us compose it with the projection $\ell_{F(\G^*)}
: T_\hbar^\omega(F(\G^*)) \to U_\hbar^\omega(F(\G^*))$. 
(Recall that $U_\hbar^\omega(F(\G^*))$ is the image
of $\ell_{F(\G^*)}$, a Hopf subalgebra of 
$\Sh^\omega(F(\G^*))$.) Then $\ell_{F(\G^*)}\circ i$
is a morphism of Hopf algebras. Moreover, for any $\xi$ 
in $\G^*$, $i(\xi)$ is an element of $T_\hbar^\omega(F(\G^*))$
in $\xi + o(\hbar)$, and the image of the latter 
element in $U_\hbar^\omega(F(\G^*))$ is again in 
$\xi + o(\hbar)$. On the other hand, 
$U_\hbar^\omega(F(\G^*))$ is a deformation of $T(\G^*)$. 
The reduction mod $\hbar$ of the image of 
$\ell_{F(\G^*)}\circ i$ contains $\G^*$, so  
$\ell_{F(\G^*)}\circ i$ is surjective. The reduction 
mod $\hbar$ of this map is the identity, so it is also 
a linear isomorphism. So we have shown that 
$\ell_{F(\G^*)}\circ i$ is a Hopf algebra isomorphism 
from $T_\hbar^\omega(\G^*)$ to $U_\hbar^\omega(\G^*)$. 
This isomorphism is clearly functorial. 
\hfill \qed\medskip 

\subsection{Proof of Theorem \ref{thm:last}}

Let us now prove Theorem \ref{thm:last}. It follows from 
Corollary \ref{cor:compatible} that the map 
$\gamma_\KK$ is a bijection from $\SSh(\KK)$ to 
$\{$universal quantization functors of the tensor algebras $T(\A)\}$. 
In Theorem \ref{quantization:QLBA} and 
Proposition \ref{prop:functor}, we constructed a map 
$\al_\KK$ from $\SSh(\KK)$ to 
$\{$universal quantization functors of the Lie bialgebras$\}$. 

Any universal quantization functor of Lie bialgebras may be restricted to 
the category of Lie bialgebras of the form $F(\A)$, $\A$ a Lie coalgebra, 
and yields
therefore a universal quantization functor of the tensor algebras 
$T(\A)$. Let us denote by $\beta'_\KK$ the corresponding map from 
$\{$universal quantization functors of the Lie bialgebras$\}$
to $\{$universal quantization functors of the tensor algebras 
$T(\A)\}$. Then the map $\beta_\KK$ defined in the Introduction 
is $\gamma^{-1}_\KK \circ \beta'_\KK$. 
Proposition \ref{can:isom} implies that 
$\beta'_\KK \circ  \al_\KK = \gamma_\KK$. Therefore 
$\beta_\KK \circ  \al_\KK = \id_{\SSh(\KK)}$.  This proves the first part
of Theorem \ref{thm:last}. 

Let us prove the second part of this Theorem. It follows from 
Proposition \ref{rem:dual} and Proposition \ref{prop:double}
that the image of $\al_\KK$ consists of quantization functors
of Lie bialgebras that are compatible with the duals and doubles. 
Conversely, let us assume that $Q$ is a quantization functor 
compatible with the duals and doubles, and let $Q_0$ be the
restriction of $Q$ to the Lie bialgebras of the form $F(\A)$
(so $Q_0 = \beta'_\KK(Q)$). 
If $\A$ is any Lie bialgebra, then the unique extension to 
$F(\A)$ of the identity map of the vector space $\A$ to the Lie 
algebra $\A$ induces a Lie bialgebra morphism $F(\A)\to\A$. In 
the same way, we have a Lie bialgebra morphism $\A\to F(\A^*)^*$, 
and $\A$ may be characterized as the image of the composed morphism 
$F(\A)\to F(\A^*)^*$. So $Q(\A)$ may be characterized as the image of the 
morphism $Q(F(\A))\to Q(F(\A^*)^*)$. Since $Q$ is compatible with 
duals and restricts to $Q_0$, this is a Hopf algebra morphism  
from $Q_0(F(\A))$ to $Q_0(F(\A^*))^\vee$. Let us set $\omega = 
\gamma_\KK^{-1}(Q_0)$, then this Hopf algebra morphism 
is the same as an element $R(\A)$ of $\Sh^\omega(\A)\widehat\otimes 
\Sh^\omega(\A^*)$, satisfying the rules (\ref{ids:Delta}). 
Therefore $R(\A)$ has the form $\cR^\omega(\sigma_\A)$, with 
$\sigma_\A\in \hbar(\A\otimes\A^*)[[\hbar]]$. 
Let us show that $R(\A)$ is the same as the image of the 
canonical $R$-matrix $R_{can}(\A)$ of $Q(D(\A))$ by the 
tensor product of 
the injections $Q(\A)\to Q_0(F(\A^*))^\vee$ and 
$Q(\A^*)\to Q_0(F(\A))^\vee$.

By the properties of the quantum double, 
the identity map $Q(\A)\to Q(\A)$ may be identified with 
the linear map $Q(\A)\to Q(\A^*)^\vee$
defined by $x\mapsto (\id\otimes x)(R_{can}(\A))$. Since 
the canonical maps $Q(F\A)\to Q(\A)$ and
$Q(\A^*)^\vee\to Q(F\A^*)^\vee$ are respectively 
surjective and injective, 
the canonical map  
$Q(F\A)\to Q(F\A^*)^\vee$ may be defined by 
$x'\mapsto (\id\otimes x')(R'_{can}(\A))$, where 
$R'_{can}(\A)$ is the image of $R_{can}(\A)$ by the tensor
product of the injections $Q(\A)\to Q_0(F(\A^*))^\vee$ and 
$Q(\A^*)\to Q_0(F(\A))^\vee$. So $R'_{can}(\A) = R(\A)$. 

One can check that $R(\A)$ has the functoriality and 
duality properties $(\Sh^\omega(\phi)\otimes\id)(R(\A)) 
= (\id\otimes\Sh^\omega(\phi^*))(R(\B))$ and  
$R(\A^*) = R(\A)^{(21)}$. It then follows from Remark 
\ref{rem:dual} that the map $\A\mapsto 
\sigma_\A$ also has functoriality and duality properties.

Since we identified $R(\A)$ with the image of the $R$-matrix of 
the double $D(Q(\A))$ of $Q(\A)$, and since $D(Q(\A))$  identifies with 
$Q(D(\A))$ and injects into $Q((D\A)^*)^\vee$, $R(\A)$ satisfies
QYBE in the latter algebra. In Proposition 
\ref{thm:unique:sol}, we defined a series $\rho^\omega = 
\sum_{n\geq 1} \rho^\omega_n$, such that if $r\in\hbar(\G\otimes\G)
[[\hbar]]$ is a solution of CYBE, then $\rho^\omega(r)$ is a
solution of the Lie QYBE. The map $r\mapsto\rho^\omega(r)$ is bijective 
and one may show that it sets up a bijection between solutions of 
CYBE and of Lie QYBE. It follows that $(\rho^\omega)^{-1}(\sigma_\A)$
is a solution of CYBE. Let $\tau_\A$ be the endomorphism of $\A[[\hbar]]$, 
such that $(\rho^\omega)^{-1}(\sigma_\A) = (\tau_\A\otimes \id)(r_\A)$. 
Since $(\rho^\omega)^{-1}(\sigma_\A)$ is  
expressed polynomially in terms of the structure constants of $\A$, 
$\tau_\A$ is a obtained by composition of tensor products of
the bracket and cobracket map of $\A$. 

\begin{lemma} \label{lemma:G0}
Let $\cH_0$ be the set of functorial assignments $(\A\mapsto
\tau_\A)$, where $\A$ runs over all finite-dimensional Lie bialgebras 
and for any $\A$, 
$\tau_\A$ belongs to $(\A\otimes\A^*)[[\hbar]]$, such that:  
1) for any Lie algebra $\A$, $\rho_\A$  is a solution of CYBE (in 
$D(\A)^{\otimes 3}[[\hbar]]$), 
2) $\rho_\A$ satisfies $\tau_{\A^*} = \tau_\A^{(21)}$, 
3) $\rho_\A$ is equal to $r_\A$ modulo $\hbar$
and is expressed polynomially in terms of the structure constants of 
$\A$. Then the rule $(\A\mapsto \rho_\A)\mapsto (\A\mapsto 
(\rho_\A\otimes\id)(r_\A))$ defines a bijection from $\cG_0$ to 
$\cH_0$ ($\cG_0$ has been defined in the Introduction).  
\end{lemma}

{\em Proof.} Let us set $r_\A = \sum_{i\in I} a_i\otimes b_i$. 
Assume that $\tau_\A\in (\A\otimes\A^*)[[\hbar]]$ is a solution of 
CYBE. Let  $\rho_\A$ be the endomorphism of $\A[[\hbar]]$ such that
$\tau_\A = (\rho_\A\otimes id)(r_\A)$. Then we have also 
$\tau_\A = (\id\otimes\rho_\A^t)(r_\A)$, so $\rho_{\A^*} = \rho_\A^t$.
We have
\begin{align*}
& \lbrack\!\lbrack \tau_\A,\tau_\A \rbrack\!\rbrack    = 
\\ &
= \sum_{i,j\in I} [\rho_\A(a_i),a_j]\otimes b_i\otimes \rho_\A^t(b_j)
+ \rho(a_i) \otimes [b_i,a_j]\otimes \rho_\A^t(b_j)
+ \rho(a_i)\otimes a_j\otimes [\rho_\A^t(b_i),b_i] 
\\ & 
= \sum_{i,j\in I} [\rho_\A(a_i),a_j]\otimes b_i\otimes \rho_\A^t(b_j)
+ \rho_\A([a_j,a_i]) \otimes b_i\otimes \rho_\A^t(b_j)
\\ & + \rho_\A(a_i) \otimes a_j\otimes \rho_\A^t([b_j,b_i])
+ \rho_\A(a_i)\otimes a_j\otimes [\rho_\A^t(b_i),b_i] 
\end{align*}
so $\tau_\A$ satisfies the CYBE iff $\rho_\A$ is such that 
$\rho_\A([x,y]) = [\rho_\A(x),y]$ for any pair $x,y$ of elements 
of $\A$ and $\delta_\A(\rho_\A(x)) = (\id\otimes\rho_\A)(\delta_\A(x))$
for any element of $\A$. The second condition is satisfied because 
$\rho_{\A^*}([\xi,\eta]) = [\rho_{\A^*}(\xi),\eta]$ 
for any pair $\xi,\eta$ of elements of $\A^*$. 
\hfill \qed\medskip 

{\em End of proof of Theorem \ref{thm:last}.}
It follows from this Lemma that $(\A\mapsto\tau_\A)$ is an 
element of $\cG_0$. Let us define $\ve'_\KK$ as the unique map 
from $\{$compatible
quantization functors of Lie bialgebras$\}$ to $\SSh(\KK)\times 
\cG_0$, such that $\ve'_\KK(Q) = (Q_0,\tau)$. It is clear how to 
reconstruct $Q$ from $(Q_0,\tau)$, and that when $\tau$ is the neutral 
element of $\cG_0$, this reconstruction coincides with the 
map $Q_0 \mapsto \al_\KK(Q_0)$. $\ve'_\KK$ is therefore 
bijective, and we set $\ve_\KK = \ve_\KK^{\prime -1}$. 
\hfill \qed\medskip

\begin{remark} \label{rem:G0}
If $(\A\mapsto \rho_\A)$ belong to $\cG_0$, and $(\A,[\ ,\, ],\delta_\A)$
is any Lie bialgebra, then 
$(\A,[\ ,\, ],(\rho_\A\otimes\rho_\A)\circ\delta_\A\circ\rho_\A^{-1})$
is again a Lie bialgebra. Indeed, set $r' = \sum_i \rho_\A(a_i)\otimes b_i$. 
Since $r'$ is a solution of CYBE in $D(\A)$, and $r'$ belongs
$(\A\otimes\A^*)[[\hbar]]$, if we set $\delta'(x) = [r',x\otimes 1 + 1\otimes x]$, 
then $(\A,[\ ,\ ]_\A,\delta')$ is a Lie bialgebra. On the other hand,
the rule $\rho_\A([x,y]) = [\rho_\A(x),y]$ implies that  
$\delta'(x) = (\rho_\A\otimes \id)(\delta_\A(x))$  and the rules 
$\rho_{\A^*}([\xi,\eta]) = [\rho_{\A^*}(\xi),\eta]$ and $\rho_{\A^*}
= \rho_\A^t$ imply that 
$(\rho_\A\otimes \id) \circ \delta_\A = \delta_\A \circ \rho_\A$.
Since the image of $\delta_\A$ is antisymmetric, we also 
have 
$(\id\otimes \rho_\A) \circ \delta_\A = \delta_\A \circ \rho_\A$, 
so $\delta' = (\rho_\A\otimes\rho_\A)\circ \delta_\A\circ \rho_\A^{-1}$.
So $\rho_\A$ changes the Lie bialgebra structure of $\A$ 
preserving both its Lie algebra and Lie coalgebra structures. 
\end{remark}

% jeru 

\newpage

\begin{appendix} 

\section{Deformation of solutions of CYBE 
(proof of Prop.\ \ref{prop:leeds:first})}
\label{proof:problem:3} 

Let $A$ be an associative algebra, and let
$r_A\in A\otimes A$ be a solution of CYBE. 

It is natural to look for a sequence $(\cR_i)_{i\geq 0}$ of 
elements of $A\otimes A$, such that $\cR_0 = 1,\cR_1 = r_A$ and
$$
\forall N\geq 0, \quad 
\sum_{p,q,r\geq 0, p+q + r = N} 
\cR_p^{12} \cR_q^{13} \cR_r^{23} 
= \sum_{p,q,r\geq 0, p+q + r = N}
\cR_r^{23} \cR_q^{13}  \cR_p^{12} .  
$$
We call such a sequence $(\cR_i)_{i\geq 0}$ a quantization of the
CYBE solution $r_A$. If   
$(\cR_i)_{i\geq 0}$ is such a quantization, then for any formal 
parameter $t$, $\sum_{i\geq 0} t^i \cR_i$ is a solution of QYBE. 
Then if $(u_i)_{i\geq 0}$ is any sequence of elements of $A$ 
such that $u_0 = 1$, the sequence $(\cR'_i)_{i\geq 0}$ defined by 
$\sum_{i \geq 0} t^i \cR'_i = 
u (\sum_{i \geq 0} t^i \cR_i) u^{-1}$, 
where $u =  \sum_{i\geq 0} t^i u_i$, 
is also a quantization of $r_A$. We will say that the sequences
$(\cR_i)_{i\geq 0}$ and $(\cR'_i)_{i\geq 0}$ are equivalent.

One then seeks to solve inductively the above system of equations. 
More precisely, it can be written as follows 
\begin{equation} \label{recursion}
\lbrack\!\lbrack r_A,\cR_{N-1} \rbrack\!\rbrack = 
- \sum_{p,q,r \leq N-2, p+q+r = N} 
\cR_{p}^{12}\cR_{q}^{13}\cR_{r}^{23} 
- \cR_{r}^{23}  \cR_{q}^{13} \cR_{p}^{12} , 
\end{equation}
(we will call this equation the equation of order $N$)
where we set 
$$
\lbrack\!\lbrack r_A,R \rbrack\!\rbrack = 
[r_A^{12},R^{13}] + [r_A^{12},R^{23}] + [r_A^{23},R^{23}] 
+ [R^{12},r_A^{13}] + [R^{12},r_A^{23}] + [R^{23},r_A^{23}] . 
$$
$R\mapsto \lbrack\!\lbrack r_A,R \rbrack\!\rbrack$ is therefore a 
linear map from $A^{\otimes 2}$
to $A^{\otimes 3}$ (the corresponding Lie algebraic map 
was called $\delta_{3,r}$ in Section \ref{sect:univ:Lie}). 

On the other hand, let $\rho\mapsto \delta(r_A|\rho)$ be the linear 
map defined by 
\begin{align*}
& \delta(r_A|\rho) = [r_A^{12} + r_A^{13} + r_A^{14}, \rho^{234}] 
- [-r_A^{12} + r_A^{23} + r_A^{24}, \rho^{134}]
\\ &
+  [-r_A^{13} - r_A^{23} + r_A^{34}, \rho^{124}] 
- [-r_A^{14} - r_A^{24} - r_A^{34}, \rho^{123}] .  
\end{align*}

Then it follows from the fact that $r_A$ is a solution of CYBE that 
$\delta(r_A|\cdot) \circ \lbrack\!\lbrack r_A,\cdot \rbrack\!\rbrack = 0$
(this 
equality appears in the Lie coalgebra cohomology complex of $A$, 
endowed with the cobracket $\kappa(a) = [r_A,a\otimes 1 + 
1 \otimes a]$). Therefore 
$$
\Imm \lbrack\!\lbrack r_A,\cdot \rbrack\!\rbrack \subset \Ker  \delta(r_A|\cdot). 
$$

We will show
\begin{thm} \label{def:theory}
Let us set 
$$
H^1(A,r_A) = \Ker \lbrack\!\lbrack r_A,\cdot \rbrack\!\rbrack / \Imm\kappa, 
\quad 
H^2(A,r_A) = \Ker \delta(r_A,|\cdot) / \Imm \lbrack\!\lbrack r_A,\cdot \rbrack\!\rbrack.
$$ 
Then the equivalence classes of quantizations of $r_A$ are obstructed by 
$H^2(A,r_A)$ and parametrized by $H^1(A,r_A)$. 
\end{thm}

In particular, for any pair $(A,r_A)$ such that the homology 
$H^2(A,r_A)$ vanishes,  the solution $r_A$ can be quantized, the
equivalence classes of  its possible quantizations being parametrized by
$\prod_{i\geq 2} H^1(A,r_A)$.  

{\em Proof of Theorem \ref{def:theory}.}
We will show that if $\cR_1,\ldots,\cR_{N-2}$ satisfy equations
(\ref{recursion}) at all orders $\leq N-1$, then the right side
of the equation (\ref{recursion}) at order $N$
is contained in $\Ker \delta(r_A|\cdot)$ (Theorem \ref{prop:leeds}).

We first prove

\begin{lemma} \label{aryeh}
For $p$ integer $\geq 0$, let   
$(R,\rho)\mapsto \delta_p(R,\rho)$ be the map from 
$A^{\otimes 2}\times A^{\otimes 3}$ to $A^{\otimes 4}$, defined by 
$\delta_0 (R,\rho)= 0$, $\delta_p(R,\rho) = 0$ if $p\geq 4$ and 
\begin{align*}
& \delta_1(R,\rho) = [R^{12} + R^{13} + R^{14}, \rho^{234}] 
- [-R^{12} + R^{23} + R^{24}, \rho^{134}]
\\ &
+  [-R^{13} - R^{23} + R^{34}, \rho^{124}] 
- [-R^{14} - R^{24} - R^{34}, \rho^{123}] , 
\end{align*}
\begin{align*}
& \delta_{2}(R,\rho) = 
(R^{12}R^{13} + R^{12}R^{14} +  R^{13}R^{14})  \rho^{234}  
- \rho^{234} (R^{14}R^{13} + R^{14}R^{12} +  R^{13}R^{12})  
\\ &
- R^{23} R^{24} \rho^{134}  - (R^{23} + R^{24})\rho^{134} R^{12} 
+ R^{12} \rho^{134} (R^{23} + R^{24})  + \rho^{134} R^{24} R^{23} 
\\ & 
- R^{23} R^{13} \rho^{124} - (R^{13} + R^{23})\rho^{124} R^{34} 
+ R^{34} \rho^{124} (R^{13} + R^{23}) + \rho^{124} R^{13} R^{23}  
\\ & 
+ (R^{34}R^{24} + R^{34}R^{14} +  R^{24}R^{14})  \rho^{123} 
-\rho^{123} (R^{14}R^{24} + R^{14}R^{34} +  R^{24}R^{34}) , 
\end{align*}
and
\begin{align*}
& \delta_3(R,\rho) = R^{12}R^{13}R^{14} \rho^{234}  
- \rho^{234}  R^{14}R^{13}R^{12}
\\ &  
- R^{23}R^{24}\rho^{134} R^{12}  
+ R^{12}  \rho^{134} R^{24} R^{23}
\\ & 
- R^{23}R^{13}\rho^{124} R^{34}  
+ R^{34}  \rho^{124} R^{13} R^{23}
\\ &
+ R^{34}R^{24}R^{14} \rho^{123} 
-  \rho^{123} R^{14}R^{24}R^{34} .  
\end{align*}
Each $\delta_p$ is linear in $\rho$ and 
homogeneous of degree $p$ in $R$. 
Moreover, these maps satisfy the identity 
$$
\delta_{p}(R,R^{12}R^{13}R^{23} - R^{23}R^{13}R^{12})
+ \delta_{p+1}(R,[R^{12},R^{13}]  + [R^{12},R^{23}]  + [R^{13},R^{23}])
 = 0
$$
for any integer $p\geq 0$. 
\end{lemma}

{\em Proof of Lemma \ref{aryeh}.} This is a direct computation: for example, 
the identity for $p = 0$ follows from the Jacobi identity. 
\hfill \qed\medskip 

As we said, the proof of Theorem \ref{def:theory} reduces
to the following statement. 

\begin{prop} \label{prop:leeds}
Let $A$ be an algebra and $r_A$ belong to $A\otimes A$.  
Let $N$ be an integer $\geq 3$ and assume that 
$\cR_1,\ldots,\cR_{N-2}$ in $A\otimes A$ satisfy $\cR_1 = r_A$ and 
\begin{equation} \label{QYBE:low}
\lbrack\!\lbrack r_A,\cR_i \rbrack\!\rbrack  
= - \sum_{p+q+r = i+1, p,q,r>0} \cR_p^{12}\cR_q^{13}\cR_r^{23} 
+  \sum_{p+q+r = i+1,p,q,r>0} \cR_r^{23} \cR_q^{13} \cR_p^{12} ,  
\end{equation}
for $i = 1,\ldots,N-2$ (in particular, $r_A$ is a solution of 
CYBE). Then 
\begin{equation} \label{test:term}
\delta(r_A|
\sum_{p+q+r = N,p,q,r>0} \cR_p^{12}\cR_q^{13}\cR_r^{23} 
-  \cR_r^{23} \cR_q^{13} \cR_p^{12})
\end{equation}
is zero. 
\end{prop}

{\em Proof of Proposition \ref{prop:leeds}.}
Let us define $\wt\delta_{\al}(R_1,\ldots,R_\al|\rho)$  
as the coefficient of $t_1\cdots t_\al$ in 
${1\over{\al!}}\delta_\al(\sum_{\beta = 1}^\al 
t_\beta R_\beta,\rho)$. 
Then $\wt\delta_{\al}(R_1,\ldots,R_\al|\rho)$ is the unique 
multilinear form in $(R_1,\ldots,R_\al)$, symmetric 
in these arguments, such that 
$\wt\delta_\al(R,\ldots,R|\rho) = \delta_\al(R|\rho)$. 

Let us define  
$C(k,N,\al)$ as the subset of 
$\{1,\ldots,k\}^{\al} \times \{0,\ldots,k\}^3$ of $(\al+3)$uples 
$(k_1,\ldots,k_\al,p,q,r)$
such that $k_1 + \ldots + k_\al + p + q + r = N$ 
and set 
$$
\cT_k = \sum_{\al\geq 0} \sum_{(k_1,\ldots,k_\al,p,q,r)\in C(k,N,\al)}
\wt\delta_{\al}(\cR_{k_1},\ldots,\cR_{k_\al}|
\cR_p^{12}\cR_q^{13}\cR_r^{23} - \cR_r^{23}\cR_q^{13}\cR_p^{12}) . 
$$
Then the fact that the $\cR_i$ satisfy equations (\ref{QYBE:low}) implies that  
$\cT_{N-3}$ is equal to  
$$ 
\delta(r_A|\sum_{p+q+r = N-1, p,q,r\leq N-3}
\cR_p^{12}\cR_q^{13}\cR_r^{23} - \cR_r^{23}\cR_q^{13}\cR_p^{12}) ,  
$$
which is (\ref{test:term}). 

For $i = 0,\ldots,3$, denote by $C_i(k,N,\al)$ the subset of $C(k,N,\al)$
of $(\al+3)$uples $(k_1,\ldots,k_\al,p,q,r)$ such that exactly $i$ of the 
integers $p,q,r$ is equal to zero. Then $C(k,N,\al)$ is the disjoint
union of the $C_i(k,N,\al)$. Set 
$$
\cT_{k,i} = \sum_{\al\geq 0} \sum_{(k_1,\ldots,k_\al,p,q,r)\in C_i(k,N,\al)}
\wt\delta_{\al}(\cR_{k_1},\ldots,\cR_{k_\al}|
\cR_p^{12}\cR_q^{13}\cR_r^{23} - \cR_r^{23}\cR_q^{13}\cR_p^{12}) ,  
$$
then we have $\cT_k = \sum_{i=0}^3 \cT_{k,i}$. 
Since $\cT_{k,2} = \cT_{k,3} = 0$, we have 
\begin{equation} \label{eq:1}
\cT_k = \cT_{k,0} + \cT_{k,1}.
\end{equation}
Now 
\begin{align} \label{expr:1}
& \cT_{k,0} = \sum_{\al\geq 0} \sum_{(k_1,\ldots,k_{\al+3})\in \{1,\ldots,k\}^{\al + 3}
| \sum_\beta k_\beta = N}
\\ & \nonumber
\wt\delta_{\al}(\cR_{k_1},\ldots,\cR_{k_\al}|
\cR_{k_{\al+1}}^{12}\cR_{k_{\al+2}}^{13}\cR_{k_{\al+3}}^{23} 
- \cR_{k_{\al+3}}^{23}\cR_{k_{\al+2}}^{13}\cR_{k_{\al+1}}^{12}) , 
\end{align}
and  
\begin{align} \label{expr:2}
& \cT_{k,1} = \sum_{\al\geq 0} \sum_{(k_1,\ldots,k_{\al+3})\in \{1,\ldots,k\}^{\al + 3}
| \sum_\beta k_\beta = N}
\\ & \nonumber
\wt\delta_{\al+1}(\cR_{k_1},\ldots,\cR_{k_{\al+1}}|
[\cR_{k_{\al+2}}^{12} , \cR_{k_{\al+3}}^{13}] 
+ [\cR_{k_{\al+2}}^{12} , \cR_{k_{\al+3}}^{23}]
+ [\cR_{k_{\al+2}}^{13} , \cR_{k_{\al+3}}^{23}]) , 
\end{align} 
where we used the fact that $\delta_0 = 0$ to change $\al$ into $\al +1$
and where the first (resp., second, third) bracket corresponds to the subset of 
$C_1(N,k,\al)$ defined by the conditions $r = 0, p\neq 0,q\neq 0$
(resp., $q = 0, p\neq 0,r\neq 0$ and $p = 0, q\neq 0,r\neq 0$).

Let us show that each $\cT_k$ is equal to zero.  
Let us define, for $\unu = (\nu_1,\cdots,\nu_k)$ a collection of integers $\geq 0$
such that $\sum_{i=1}^k =\nu_i = \al +3$,  
$C(k,N,\unu)$ as the subset of
$C(k,N,\al)$ 
of $(\al+3)$uples $(k_1,\ldots,k_{\al+3})$
of integers such that for any $i = 1,\ldots,k$, card$\{\beta|k_\beta = i\} = \nu_i$. 
Set 
\begin{align*}
& \cT_k^{\unu} = 
\sum_{(k_1,\ldots,k_{\al+3})\in C(k,N,\al,\unu)}
\Big(  
\wt\delta_{\al}(\cR_{k_1},\ldots,\cR_{k_\al}|
\cR_{k_{\al + 1}}^{12}\cR_{k_{\al + 2}}^{13}\cR_{k_{\al + 3}}^{23} 
- \cR_{k_{\al + 3}}^{23} \cR_{k_{\al + 2}}^{13}\cR_{k_{\al + 1}}^{12})
\\ & 
+ \wt\delta_{\al + 1}(\cR_{k_1},\ldots,\cR_{k_{\al+1}}|
[\cR_{k_{\al + 2}}^{12},\cR_{k_{\al + 3}}^{13}]
+ [\cR_{k_{\al + 2}}^{12},\cR_{k_{\al + 3}}^{23}] 
+ [\cR_{k_{\al + 2}}^{13},\cR_{k_{\al + 3}}^{23}]) \Big) . 
\end{align*}
Then it follows from (\ref{eq:1}) and the expressions (\ref{expr:1})   and
(\ref{expr:2}) for $\cT_{k,0}$    and $\cT_{k,1}$ that   
$\cT_k = \sum_{\unu\in\NN^k} \cT_k^{\unu}$. 
Define the sequence $(k^\unu_1,\ldots,k^\unu_{\al + 3})$ by the conditions 
that it is non-decreasing and that it belongs to  $C(k,N,\unu)$. 
Then $\cT_k^{\unu}$ is proportional to the symmetrization 
\begin{align} \label{symm}
& \sum_{\sigma\in \SG_{\al + 3}}
\wt\delta_{\al}(\cR_{k^\unu_{\sigma(1)}},\ldots,\cR_{k^\unu_{\sigma(\al)}}|
\cR_{k^\unu_{\sigma(\al + 1)}}^{12}\cR_{k^\unu_{\sigma(\al + 2)}}^{13}
\cR_{k^{\unu}_{\sigma(\al + 3)}}^{23} 
- \cR_{k^\unu_{\sigma(\al + 3)}}^{23} \cR_{k^\unu_{\sigma(\al + 2)}}^{13}\cR_{k^\unu_{\sigma(
\al + 1)}}^{12}) 
\\ & \nonumber 
+ \wt\delta_{\al+1}(
\cR_{k^\unu_{\sigma(1)}},\ldots,\cR_{k^\unu_{\sigma(\al+1)}}|
[\cR_{k^\unu_{\si(\al + 2)}}^{12},\cR_{k^\unu_{\si(\al + 3)}}^{13}]
+ [\cR_{k^\unu_{\si(\al + 2)}}^{12},\cR_{k^\unu_{\si(\al + 3)}}^{23}] 
+ [\cR_{k^\unu_{\si(\al + 2)}}^{13},\cR_{k^\unu_{\si(\al + 3)}}^{23}]) . 
\end{align}

We have now 
\begin{lemma}
If $\al$ is an integer and
$R_1,\ldots,R_{\al + 3}$ belong to $A\otimes A$, 
we have 
\begin{align*}
& 
\sum_{\sigma\in \SG_{\al + 3}}
\wt\delta_{\al}(R_{\sigma(1)},\ldots,R_{\sigma(\al)}|
R_{\sigma(\al + 1)}^{12}R_{\sigma(\al + 2)}^{13}
R_{\sigma(\al + 3)}^{23} 
- R_{\sigma(\al + 3)}^{23} R_{\sigma(\al + 2)}^{13}
R_{\sigma(\al + 1)}^{12}) 
\\ & 
+ \wt\delta_{\al + 1}( 
R_{\sigma(1)},\ldots,R_{\sigma(\al+1)}
| [R_{{\si(\al + 2)}}^{12},R_{{\si(\al + 3)}}^{13}]
+ [R_{{\si(\al + 2)}}^{12},R_{{\si(\al + 3)}}^{23}] 
+ [R_{{\si(\al + 2)}}^{13},R_{{\si(\al + 3)}}^{23}]) = 0. 
\end{align*}
\end{lemma}

{\em Proof of Lemma.} The left side is 
multilinear and symmetric in variables $R_1,\ldots,R_{\al +3}$; moreover, 
its value for $R_1 = \ldots = R_{\al +3} = R$ is equal to 
$$
(\al+3)! \big( \delta_\al(R, R^{12}  R^{13} R^{23}  - R^{23} R^{13} R^{12})   
+ \delta_{\al+1}(R,[R^{12},R^{13}] + [R^{12},R^{23}] + [R^{13},R^{23}] ) 
\big), 
$$
which is equal to zero. 
Therefore the left side is equal to zero. \hfill \qed\medskip 

{\em End of proof of Proposition \ref{prop:leeds}.}
Since $\cT_k^\unu$ is proportional to (\ref{symm}), the above Lemma  
implies that it is zero. It follows that for any $k$, $\cT_k$ is
equal to zero. In particular, $\cT_{N-3} = 0$, which proves the Proposition. 
\hfill \qed\medskip  

\begin{remark}
When $N = 3$, equation (\ref{recursion}) is written 
\begin{equation} \label{part:case}
\lbrack\!\lbrack r_A,\cR_2 \rbrack\!\rbrack = - r_A^{12}r_A^{13}r_A^{23} 
+  r_A^{23}r_A^{13} r_A^{12} . 
\end{equation}
One checks that the fact that $r_A$ is a solution of CYBE implies that 
$\cR_2 = {1\over 2} r_A^2$ is a solution of (\ref{part:case}). Then if $x$ is 
any element of $A$, then $\cR_2 = {1\over 2} r_A^2 + \kappa(x)$ is also a 
solution of (\ref{part:case}). 
 
In Theorem \ref{big:thm:QYBE}, we construct a solution $\cR(\rho(r_\A))$
of QYBE. In that case, $A = \Sh(\A)[[\hbar]]$, $r_A$ is the element $r_\A
\in\A\otimes\A$, viewed as an element of $A\otimes A$, and $\cR_2 = 
{1\over 2}r_\A^2 + \kappa(x)$, where $x = {1\over 2}\sum_{i\in I} (a_i b_i)$, 
and we write $r_\A = \sum_{i\in I} a_i \otimes b_i$. 
\end{remark}
\newpage 

\section{Construction of $\mu^{p,q,r}_{Lie}$ 
(proof of Prop.\ \ref{prop:key})} 
\label{proof:prop:key}

Let us denote by $FA_n$ the part of the
free algebra in $n$ generators, homogeneous of degree one
in each generator. In this Section, we also denote by $FL_n$ 
the subspace of $FA_n$ consisting of Lie elements (so $FL_n
= Free_n$). 

\subsection{Definition of the algebra $(F^{(3)},m_{F^{(3)}})$}

When $p,q,r$ are integers $\geq 0$, we set 
$$
F^{(3)}_{pqr} = \big( FA_{q+r} \otimes (FA_p \otimes 
FA_q) \otimes FA_{p+r} \big)_{\SG_p\times\SG_q \times\SG_r} 
\quad \on{and} \quad F^{(3)} = \oplus_{p,q,r\geq 0} F^{(3)}_{pqr} .  
$$
Let us denote by $x^{(2)}_1,\ldots,x^{(2)}_q,x^{(3)}_1,\ldots,x^{(3)}_r$ the generators
of $FA_{q+r}$, by $x^{(1)}_1,\ldots,x^{(1)}_p$ the generators of $FA_p$, by 
$y^{(2)}_1,\ldots,y^{(2)}_q$ the generators of $FA_q$, and by 
$y^{(1)}_1,\ldots,y^{(1)}_p$, $y^{(3)}_1,\ldots,y^{(3)}_r$ the generators of 
$FA_{p+r}$. Then the first and last set of generators are 
split in two subsets (e.g., the first subset of generators of 
$FA_{q+r}$ is $x^{(2)}_1,\ldots,x^{(2)}_q$). 
The symmetric group $\SG_p$  
(resp., $\SG_q$ and $\SG_r$) acts on 
$FA_{q+r} \otimes (FA_p \otimes 
FA_q) \otimes FA_{p+r}$
by simultaneously  permuting the 
variables $(x^{(1)}_i)_{i = 1,\ldots,p}$ and $(y^{(1)}_i)_{i = 1,\ldots,p}$
(resp., the variables $(x^{(2)}_i)_{i = 1,\ldots,q}$ and $(y^{(2)}_i)_{i = 1,\ldots,q}$ 
and the variables $(x^{(3)}_i)_{i = 1,\ldots,r}$ and $(y^{(3)}_i)_{i = 1,\ldots,r}$). 

If $n,n'$ are integers $\geq 0$, let us define $M_{n,n'}$
as the set of pairs of maps $(c,c')$, where $c$ is a map 
from $\{1,\ldots,n'\}$ to $\{-\infty,1,\ldots,n\}$ and 
$c'$ is a map from $\{1,\ldots,n\}$ to $\{1,\ldots,n',\infty\}$, 
such that for any $k\in \{1,\ldots,n\}$, and any 
$k'\in \{1,\ldots,n'\}$, the inequalities 
\begin{equation} \label{ineqs}
c(c'(k)) < k \quad \on{and}\quad  c'(c(k')) > k'
\end{equation}
hold whenever the left sides are defined. By convention, if $k$
is any integer, $k < \infty$ and $-\infty < k$. 

We are going to define a bilinear map $m_{F^{(3)}} : F^{(3)}
\otimes F^{(3)} \to F^{(3)}$. 
If $p,q$ and $r$ are integers $\geq 0$, 
there is a unique 
linear isomorphism $\xi\mapsto a_\xi$ from $FA_{q+r}\otimes FA_{p+r}$  
to $F^{(3)}_{pqr}$, such that if $P\in FA_{q+r},Q\in FA_{p+r}$ and $(\sigma,\tau)
\in \SG_{p}\times \SG_q$,  
then 
\begin{align*}
a_{P\otimes Q} = & P(x^{(2)}_1,\ldots,x^{(2)}_q,x^{(3)}_1,\ldots,x^{(3)}_r) 
\otimes x^{(1)}_{1}\ldots x^{(1)}_{p} 
y^{(2)}_{1}\ldots y^{(2)}_{q}
\\ & \otimes Q(y^{(1)}_1,\ldots,y^{(1)}_p,y^{(3)}_1,\ldots,y^{(3)}_r) . 
\end{align*}

If $p,q,r$, and $p',q',r'$ are integers $\geq 0$, and if 
$P\in FA_{q+r}$, $Q\in FA_{p+r}$, $P'\in FA_{q'+r'}$, 
$Q'\in FA_{p'+r'}$, 
let us set 
$$
a_{(P,Q),(P',Q')} = \sum_{(c,c')\in M_{q,p'}}
a_{(P,Q),(P',Q')}(c,c') , 
$$
where for any $(c,c') \in M_{q,p'}$, we set 
$\al = \card (c')^{-1}(\infty)$ and 
$\al' = \card c^{-1}(-\infty)$  and define 
$a_{(P,Q),(P',Q')}(c,c')$
as the element of $F^{(3)}_{p+\al',q'+\al,p'+q+r+r' - \al - \al'}$ 
given by 
\begin{align*}
&
a_{(P,Q),(P',Q')}(c,c') = 
\\ & 
\big( P(\prod^{\nearrow}_{j\in c^{-1}(1)}\ad'(x^{(1)}_{p+j})
(x^{(2)}_1), 
\ldots, 
\prod^{\nearrow}_{j\in c^{-1}(q)}\ad'(x^{(1)}_{p+j})(x^{(2)}_q), 
x^{(3)}_1,\ldots,x^{(3)}_r) 
\\ & 
P'(x^{(2)}_{q+1},\ldots,x^{(2)}_{q+q'},
x^{(3)}_{r+1},\ldots,x^{(3)}_{r+r'}) \big) 
\\ & 
\otimes  
\big( \prod^{\nearrow}_{i\in \{1,\ldots,p\}} x^{(1)}_{i}  
\prod^{\nearrow}_{i'\in c^{-1}(\infty)} x^{(1)}_{p+i'} 
\otimes 
\prod^{\nearrow}_{j\in (c')^{-1}(\infty)} y^{(2)}_{j} 
\prod^{\nearrow}_{j'\in \{1,\ldots,q'\}} y^{(2)}_{q+j'} 
\big) 
\\ & \otimes 
\big( Q(y^{(1)}_1,\ldots,y^{(1)}_p, y^{(3)}_1,\ldots,y^{(3)}_r) 
\\ & 
Q'(\prod^{\searrow}_{j\in (c')^{-1}(1)} 
\ad'(y^{(2)}_{j})(y^{(1)}_{p+1}), 
\ldots, 
\prod^{\searrow}_{j\in (c')^{-1}(p')} 
\ad'(y^{(2)}_{j})(y^{(1)}_{p+p'}), 
y^{(3)}_{r+1},\ldots,y^{(3)}_{r+r'}) \big) . 
\end{align*}

Let us explain the notation in this formula.  The generators of the
components of  $F^{(3)}_{p+\al',q'+\al,p'+q+r+r' - \al - \al'}$ are denoted as
follows: 

-- generators of $FA_{p'+q+q'+r+r'-\al'}$. 
First subset of generators: $x^{(2)}_i$, $i\in (c')^{-1}(\infty)$
and $x^{(2)}_{q+j'}$, $j'\in \{1,\ldots,q'\}$.  
Second subset of generators: $x^{(1)}_{p+i'}$, 
$i'\in c^{-1}(\{1,\ldots,q\})$, 
$x^{(2)}_j$, $j\in (c')^{-1}(\{1,\ldots,q'\})$
and $x^{(3)}_k$, $k\in \{1,\ldots,r+r'\}$

-- the generators of $FA_{p+\al'}$ are  
$x^{(1)}_{i}$, $i\in  \{1,\ldots,p\}$ 
and $x^{(1)}_{p+i'}$, $i'\in c^{-1}(-\infty)$ 

-- the generators of  $FA_{p'+\al}$ are  
$y^{(2)}_j$, $j\in
\{q+1,\ldots,q+q'\}$ and $y^{(2)}_{j'}$, $j'\in (c')^{-1}(\infty)$

-- generators of  $FA_{p+p'+q+r+r'-\al}$. 
First subset of generators:  
$y^{(1)}_i$, $i\in \{1,\ldots,p\}$, $y^{(2)}_{p+i'}$, 
$i'\in c^{-1}(-\infty)$. 
Second subset of generators: 
$y^{(1)}_{p+i'}$, $i'\in c^{-1}(\{1,\ldots,q\})$, 
$y^{(2)}_{j}$, $j\in (c')^{-1}(\{1,\ldots,p'\})$, 
$y^{(3)}_{k}$, $k\in \{1,\ldots,r+r'\}$. 

For $x,y$ elements of an associative algebra, we set
$\ad'(x)(y) = yx - xy$ (so $\ad' = - \ad$). If $J$
is a finite ordered set of indices, and $J = \{j_1,\ldots,j_r\}$, with 
$j_1<\ldots <j_r$, then $\prod^{\nearrow}_{j\in J} a_j$ 
and $\prod^{\searrow}_{j\in J} a_j$ denote the products of elements 
$a_{j_1} a_{j_2}\cdots a_{j_r}$ and $a_{j_r} a_{j_{r-1}}\cdots a_{j_1}$. 
So $\prod^{\nearrow}_{i\in\{1,\ldots,n\}}\ad'(a_i)(a) 
= [[[a,a_n],a_{n-1}],\ldots,a_1]$. 

Let us set $\wt F_{pqr} = FA_{q+r} \otimes 
(FA_p\otimes FA_q) \otimes FA_{p+r}$ and $\wt F = 
\oplus_{p,q,r\geq 0} \wt F_{pqr}$. Then the rule 
$(P\otimes Q)\otimes (P'\otimes Q') \mapsto a_{(P,Q),(p',Q')}$
defines a linear map from $(FA_{q+r}\otimes Fa_{p+r})
\otimes(FA_{q'+r'}\otimes FA_{p'+r'})$ to $\wt F$, which is covariant
with respect to the action of $\SG_r\times\SG_{r'}$. Therefore, 
its induces a linear map from $(FA_{q+r}\otimes Fa_{p+r})_{\SG_r}
\otimes (FA_{q'+r'}\otimes FA_{p'+r'})_{\SG_{r'}}$ 
(which we identified with $F^{(3)}_{pqr}\otimes F^{(3)}_{p'q'r'}$)
to $F^{(3)}$. This map may be extended by linearity in a unique way 
to  a linear  map $m_{F^{(3)}}$ from $F^{(3)}\otimes F^{(3)}$ 
to $F^{(3)}$.

\subsection{Associativity of $m_{F^{(3)}}$}

In this Section, we prove that $m_{F^{(3)}}$ is associative.  
For this, we first define composition operations on the sets of 
maps $M_{nn'}$ introduced above. 

\subsubsection{The operations $\Comp^{12,3}$ and $\Comp^{1,23}$}  

For any pair of integers $(\al,\al')$, define 
$M_{nn'}^{\al\al'}$  as the subset of $M_{n,n'}$ of all pairs 
$(c,c')$ such that $\card c^{-1}(-\infty) = \al'$ and $\card (c')^{-1}(\infty)
= \al$. 

For any quadruple of integers $(n,n',n'',n''')$, let us also define 
$M_{n,n'|n'',n'''}$ as the subset of $M_{n+n',n''+n'''}$
of all pairs $(\wt c,\wt c')$ such that 
$$
\wt c(\{1,\ldots,n'\}) \subset \{-\infty,1,\ldots,n\}
\quad \on{and} \quad 
\wt c'(\{n+1,\ldots,n+n''\}) \subset \{n'+1,\ldots,n'+n''',\infty\} . 
$$
We define then two maps 
$$
\Comp^{12,3} : \coprod_{(\al,\al')\in\NN^2}
( M_{nn'}^{\al\al'} \times M_{\al + n'',n'''})
\to M_{n,n''|n',n'''} 
$$
and 
$$
\Comp^{1,23} : \coprod_{(\al'',\al''')\in\NN^2}
( M_{n,n'+\al'''} \times M_{n''n'''}^{\al''\al'''}) 
\to M_{n,n''|n',n'''} 
$$
in the following way. If $(c,c')\in M_{nn'}^{\al\al'}$ and 
$(c'',c''')\in M_{\al+n'',n'''}$, then 
$$
\Comp^{12,3}((c,c'),(c'',c''')) = (\wt c,\wt c'),
$$ 
where $\wt c$ and $\wt c'$ are defined as follows. Let us denote 
the elements of $(c')^{-1}(\infty)$ by $i_1,\ldots,i_\al$, 
where the sequence $(i_\beta)_{\beta = 1,\ldots,\al}$ is increasing. 
Then if $i'\in \{1,\ldots,n'\}$, then  
$\wt c(i') = c(i')$; if $i'''\in \{1,\ldots,n'''\}$
and  $c''(i''')\in \{1,\ldots,\al\}$, then $\wt c(n' + i''') = i_{c''(i''')}$; and 
if $c''(i''')\in \{\al + 1,\ldots,\al + n'''\}$, then  
 $\wt c(n' + i''') = (n - \al ) + c''(i''')$. On the other hand, 
if $i\in \{1,\ldots,n\} - (c')^{-1}(\infty)$, then $\wt c'(i) = c'(i)$, 
and for $\beta\in \{1,\ldots,\al\}$, $\wt c'(i_\beta) = n' + c'''(\beta)$; 
and if $i''\in \{1,\ldots,n''\}$, then $\wt c'(n + i'') = n' + c'''(\al + i'')$. 

By convention, if $k$ is any integer, then $k + \infty = \infty$ and
$-\infty + k = - \infty$.  

In the same way, if $(d,d')\in M_{n,n'+\al''}$ and $(d'',d''')\in 
M_{n''n'''}^{\al''\al'''}$, then we set $\Comp^{1,23}((d,d'),(d'',d''')) = (\wt d,\wt d')$, 
where $\wt d$ and $\wt d'$ are defined as follows. Denote the elements of 
$(d'')^{-1}(-\infty)$ by $i_1,\ldots,i_{\al'''}$, where the sequence 
$(i_\beta)_{\beta = 1,\ldots,\al'''}$ is increasing. Then 
$\wt d(i') = d(i')$ if $i'\in\{1,\ldots,n'\}$, $\wt d(n' + i_\beta) = d(n' + \beta)$
if $\beta\in\{1,\ldots,\al'''\}$ and $\wt d(n' + i''') = n + d''(i''')$ if 
$i'''\in\{1,\ldots,n'''\}$ and $i'''\notin (d''')^{-1}(\infty)$. On the other hand, 
if $i\in \{1,\ldots,n\}$, then $\wt d'(i) = d'(i)$ 
if $d'(i)\in \{1,\ldots,n'\}$, and  
$\wt d'(i) = n' + i_{d'(i) - n'}$ if $d'(i)\in \{n' + 1,\ldots,n'+ \al'''\}$; and 
if $i''\in\{1,\ldots,n''\}$, then $\wt d'(n+i'') = n' + d'''(i'')$. 

\begin{lemma}
The maps $\Comp^{12,3}$ and $\Comp^{1,23}$ are both bijective. 
\end{lemma}

{\em Proof.} One first checks that the pairs $(\wt c,\wt c')$ and 
$(\wt d,\wt d')$ defined above actually belong to $M_{n,n''|n',n'''}$. 
This is a direct verification. One also checks that each pair 
$(\wt c,\wt c')$ has a unique preimage by $\Comp^{12,3}$ and $\Comp^{1,23}$. 
For example, let us describe $(\Comp^{12,3})^{-1}(\wt c,\wt c') = 
((c,c'),(c'',c'''))$. This 
is an element of $M_{nn'}^{\al\al'}\times M_{\al + n'',n'''}$, where 
$$
\al = \card \big( \wt c^{-1}(\{n'+1,\ldots,n'+n''',-\infty\}) \cap \{1,\ldots,n\} \big) 
$$
and $\al' = \card(\wt c^{-1}(-\infty)\cap \{1,\ldots,n'\})$. 
Let us denote by $j_1,\ldots,j_\al$ the elements of 
$c^{-1}(\{n'+1,\ldots,n'+n''',-\infty\}) \cap \{1,\ldots,n\}$, 
where the sequence $(j_\beta)_{\beta = 1,\ldots,\al}$ is increasing.

Then the pairs $(c,c')$ and $(c'',c''')$ are obtained as follows. 
If $j' \in \{1,\ldots,n'\}$, then $c(j') = \wt c(j')$. If 
$j'''\in \{1,\ldots,n'''\}$ and $\wt c(n' + j''')\in \{1,\ldots,n\}$, 
then since $\wt c'(\wt c(n'+j''')) > n'+j'''\geq n'+1$,  
$\wt c(n' + j''')$ belongs to $\{j_1,\ldots,j_\al\}$,  and 
we define $c''(j''')$ as the index $\beta$ such that $\wt c(n'+j''') = j_\beta$. 
If  $j'''\in \{1,\ldots,n'''\}$ and $\wt c(n' + j''')\in \{n+1,\ldots,n+n'',-\infty\}$,
then we set $c''(j''') = \wt c(n'+j''') + \al - n$. 

If $j\in \{1,\ldots,n\}$, we define $c'(j)$ as $\wt c'(j)$ if $\wt c'(j)\in 
\{1,\ldots,n'\}$ and as $\infty$ else. For any $\beta\in \{1,\ldots,\al\}$, 
we define $c'''(\beta)$ as $\wt c'(j_\beta) - n'$, and if $j''\in 
\{\al+1,\ldots,\al+n''\}$, we define $c'''(j'')$ as 
$\wt c'(j'' + n - \al) - n'$.   
\hfill \qed\medskip 

\subsubsection{Associativity of $m_{F^{(3)}}$}

\begin{thm}
$m_{F^{(3)}}$ is associative. 
\end{thm}

{\em Proof.}
Let $(p,q,r)$, $(p',q',r')$ and $(p'',q'',r'')$ be arbitrary 
triples of integers $\geq 0$. Let $P,P',P''$ be 
elements of $FA_{q+r},FA_{q'+r'},FA_{q''+r''}$, let  
$Q,Q',Q''$ be elements of  
$FA_{p+r},FA_{p'+r'},FA_{p''+r''}$. 
It will be enough to  prove that 
\begin{equation} \label{eq:assoc}
m_F(a_{(P\otimes Q} 
\otimes a_{(P',Q'),(P'',Q'')})  
=  m_F( a_{(P,Q),(P',Q')} \otimes 
a_{P''\otimes Q''})  . 
\end{equation}
The left side is a sum indexed by $\coprod_{(\al,\al')\in\NN^2} 
M_{qp'}^{\al\al'} \times M_{\al+q',p''}$ 
and the right side is a 
sum indexed by $\coprod_{(\al'',\al''')\in\NN^2} M_{q,p'+\al'''}
\times M_{q'p''}^{\al''\al'''}$. 
Using maps $\Comp^{12,3}$ and $\Comp^{1,23}$, 
we transform both sums into the same expression 
$$
\sum_{(\wt c,\wt c')\in M_{q,q'|p',p''}} 
a_{(P,Q),(P',Q'),(P'',Q'')}
(\wt c,\wt c'), 
$$ 
where we set $\beta = \card (\wt c')^{-1}(\infty)$, $\beta' = 
\card \wt c^{-1}(-\infty)$, 
and 
$$
a_{(P,Q),(P',Q'),(P'',Q'')}(\wt c,\wt c')
$$ 
is the element of
$F^{(3)}_{p+\beta,q''+\beta',p'+p''+q+q'+r+r'+r''-\beta - \beta'}$ 
given by 
\begin{align} \label{big:expr}
& a_{(P,Q),(P',Q'),
(P'',Q'')}(\wt c,\wt c') = 
\\ & \nonumber  
\big( 
P(\prod^{\nearrow}_{i\in \wt c^{-1}(1)} 
\ad'(x^{(1)}_{p+i})(x^{(2)}_1), \ldots,
\prod^{\nearrow}_{i\in \wt c^{-1}(q)} 
\ad'(x^{(1)}_{p+i})(x^{(2)}_{q}), 
x^{(3)}_1,\ldots,x^{(3)}_r)
\\ & \nonumber 
P'(\prod^{\nearrow}_{i\in \wt c^{-1}(q+1)} 
\ad'(x^{(1)}_{p+i})(x^{(2)}_{q+1}), \ldots,
\prod^{\nearrow}_{i\in \wt c^{-1}(q+q')} 
\ad'(x^{(1)}_{p+i})(x^{(2)}_{q+q'}), 
x^{(3)}_{r+1},\ldots,x^{(3)}_{r+r'})
\\ & \nonumber 
P''(x^{(2)}_{q+q'+1},\ldots,x^{(2)}_{q+q'+q''}, 
x^{(3)}_{r+r'+1},\ldots,x^{(3)}_{r+r'+r''}) \big) 
\\ & \nonumber 
\otimes \big( \prod^{\nearrow}_{i\in\{1,\ldots,p\}} x^{(1)}_{i}
\prod^{\nearrow}_{i'\in \wt c^{-1}( - \infty)} 
x^{(1)}_{p+i'}
\otimes \prod^{\nearrow}_{j\in (\wt c')^{-1}(\infty)} 
y^{(2)}_{j}
\prod^{\nearrow}_{j''\in\{1,\ldots,q''\}} y^{(2)}_{q+q'+j''}
\big) 
\\ & \nonumber 
\otimes \big( 
Q(y^{(1)}_1,\ldots,y^{(1)}_{p},y^{(3)}_1,\ldots,y^{(3)}_r) 
\\ & \nonumber 
Q'(\prod^{\searrow}_{j\in (\wt c')^{-1}(1)}
\ad'(y^{(2)}_{j}) (y^{(1)}_{p+1}),\ldots,
\prod^{\searrow}_{j\in (\wt c')^{-1}(p')}
\ad'(y^{(2)}_{j})(y^{(1)}_{p+p'}),
y^{(3)}_{r+1},\ldots,y^{(3)}_{r+r'}) 
\\ & \nonumber 
Q''( \prod^{\searrow}_{j\in (\wt c')^{-1}(p'+1)}
\ad'(y^{(2)}_{j})
(y^{(1)}_{p+p'+1}),\ldots,
\prod^{\searrow}_{j\in (\wt c')^{-1}(p'+p'')}
\\ & \nonumber 
\ad'(y^{(2)}_{j}) (y^{(1)}_{p+p'+p''}), 
y^{(3)}_{r+r'+1},\ldots,y^{(3)}_{r+r'+r''}) 
\big) . 
\end{align} 
In this formula, the generators of the factors of 
$F^{(3)}_{p+\beta,q''+\beta',p'+p''+q+q'+r+r'+r''-\beta - \beta'}$ 
are denoted as follows: 

-- generators of $F^{(3)}_{p'+p''+q+q'+q''+r+r'+r''-\beta}$. 
First subset of generators: 
$x^{(2)}_j$, $j\in (\wt c')^{-1}(\infty)$,  
$x^{(2)}_{q+q'+j''}$, $j''\in \{1,\ldots,q''\}$.  
Second subset of generators: 
$x^{(1)}_{p+j'}$, $j'\in \wt c^{-1}(\{1,\ldots,q+q'\})$, 
$x^{(2)}_j$, $j\in (\wt c')^{-1}(\{1,\ldots,p'+p''\})$, $x^{(3)}_k$, 
$k\in \{1,\ldots,r+r'+r''\}$. 
 
-- generators of $F^{(3)}_{p+\beta}$: $x^{(1)}_{i}$, $i\in \{1,\ldots,p\}$
and $x^{(1)}_{p+i'}$, $i'\in \wt c^{-1}(\infty)$

-- generators of $F^{(3)}_{q''+\beta'}$: 
$y^{(2)}_{j}$, $j\in (\wt c')^{-1}(\infty)$, 
and $y^{(2)}_{q+q'+j''}$, $j''\in \{1,\ldots,q''\}$

-- generators of $F^{(3)}_{p+p'+p''+q+q'+r+r'+r''-\beta'}$. 
First subset of generators: $y^{(1)}_i$, $i\in \{1,\ldots,p\}$, 
and $y^{(1)}_{p+i'}$, $i'\in \wt c^{-1}(\infty)$.  
Second subset of generators: 
$y^{(1)}_{p+i'}$, $i'\in \wt c^{-1}(\{1,\ldots,q+q'\})$, 
$y^{(2)}_{i}$, $i\in (\wt c')^{-1}(\{1,\ldots,p'+p''\})$, 
and $y^{(3)}_i$, $i\in \{1,\ldots,r+r'+r''\}$. 

Since each side of (\ref{eq:assoc}) 
is equal to (\ref{big:expr}), equation (\ref{eq:assoc}) is satisfied. 
\hfill \qed\medskip

\subsection{Universal properties of $(F^{(3)},m_{F^{(3)}})$}

If $p,q,r$ are integers $\geq 0$, a basis of $F^{(3)}_{pqr}$ consists of the 
$$
\prod_{k = 1}^{q+r} z_{\sigma(k)} \otimes 
(\prod_{i = 1}^p x^{(1)}_{i} \otimes \prod_{j = 1}^q 
y^{(2)}_{j})
\otimes \prod_{l = 1}^{p+r} t_{\tau(l)} , 
$$
where $\sigma\in \SG_{q+r}$, $\tau\in \SG_{p+r}$, 
$(z_1,\ldots,z_{q+r}) = (x^{(2)}_1,\ldots,x^{(2)}_q,x^{(3)}_1,
\ldots,x^{(3)}_r)$, 
$(t_1,\ldots,t_{p+r}) = (y^{(1)}_1,\ldots,y^{(1)}_p,y^{(3)}_1,\ldots,
y^{(3)}_r)$, and $\sigma$ preserves the order of the $r$ last elements of 
$\{1,\ldots,q+r\}$. 

If $\A$ is a Lie algebra and $r_\A\in \A\otimes\A$  is a solution of 
CYBE, then there are unique maps 
$$
\kappa^{\tensor}_{\A,r_\A} : F^{(3)} \to T\A\otimes U\A\otimes T\A
\quad \on{and} \quad
\kappa_{\A,r_\A} : F^{(3)} \to U\A^{\otimes  3},
$$ 
where $T\A$ is the tensor algebra of $\A$, 
such that if $r_\A = \sum_{i\in I} a_i \otimes b_i$, then 
$\kappa^{\tensor}_{\A,r_\A}$ maps $F^{(3)}_{pqr}$ to $\A^{\otimes q+r}
\otimes U\A\otimes \A^{\otimes p+r}$ 
in such a way that 
\begin{align*}
& \kappa^{\tensor}_{\A,r_\A}( 
\prod_{k = 1}^{q+r} z_{\sigma(k)} \otimes 
(\prod_{i = 1}^p x^{(1)}_{i} \otimes 
\prod_{j = 1}^q y^{(2)}_{j})
\otimes \prod_{l = 1}^{p+r} t_{\tau(l)} )
\\ & = \sum_{\al_1,\ldots,\al_{p+q+r}\in I}
\big( \bigotimes_{k = 1}^{q+r} a(\al_{p+\sigma(k)}) \big)  
\otimes \prod_{i = 1}^p a(\al_{i}) \prod_{j = 1}^q 
b(\al_{p+j})
\otimes
\big( \bigotimes_{j = 1}^{p+r} b(\beta(\tau,l)\big) ) ,  
\end{align*}
where $\beta(\tau,l)$ is equal to $\tau(l)$ if $\tau(l)\leq p$
and to $q + \tau(l)$ else, and $\kappa_{\A,r_\A}$ is the composition 
of $\kappa_{\A,r_\A}^{\tensor}$ with the projection map 
$T\A\to U\A$. 

\begin{prop} \label{prop:univ:pty}  
For any pair $(\A,r_\A)$ of
a Lie algebra and a solution of CYBE, $\kappa^{\tensor}_{\A,r_\A}$ is an 
algebra morphism from $(F^{(3)},m_{F^{(3)}})$ to $T\A\otimes U\A\otimes T\A$, 
and $\kappa_{\A,r_\A}$ is an algebra morphism from 
$(F^{(3)},m_{F^{(3)}})$ to $U\A^{\otimes 3}$. 
\end{prop}

{\em Proof.} The proof of the statement on $\kappa^{\tensor}_{\A,r_A}$
is by a double induction on $(q,q')$. The statement on 
$\kappa_{\A,r_A}$ follows immediately. 
\hfill \qed\medskip

\subsection{The normal ordering map}

For any triple of integers $p,q,r\geq 0$, let us form the tensor product 
$FA_{q+r}\otimes FA_{p+q}\otimes FA_{p+r}$; we denote the generators of its
factors as follows

-- generators of $FA_{q+r}$: 
$x^{(2)}_1,\ldots,x^{(2)}_q$,  $x^{(3)}_1,\ldots,x^{(3)}_r$

-- generators of $FA_{p+q}$: 
$y^{(2)}_1,\ldots,y^{(2)}_q$, $x^{(1)}_1,\ldots,x^{(1)}_p$  

-- generators of $FA_{p+r}$: 
$y^{(1)}_1,\ldots,y^{(1)}_p$,  $y^{(3)}_1,\ldots,y^{(3)}_r$. 

Each set of generators is split in two subsets (e.g., the first
subset of generators of $FA_{q+r}$ is  
$x^{(2)}_1,\ldots,x^{(2)}_q$ and the second subset 
is $x^{(3)}_1,\ldots,x^{(3)}_r$).

Then the group $\SG_p\times\SG_q \times\SG_r$
acts on this tensor product by simultaneous permutations of the 
three pairs of subsets of generators. We will set 
$$
G^{(3)}_{pqr} = \big( FA_{q+r}\otimes FA_{p+q}\otimes FA_{p+r} 
\big)_{\SG_p\times\SG_q\times \SG_r}. 
$$

Let $\Part_{p,q}$ be the set of pairs of partitions $(\up,\uq)$, 
where   $\up = (p_1,\ldots,p_\la)$ and $\uq = (q_1,\ldots,q_\la)$ 
are partitions of $p$ and $q$, such that   
$p_\la\geq 0$, $q_\la\geq 0$ and $p_\nu>0$, $q_{\nu'}>0$ when 
$\nu\neq \la$
and $\nu'\neq 1$. For any $(\up,\uq)\in\Part_{p,q}$, let us 
denote by $G^{\up,\uq}_{pqr}$ the image of   
$$
FA_{q+r} \otimes \prod_{\nu = 1}^{\la}
\big( \prod_{j = 1}^{q_\nu}  y^{(2)}_{q_1 + \cdots + q_{\nu - 1} + j} 
\prod_{i = 1}^{p_\nu}  x^{(1)}_{p_1 + \cdots + p_{\nu - 1} + i}
\big) \otimes FA_{p+r}   
$$
in $G^{(3)}_{pqr}$. Then we have 
$$
G^{(3)}_{pqr} = \bigoplus_{(\up,\uq)\in\Part_{p,q}} G^{\up,\uq}_{pqr} . 
$$

Let us define $M_{\uq,\up}$ as the subset of $M_{qp}$
consisting of all pairs of maps $(c,c')$, where  
$c : \{1,\ldots,q\}\to\{-\infty,1,\ldots,p\}$ and  
$c' : \{1,\ldots,p\}\to\{1,\ldots,q,\infty\}$ satisfy 
(\ref{ineqs}) and are such that for any $\nu = 1,\ldots,\la$, 
$$
c(\{1,\ldots,\sum_{i = 1}^{\nu} q_i \}) 
\subset \{-\infty,1,\ldots,\sum_{i = 1}^{\nu} p_i\} 
$$
and
$$
c'(\{\sum_{i = \nu}^{\la} p_i + 1,\ldots,p\}) \subset 
\{\sum_{i = \nu}^{\la} q_i + 1,\ldots,q,\infty\} . 
$$

Let us define the map 
$$
\mu^{\up,\uq}_{pqr} : G_{pqr}^{\up,\uq} \to F^{(3)}
$$
as follows. If $P\in FA_{q+r}$, $Q\in FA_{p+r}$, let us set 
$$ 
\mu_{pqr}^{\up,\uq}(P\otimes 
\prod_{\nu = 1}^\la 
\big( \prod_{j = 1}^{q_\nu}  y^{(2)}_{q_1 + \cdots + q_{\nu - 1} + j} 
\prod_{i = 1}^{p_\nu}  x^{(1)}_{p_1 + \cdots + p_{\nu - 1} + i}
\big)   \otimes  Q) = 
\sum_{(c,c')\in M_{\uq,\up}} 
\al_{p,q}^{\up,\uq}(P,Q) (c,c') ;   
$$
we set $\beta = \card (c')^{-1}(\infty)$ and $\beta' = \card
c^{-1}(-\infty)$, and define $\al_{p,q}^{\up,\uq}(P,Q,\gamma,\delta) (c,c')$
as the element of $F^{(3)}_{p-\beta,q-\beta',r+\beta + \beta'}$
equal to 
\begin{align*}
& \al_{p,q}^{\up,\uq}(P,Q) (c,c')
\\ & = 
P( 
\prod^{\nearrow}_{i_1\in (c')^{-1}(1)} \ad'(x^{(1)}_{i_1})(x^{(2)}_1),
\ldots,
\prod^{\nearrow}_{i_q\in (c')^{-1}(q)} \ad'(x^{(1)}_{i_q})(x^{(2)}_q),
x^{(3)}_1,\ldots,x^{(3)}_r) 
\\ & 
\otimes 
( \prod^{\nearrow}_{i\in (c')^{-1}(\infty)} x^{(1)}_{i}
\otimes 
\prod^{\nearrow}_{j\in c^{-1}(-\infty)} y^{(2)}_{j} )
\\ & 
\otimes 
Q( \prod^{\searrow}_{j_1\in c^{-1}(1)} 
\ad'(y^{(2)}_{j_1})(y^{(1)}_1),
\ldots,
\prod^{\searrow}_{j_p\in c^{-1}(p)} 
\ad'(y^{(2)}_{j_p})(y^{(1)}_p),
y^{(3)}_1,\ldots,y^{(3)}_r) 
\end{align*} 
where the generators of the components of 
$F^{(3)}_{p-\beta,q-\beta',r+\beta + \beta'}$ are denoted as follows

-- generators of $FA_{q+r+\beta}$. 
First set of generators: $x^{(2)}_j$, $j\in c^{-1}(-\infty)$. 
Second set of generators: $x^{(1)}_{j}$, 
$j\in (c')^{-1}(\{1,\ldots,q\})$, 
$x^{(2)}_i$, $i\in c^{-1}(\{1,\ldots,p\})$, 
and 
$x^{(3)}_1,\ldots,x^{(3)}_r$  

-- generators of $FA_{p-\beta}$: 
$x^{(1)}_{i}$, $i\in (c')^{-1}(\infty)$.  

-- generators of $FA_{q-\beta'}$:  
$y^{(2)}_{j}$, $j\in c^{-1}(-\infty)$. 

-- generators of $FA_{p+r+\beta'}$. 
First set of generators: $y^{(1)}_i$, $i\in (c')^{-1}(\infty)$. 
Second set of generators: $y^{(1)}_j$, $j\in (c')^{-1}(\{1,\ldots,q\})$, 
$y^{(2)}_{i}$, $i\in c^{-1}(\{1,\ldots,p\})$, 
and $y^{(3)}_1,\ldots,y^{(3)}_r$.   

Let $\A$ be a Lie algebra and let $r_\A\in \A\otimes \A$ be a solution of 
CYBE. Elements of the form 
$$
z_{\sigma(1)} \cdots z_{\sigma(q+r)}
\prod_{\nu = 1}^\la \big( \prod_{j = 1}^{q_\nu} 
y^{(2)}_{q_1 + \cdots + q_{\nu - 1} + j}
\prod_{i = 1}^{p_\nu} 
x^{(1)}_{p_1 + \cdots + p_{\nu - 1} + i} \big) 
\otimes 
w_{\tau(1)} \cdots w_{\tau(p+r)}, 
$$
form a generating family of $G^{(3)}_{pqr}$,  where  
we set 
$(z_1,\ldots,z_{q+r}) = (x^{(2)}_1,\ldots,x^{(2)}_q,
x^{(3)}_1,\ldots,x^{(3)}_r)$, 
$(w_1,\ldots,w_{p+r}) = (y^{(1)}_1,\ldots,y^{(1)}_p,
y^{(3)}_1,\ldots,y^{(3)}_r)$, 
and $(\up,\uq)$ belongs to $\Part_{p,q}$ and $(\sigma,\tau)$
belongs to $\SG_{q+r} \times \SG_{p+r}$. 
Then there is a unique map
$$
\gamma^{(\A,r_\A)}_{pqr} : G^{(3)}_{pqr} \to T\A\otimes U\A\otimes T\A,  
$$ 
such that for any $(\up,\uq)\in \Part_{p,q}$ and any 
$\sigma\in \SG_{q+r}$, $\tau\in \SG_{p+r}$, 
\begin{align*}
& \gamma^{(\A,r_\A)}_{pqr}(z_{\sigma(1)} \cdots z_{\sigma(q+r)}
\otimes \prod_{\nu = 1}^\la \big( \prod_{j = 1}^{q_\nu} 
y^{(2)}_{q_1 + \cdots + q_{\nu - 1} + j}
\prod_{i = 1}^{p_\nu} 
x^{(1)}_{p_1 + \cdots + p_{\nu - 1} + i} \big) 
\otimes 
w_{\tau(1)} \cdots w_{\tau(p+r)})
\\ & 
= \sum_{\al(1),\ldots,\al(p+q+r)\in I}
a(\al(p + \sigma(1))) \cdots a(\al(p + \sigma(q+r)))
\\ & 
\otimes  \prod_{\nu = 1}^\la 
\big( \prod_{j =1}^{q_\nu} b(\al(q_1 + \cdots + q_{\nu - 1}+ j))
\prod_{i =1}^{p_\nu} a(\al(p_1 + \cdots + p_{\nu - 1} + i)) \big) 
\otimes b(\epsilon(1)) \cdots b(\epsilon(p+r)) , 
\end{align*}
where $\epsilon(i) = \tau(j)$ if $\tau(j)\leq p$
and $\epsilon(j) = q + \tau(j)$ if $\tau(j)>p$, and we set 
$r_\A = \sum_{\al\in I} a(\al)\otimes b(\al)$. 

Let us define $\mu_{pqr} : G^{(3)}_{pqr} \to F^{(3)}$ as the direct sum of the maps
$\mu_{pqr}^{\up,\uq}$. Then $\mu_{pqr}$ has the following 
property. 

\begin{prop} \label{prop:normal:order}
We have $\kappa_{\A,r_\A}\circ\mu_{pqr} = \gamma^{(\A,r_\A)}_{pqr}$. 
\end{prop}

{\em Proof.} By induction. \hfill \qed\medskip

\begin{remark}  
The direct sum $\mu^{(3)} = \oplus_{(p,q,r)\in\NN^3}\mu_{pqr}$ may therefore 
be viewed as a universal version of the 
normal ordering of expressions in a solution $r_\A$ of CYBE. 
We will sometimes identify elements of $G_{pqr}$
with their images in $F^{(3)}$ by $\mu^{(3)}$.  
\end{remark}

\subsection{The CYBE identity in $F^{(3)}$}

We have the following identity in $F^{(3)}$ 
$$
\mu^{(3)}( x_1^{(1)} \otimes [y_1^{(1)},x^{(2)}_1] \otimes x^{(2)}_1)
= - \mu^{(3)}( [x^{(1)}_1,x^{(3)}_1]\otimes y^{(1)}_1\otimes y^{(3)}_1)
- \mu^{(3)}(x^{(3)}_1 \otimes x^{(2)}_1 \otimes [x^{(2)}_1,y^{(3)}_1])  , 
$$
in which the first expression belongs to $\mu^{(3)}(G^{(3)}_{110})$, the second to 
$\mu^{(3)}(G^{(3)}_{101})$ and the third to $\mu^{(3)}(G^{(3)}_{011})$. 
(The image of this identity by any map $\kappa_{\A,r_\A}$ 
simply expresses the fact that $r_\A$ satisfies CYBE.)

Let us define a map $\conc_{G^{(3)}} : G^{(3)}\otimes G^{(3)} \to G^{(3)}$ as follows. 
$\conc_{G^{(3)}}$ maps $G^{(3)}_{pqr} \otimes G^{(3)}_{p'q'r'}$ to 
$G^{(3)}_{p+p',q+q',r+r'}$; 
if $P,Q,R,P',Q',R'$ belong to $FA_{q+r}$, $FA_{p+q}$, $FA_{p+r}$, 
$FA_{q'+r'}$, $FA_{p'+q'}$, $FA_{p'+r'}$, then 
\begin{align*}
& \conc_{G^{(3)}}((P\otimes Q\otimes R)\otimes (P'\otimes Q'\otimes R'))
\\ & = 
P(x^{(2)}_1,\ldots,x^{(2)}_q,x^{(3)}_1,\ldots,x^{(3)}_r)
P'(x^{(2)}_{q+1},\ldots,x^{(2)}_{q+q'},x^{(3)}_{r+1},\ldots,x^{(3)}_{r+r'})
\\ & 
\otimes 
Q(x^{(1)}_1,\ldots,x^{(1)}_p,y^{(2)}_1,\ldots,y^{(2)}_q) 
Q'(x^{(1)}_{p+1},\ldots,x^{(1)}_{p+p'},y^{(2)}_{q+1},\ldots,y^{(2)}_{q+q'})
\\ & 
\otimes
R(y^{(1)}_1,\ldots,y^{(1)}_p,y^{(3)}_1,\ldots,y^{(3)}_r)
R'(y^{(1)}_{p+1},\ldots,y^{(1)}_{p+p'},y^{(3)}_{r+1},\ldots,y^{(3)}_{r+r'}) . 
\end{align*} 

Then we have 

\begin{lemma} \label{jury}
$(G^{(3)},\conc_{G^{(3)}})$ is an algebra and $\mu$ is an algebra morphism 
from $(G^{(3)},\conc_{G^{(3)}})$ to $(F^{(3)},m_{F^{(3)}})$.  
\end{lemma}

This follows from analysis of the behavior the spaces of maps 
$M_{pqr}^{\up\uq}$.  

We have then 

\begin{lemma} \label{lemma:CYBE:univ}
If $P\in FA_{q+r}$, $Q\in FA_{p+q-1}$ and $R\in FA_{p+r}$, we have 
the identity 
\begin{align*}
& \mu^{(3)}\big( P(x^{(2)}_1,\ldots,x^{(2)}_q,x^{(3)}_1,\ldots,x^{(3)}_r) 
\otimes Q(y^{(2)}_1,\ldots, [y^{(2)}_j,x^{(1)}_i],\ldots,y^{(2)}_q,
x^{(1)}_1,\ldots,\check{i},\ldots, x^{(1)}_p) 
\\ & \otimes R(y^{(1)}_1,\ldots,y^{(1)}_p,y^{(3)}_1,\ldots,y^{(3)}_q) \big)
\\ & 
= 
- \mu^{(3)} \big( P(x^{(2)}_1,\ldots,[x^{(1)}_i,x^{(2)}_j]\ldots,x^{(2)}_q,
x^{(3)}_1,\ldots,x^{(3)}_r) 
\otimes Q(y^{(2)}_1,\ldots,y^{(2)}_q,
x^{(1)}_1,\ldots,\check{i},\ldots, x^{(1)}_p) 
\\ & \otimes R(y^{(1)}_1,\ldots,y^{(1)}_p,y^{(3)}_1,\ldots,y^{(3)}_q) 
\big)
\\ & 
- \mu^{(3)} \big( P(x^{(2)}_1,\ldots,x^{(2)}_q,x^{(3)}_1,\ldots,x^{(3)}_r) 
\otimes Q(y^{(2)}_1,\ldots,x^{(1)}_i,\ldots,y^{(2)}_q,
x^{(1)}_1,\ldots,\check{i},\ldots, x^{(1)}_p) 
\\ & \otimes R(y^{(1)}_1,\ldots,[y^{(1)}_i,y^{(2)}_j],
\ldots,y^{(1)}_p,y^{(3)}_1,\ldots,y^{(3)}_q) \big), 
\end{align*}
where the arguments of $\mu$ belong to $G^{(3)}_{pqr}$, $G^{(3)}_{p-1,q,r+1}$ and 
$G^{(3)}_{p,q-1,r+1}$. 
\end{lemma}

{\em Proof.} If $x'$ and $x''$ are elements of $G^{(3)}$, then the image by 
$\mu^{(3)}$ of the product $\conc_{G^{(3)}}(x'\otimes \big(
x_1^{(1)} \otimes [y_1^{(1)},x^{(2)}_1] \otimes x^{(2)}_1
+ [x^{(1)}_1,x^{(3)}_1]\otimes y^{(1)}_1\otimes y^{(3)}_1
+ x^{(3)}_1 \otimes x^{(2)}_1 \otimes [x^{(2)}_1,y^{(3)}_1]  
\big) \otimes x'')$ is zero by Lemma \ref{jury}.  
\hfill \qed\medskip

\subsection{The maps $\mu_{Lie}^{p,q,r}$}

If $p$ and $q$ are integers $>0$ and $\al$ (resp., $\beta$) 
is an integer in $\{1,\ldots,p\}$ (resp., in $\{1,\ldots,q\}$), 
we will denote by $\wt M_{pq}^{\al 0}$  (resp., 
$\wt M_{pq}^{0 \beta}$) the set of all triples
$(c,c',\omega)$
consisting of a pair $(c,c')$ in  $M_{pq}^{\al 0}$  
(resp., in $M_{pq}^{0 \beta}$) and of order relations
on the sets $c^{-1}(1),\ldots,c^{-1}(q),(c')^{-1}(1),\ldots,
(c')^{-1}(p)$.

\begin{prop}
Let $(p,q)$ be a pair of integers $> 0$ and let 
$\al\in \{1,\ldots,p\}$ and $\beta\in \{1,\ldots,q\}$. 
Then there exist families of linear maps $(\Phi'(c,c',\omega))_{(c,c',\omega)\in 
\wt M_{pq}^{\al 0}}$ 
and $(\Phi''(c,c',\omega))_{(c,c',\omega)\in \wt M_{pq}^{0\beta}}$, indexed by 
$\wt M_{pq}^{\al 0}$  and $\wt M_{pq}^{0\beta}$  respectively, such that 
$\Phi'(c,c',\omega)$ maps $FL_{p+q}$ to $FL_{p - \al}$ and $\Phi'(c,c',\omega)$
maps $FL_{p+q}$ to $FL_{q-\beta}$, and such that the following identities
are satisfied 
\begin{align} \label{Lie:id}
& P(x^{(2)}_1,\ldots,x^{(2)}_q,x^{(3)}_1,\ldots,x^{(3)}_r)
\otimes Q(x^{(1)}_1,\ldots,x^{(1)}_p, y^{(2)}_1,\ldots,y^{(2)}_q)
\\ & \nonumber
\otimes R(y^{(1)}_1,\ldots,y^{(1)}_p,y^{(3)}_1,\ldots,y^{(3)}_r) 
\\ & \nonumber 
=  \sum_{\al = 1}^{p} \sum_{(c,c',\omega)\in \wt M_{pq}^{\al 0}}
P( \prod_{k_1\in (c')^{-1}(1)} \ad'(x^{(1)}_{k_1})(x^{(2)}_1), 
\ldots,
\prod_{k_q\in (c')^{-1}(q)} \ad'(x^{(1)}_{k_q})(x^{(2)}_q),
x^{(3)}_1,\ldots,x^{(3)}_r) 
\\ & \nonumber  
\otimes \Phi'(c,c',\omega)(Q)(x^{(1)}_j,j\in (c')^{-1}(\infty)) 
\\ & \nonumber
\otimes 
R(\prod_{l_1\in c^{-1}(1)} \ad'(y^{(2)}_{l_1})
(y^{(1)}_1), \ldots,
\prod_{l_p\in c^{-1}(p)} \ad'(y^{(2)}_{l_p})
(y^{(1)}_p), y^{(3)}_1,\ldots,y^{(3)}_r) 
\\ & \nonumber  
+  \sum_{\beta = 1}^{q} \sum_{(c,c')\in \wt M_{pq}^{0\beta}}
P( \prod_{k_1\in (c')^{-1}(1)} \ad'(x^{(1)}_{k_1})(x^{(2)}_1), 
\ldots,
\prod_{k_q\in (c')^{-1}(q)} \ad'(x^{(1)}_{k_q})(x^{(2)}_q),
x^{(3)}_1,\ldots,x^{(3)}_r) 
\\ & \nonumber 
\otimes \Phi''(c,c',\omega)(Q)(y^{(2)}_k,k\in c^{-1}(\infty)) 
\\ & \nonumber 
\otimes 
R(\prod_{l_1\in c^{-1}(1)} \ad'(y^{(2)}_{l_1})
(y^{(1)}_1), \ldots,
\prod_{l_p\in c^{-1}(p)} \ad'(y^{(2)}_{l_p})
(y^{(1)}_p), y^{(3)}_1,\ldots,y^{(3)}_r) 
\end{align}
for any integer $r\geq 0$ and any pair of elements 
$(P,R)\in FA_{q+r} \times FA_{p+r}$. 
\end{prop}

{\em Proof.} Let us construct the maps $\Phi(c,c',\omega)$ and $\Phi''(c,c',\omega)$
by induction on $p+q$. Assume that these maps are constructed 
for all pairs of integers $(p',q')$ with $p'+q'<p+q$. 
Let $(Q_\al)_\al$ be a basis of $FL_{p+q}$. We may assume that 
each $Q_\al$ is of the form $[Q'_\al,x^{(1)}_i]$, where
$i\in \{1,\ldots,p\}$, or $[Q'_\al,y^{(2)}_j]$, where $j\in 
\{1,\ldots,q\}$, and $Q'_\al\in FL_{p+q-1}$.  
Let us treat the latter case.  
Then 
\begin{align*} 
& P(x^{(2)}_1,\ldots,x^{(2)}_q,x^{(3)}_1,\ldots,x^{(3)}_r)
\otimes [Q'_\al(x^{(1)}_1,\ldots,x^{(1)}_p, 
y^{(2)}_1,\ldots,\check y^{(2)}_j, \ldots,y^{(2)}_q),y^{(2)}_j]
\\ & 
\otimes R(y^{(1)}_1,\ldots,y^{(1)}_p,y^{(3)}_1,\ldots,y^{(3)}_r) 
\end{align*}
may be rewritten as follows 
\begin{align*}
& \sum_{\al = 1}^{p} \sum_{(c,c',\omega)\in \wt M_{p,q-1}^{\al 0}}
P( \prod_{k_1\in (c')^{-1}(1)} \ad'(x^{(1)}_{k_1})(x^{(2)}_1), 
\ldots,
\prod_{k_q\in (c')^{-1}(q)} \ad'(x^{(1)}_{k_q})(x^{(2)}_q),
x^{(3)}_1,\ldots,x^{(3)}_r) 
\\ & 
\otimes [\Phi'(c,c',\omega)(Q'_\al)(x^{(1)}_j,j\in (c')^{-1}(\infty)),y^{(2)}_j] 
\\ & 
\otimes 
R(\prod_{l_1\in c^{-1}(1)} \ad'(y^{(2)}_{l_1})
(y^{(1)}_1), \ldots,
\prod_{l_p\in c^{-1}(p)} \ad'(y^{(2)}_{l_p})
(y^{(1)}_p), y^{(3)}_1,\ldots,y^{(3)}_r) 
\\ & + 
\sum_{\beta = 1}^{q-1} \sum_{(c,c',\omega)\in \wt M_{p,q-1}^{0\beta}}
P( \prod_{k_1\in (c')^{-1}(1)} \ad'(x^{(1)}_{k_1})(x^{(2)}_1), 
\ldots,
\prod_{k_q\in (c')^{-1}(q)} \ad'(x^{(1)}_{k_q})(x^{(2)}_q),
x^{(3)}_1,\ldots,x^{(3)}_r) 
\\ & 
\otimes [\Phi''(c,c',\omega)(Q'_\al)(y^{(2)}_k,k\in c^{-1}(\infty)),y^{(2)}_j] 
\\ & 
\otimes 
R(\prod_{l_1\in c^{-1}(1)} \ad'(y^{(2)}_{l_1})
(y^{(1)}_1), \ldots,
\prod_{l_p\in c^{-1}(p)} \ad'(y^{(2)}_{l_p})
(y^{(1)}_p), y^{(3)}_1,\ldots,y^{(3)}_r) , 
\end{align*}
applying the identity (\ref{Lie:id}) to  $P\otimes Q'_\al \otimes R y^{(2)}_j$
and transferring the extreme right term $y^{(2)}_j$ to the middle tensor factor 
via the adjoint action. The second sum is of the desired form; we transform
the first sum writing that 
\begin{align*}
& [\Phi'(c,c',\omega)(Q'_\al)(x^{(1)}_{j'},j'\in (c')^{-1}(\infty))
,y^{(2)}_j] 
\\ & = \sum_{j'\in (c')^{-1}(\infty))}
\Phi'(c,c',\omega)(Q'_\al)([x^{(1)}_{j''},x^{(1)}_{j'}
,y^{(2)}_j] ,j''\in (c')^{-1}(\infty) - \{j'\}) . 
\end{align*}
and using the CYBE identity of Lemma \ref{lemma:CYBE:univ} 
to lower the degree of the 
middle term. This procedure and the condition that 
(\ref{Lie:id}) holds defines uniquely 
the $\Phi'(c,c',\omega)(Q_\al)$  and $\Phi''(c,c',\omega)(Q_\al)$. Of course, 
the maps $\Phi'(c,c',\omega)$ and  $\Phi''(c,c',\omega)$ are far from 
unique because of the many possible systems of bases of the $FL_{p+q}$. 
\hfill \qed\medskip 

Recall that we defined 
$$
F^{(aab)} 
= \oplus_{p + q\geq 0} (FL_p \otimes FL_q \otimes FL_{p+q})_{\SG_p\times\SG_q}, 
\quad  
F^{(abb)}  = 
\oplus_{p + q \geq 0} (FL_n \otimes FL_p \otimes FL_q)_{\SG_p\times\SG_q}  
$$

\begin{lemma}
The natural inclusions of $FL_k$ in $FA_k$ ($k = p,q,p+q$) induce  
an inclusion of $F^{(aab)}\oplus F^{(abb)}$ in $F^{(3)}$. 
\end{lemma}

{\em Proof.} Each $FL_k$ is a $\SG_k$-submodule of $FA_k$, therefore 
$FL_p \otimes FL_q \otimes FL_{p+q}$ is a $\SG_p\times\SG_q$-submodule of 
$FA_p \otimes FA_q \otimes FA_{p+q}$. The result now follows from the 
fact that if $\Gamma$ is any finite group and $N\subset M$ is an inclusion of 
$\Gamma$-modules, then the natural map $N_\Gamma\to M_\Gamma$ is injective. 
(This fact is proven as follows: the  representations of $\Gamma$ are completely 
reducible, so  
$M_\Gamma$ identifies with the multiplicity space of the trivial representation. 
Then $N\subset M$ corresponds to embeddings of the multiplicity spaces of each 
simple $\Gamma$-module, which shows that $N_\Gamma\subset M_\Gamma$.)
\hfill \qed\medskip 

\begin{cor} \label{cor:Lie}
Let $p,q,r$ be integers $>0$ and let 
$P,Q,R$ be elements of $FL_{q+r},FL_{p+q}$ and $FL_{p+r}$. 
Then the element of $F^{(3)}$ equal to   
\begin{align*}  
& P(x^{(2)}_1,\ldots,x^{(2)}_q,x^{(3)}_1,\ldots,x^{(3)}_r)\otimes 
Q(x^{(1)}_1,\ldots,x^{(1)}_p,y^{(2)}_1,\ldots,y^{(2)}_q) 
\\ & 
\otimes R(y^{(1)}_1,\ldots,y^{(1)}_p,y^{(3)}_1,\ldots,y^{(3)}_r)
\end{align*}
of $\mu(G_{pqr})$
belongs to $F^{(aab)}\oplus F^{(abb)}$. 
\end{cor}

{\em Proof.} If $P$ and $R$ are Lie polynomials, then the expressions in 
the first and third tensor factors of the right side of (\ref{Lie:id})
are also Lie polynomials. 
\hfill \qed\medskip  

$\mu^{p,q,r}_{Lie}(P\otimes Q\otimes R)$ is then defined as this element. 
\hfill \qed\medskip

\newpage

\section{Construction and properties of $\delta_4^{(F)}$ (proof of Props.\  
\ref{prop:delta:4} and \ref{prop:delta:diff})}
\label{proof:prop:delta:4}

In this Section, we first generalize the construction and the 
properties of the algebra 
$F^{(3)}$ to integers $n\geq 2$. The statements on $\delta_4^{(F)}$
will be immediate corollaries of these constructions. 

\subsection{Definition of $F^{(n)}$}

If $n$ is an integer $\geq 2$, let 
us define 
$$
P_n   = \{\up = (p_{ij})_{1\leq i < j\leq n} | p_{ij}\geq 0\}. 
$$   
When $\up$ belongs to $P_n$, let us set 
$$
F^{(n)}_\up = \big(\bigotimes_{i = 1}^n (FA_{\sum_{j|j<i} p_{ji}}
\otimes FA_{\sum_{j|j>i} p_{ij}}) \big)_{\prod_{1\leq i<j\leq n} 
\SG_{p_{ij}}} .  
$$
Here the generators of the first part of the $i$th factor
are denoted $x^{(ji)}_{\al}$, where $(j,\al)$ are such that 
$j<i$ and $1\leq \al\leq p_{ji}$, and the generators of the 
second part of the $i$th factor are denoted 
$y^{(ij)}_{\al}$, where $(j,\al)$ are such that 
$j>i$ and $1\leq \al\leq p_{ij}$. There are therefore 
$n(n-1)$ sets of generators, indexed by the pairs $(i,j)$
such that $(i,j)\in\{1,\ldots,n\}^2$ and $i\neq j$. 
Let us put $S_{ij} = \{y^{(ij)}_\al|\al = 1,\ldots,p_{ij}\}$ 
when $i<j$ and  $S_{ij} = \{x^{(ji)}_\al|\al = 1,\ldots,p_{ji}\}$ 
when $i>j$. If $i<j$, there is a bijection between 
$S_{ij}$ and $S_{ji}$, sending each $y^{(ij)}_\al$ to $x^{(ij)}_\al$. 
Then $\SG_{p_{ij}}$ acts by simultaneaous permutation of the 
sets $S_{ij}$ and $S_{ji}$ and therefore also on 
$\bigotimes_{i = 1}^n (FA_{\sum_{j|j<i} p_{ji}}
\otimes FA_{\sum_{j|j>i} p_{ij}})$. $F^{(n)}_\up$ is then 
defined as the space of coinvariants of this action. 
We set 
$$
F^{(n)} = \bigoplus_{\up\in P_n} F^{(n)}_\up . 
$$

\subsection{Basis of $F^{(n)}$}

When $\us\in\NN^n$, let us set $I_\us = \{(i,\al)|1\leq i\leq n,
1\leq \al\leq s_i\}$. Let us denote by $\ind$ the map from 
$I_\us$ to $\{1,\ldots,n\}$ such that $\ind(i,\al) = i$. We set
$I_{\us,i} = \ind^{-1}(\{i\})$.  

Let us define $\Phi_n$ as the set of all triples $(\us,\ut,\phi)$, 
where $\us$ and $\ut$ belong to $\NN^n$ and $\phi$ is a 
bijection from $I_\us$ to $I_\ut$, such that for any $(i,\al)\in 
I_\us$, $\ind(\phi(i,\al))>i$. (If $(\us,\ut,\phi)\in \Phi_n$, 
we have therefore $\sum_{i = 1}^n s_i = \sum_{i = 1}^n t_i$ 
and $s_1 = t_n = 0$.)

There is a unique map $\pi_n$ from $\Phi_n$ to $P_n$, such that 
for any $(\us,\ut,\phi)$ in $\Phi_n$, $\pi_n(\us,\ut,\phi)$
is the element $\up = (p_{ij})_{1\leq i < j \leq n}$ of $P_n$
such that $p_{ij} = \card(\phi(I_{\us,i})\cap I_{\ut,j})$. 
In particular, if $\up = \pi_n(\us,\ut,\phi)$, then 
we have $s_i = \sum_{j|j > i} p_{ij}$    and 
$t_i = \sum_{j|j < i} p_{ji}$. 

If $(\us,\ut,\phi)\in \Phi_n$, let us set 
$$
A_{ij} = I_{\us,i}\cap \phi^{-1}(I_{\ut,j}), \quad 
B_{ij} = \phi(I_{\us,i}) \cap I_{\ut,j} . 
$$    
Then $A_{ij}$ and $B_{ij}$) are empty when $j\leq i$;  
$(A_{i\al})_{\al | \al > i}$is a partition of $I_{\us,i}$ and  
$(B_{\al i})_{\al | \al < i}$ is a partition of $I_{\ut,i}$. 

If  $(\us,\ut,\phi)\in P_n$, define $z(\us,\ut,\phi)$ in 
$F^{(n)}_{\pi_n(\us,\ut,\phi)}$ as 
$$
z(\us,\ut,\phi) = \bigotimes_{i = 1}^n 
(z^{(i)}_1\cdots z^{(i)}_{s_i} ) \otimes ( w^{(i)}_1\cdots 
w^{(i)}_{t_i}) , 
$$  
where if $(a_{ij,\al})_{\al = 1,\ldots,p_{ij}}$ is the 
increasing sequence such that $A_{ij} = \{a_{ij,1},\ldots,
a_{ij,p_{ij}}\}$, then 
$$
z^{(i)}_{a_{ij,\al}} = x^{(ij)}_\al \quad \on{and} \quad  
w^{(j)}_{\phi(a_{ij,\al})} = y^{(ij)}_\al. 
$$

\begin{lemma}
The family $(z(\us,\ut,\phi))_{(\us,\ut,\phi)\in P_n}$
is a basis of $F^{(n)}$. 
\end{lemma}

The families $(z^{(i)}_\al)$ and $(w^{(i)}_\beta)$ may also be
defined as follows. Let $(b_{ij,\al})_{\al = 1,\ldots,p_{ij}}$
be the increasing sequence such that $B_{ij} = 
\{b_{ij,1},\ldots,b_{ij,p_{ij}}\}$, then 
$$
z^{(i)}_{\phi^{-1}(b_{ij,\al})} = x^{(\al)}_{ij} 
\quad \on{and} \quad 
w^{(j)}_{b_{ij,\al}} = y^{(ij)}_\al .  
$$

\subsection{Product in $F^{(n)}$}

If $(\us,\ut,\phi)$ and $(\us',\ut',\phi')$ are elements of $P_n$, 
and if for each $\al = 1,\ldots,n$,  
$(c_\al,c'_\al)$ is an element of $M_{t_\al,s'_\al}$ let us define  
$\up((\us,\ut,\phi),(\us',\ut',\phi'),(c_\al,c'_\al)_{\al = 1,\ldots,n})$
as the element of $P_n$ such that if  $1\leq i < j\leq n$, then 
\begin{align*}
& \up ((\us,\ut,\phi),(\us',\ut',\phi'),(c_\al,c'_\al)_{\al = 1,\ldots,n}) 
\\ &  = 
\card(c_i^{-1}(-\infty)\cap A'_{ij}) + \card(c_j^{\prime -1}(\infty)
\cap B_{ij} ) 
\\ & 
+ \sum_{k = i + 1}^{j - 1} \card(c_k^{\prime -1}(A'_{kj})\cap B_{ik}) 
+ \card(c_k^{-1}(B_{ik}) \cap A'_{kj}) 
\end{align*} 
where $(A'_{ij})_{j|j>i}$ and $(B'_{ij})_{j|j<i}$ denote the partitions 
of $I_{\us',i}$ and $I_{\ut',i}$ associated to $(\us',\ut',\phi')$.    
Let us define now $z((\us,\ut,\phi),(\us',\ut',\phi'),
(c_\al,c'_\al)_{\al = 1,\ldots,n})$ as the element of 
$F^{(n)}_{\up ((\us,\ut,\phi),(\us',\ut',\phi'),(c_1,c'_1), 
\ldots,(c_n,c'_n))}$ equal to 
\begin{align} \label{asnieres}
& z((\us,\ut,\phi),(\us',\ut',\phi'),
(c_\al,c'_\al)_{\al = 1,\ldots,n})
\\ & \nonumber 
= \bigotimes_{i = 1}^n 
\big( \prod^{\nearrow}_{\al\in c_{j_{i,1}}^{-1}(\phi(i,1))}
\ad'(z_{\al}^{\prime (j_{i,1})})(z^{(i)}_1)
\cdots 
\prod^{\nearrow}_{\al\in c_{j_{i,s_i}}^{-1}(\phi(i,s_i))}
\ad'(z_{\al}^{\prime (j_{i,s_i})})(z^{(i)}_{s_i})
\prod^{\nearrow}_{\al\in c_i^{-1}(-\infty)} z_\al^{\prime(i)}
\big)
\\ & \nonumber 
\otimes
\big( \prod^{\nearrow}_{\al\in c_i^{\prime -1}(\infty)} w_\al^{(i)}
\prod^{\searrow}_{\al\in c^{\prime -1}_{k_i,1}(\phi^{-1}(i,1))} 
\ad'( w_\al^{(k_{i,1})} )(w^{\prime(i)}_1) 
\cdots 
\prod^{\searrow}_{\al\in c^{\prime -1}_{k_i,t_i}(\phi^{-1}(i,t_i))} 
\ad'( w_\al^{(k_{i,t_i})} )(w^{\prime(i)}_{t_i}) 
\big) , 
\end{align}
where we set $j_{i,1} = \ind(\phi(i,1)), \ldots$, $j_{i,s_i} = 
\ind(\phi(i,s_i))$ and $k_{i,1} = \ind(\phi^{-1}(i,1)), \ldots$,
$k_{i,t_i} = \ind(\phi^{-1}(i,t_i))$. 
 
In the right side of (\ref{asnieres}), the sets of generators of 
$F^{(n)}_{\up ((\us,\ut,\phi),(\us',\ut',\phi'),(c_1,c'_1), 
\ldots,(c_n,c'_n))}$ are the following: 

-- if $i<j$, $S_{ij}$ consists of the $z^{\prime (i)}_\al$, 
$\al\in c_i^{-1}(-\infty)\cap A'_{ij}$, the $z^{(i)}_\beta$, 
where $(i,\beta)$ belongs to $\phi^{-1}(c_j^{\prime -1}(\infty)\cap B_{ij})$, 
the $z^{\prime (k)}_\gamma$, where $\gamma\in c_k^{-1}(B_{ik})\cap A'_{kj}$
and the $z^{(i)}_\delta$, where $(i,\delta)\in 
\phi^{-1}\big((c'_k)^{-1}(A'_{kj}\cap B_{ik})\big)$

-- if $i<j$, $S_{ji}$ consists of the $w^{\prime(j)}_\beta$, 
where $\beta\in \phi(c_i^{-1}(\infty)\cap A'_{ij})$, 
the $w^{(j)}_\al$, where $\al\in (c'_j)^{-1}(\infty)\cap B_{ij}$, 
the $w^{\prime (j)}_\delta$, where $(j,\delta)\in \phi\big( c_k^{-1}(B_{ik})
\cap A'_{kj}\big)$, and the $w^{(k)}_\gamma$, where 
$\gamma\in (c'_k)^{-1}(A'_{kj})\cap B_{ik}$.   

The bijection from $S_{ij}$ to $S_{ji}$ maps each $z^{\prime (i)}_\al$
to $w^{\prime (j)}_{\beta}$, where $(j,\beta) = \phi(i,\al)$, 
each $z^{(i)}_\beta$ to $w^{(j)}_\al$, where $(j,\al) = \phi(i,\beta)$, 
each $z^{\prime (k)}_\gamma$ to $w^{\prime (j)}_\delta$, where 
$(j,\delta) = \phi(k,\gamma)$ and each $z^{(i)}_\delta$ to 
$w^{(k)}_\gamma$, where $(k,\gamma) = \phi(i,\delta)$. 

Define $m_{F^{(n)}}$ as the unique linear map from $F^{(n)}\otimes F^{(n)}$
to $F^{(n)}$, such that for any $(\us,\ut,\phi)$ and $(\us',\ut',\phi')$
in $P_n$, we have 
$$
m_{F^{(n)}} (z(\us,\ut,\phi) \otimes z(\us',\ut',\phi') ) = 
\sum_{(c_1,c'_1)\in M_{t_1,s'_1}, 
\ldots, (c_n,c'_n)\in M_{t_n,s'_n}} 
$$

A construction analogous to that of Appendix \ref{proof:prop:key}  
shows 

\begin{prop}
$(F^{(n)}, m_{F^{(n)}})$ is an associative algebra. 
\end{prop}

\begin{remark}
When $n = 3$, the connection between this presentation of $F^{(3)}$ and that of 
Appendix \ref{proof:prop:key}   relies on the identification 
$(p_{12},p_{13},p_{23}) = (q,r,p)$. 

On the other hand, $F^{(2)}$ coincides with the direct sum 
$\oplus_{n\geq 0} \KK\SG_n$, where the product is the linear
extension of the concatenation of permutations. 
\end{remark}

\subsection{Universal property of $F^{(n)}$} \label{univ:pties:Fn}

Let $(\A,r_\A)$ be the pair of a Lie algebra $\A$ and a solution $r_\A
= \sum_{i\in I} a_i\otimes b_i\in 
A\otimes A$ of CYBE. Then there is a unique map 
$$
\kappa^{(n)}_{\A,r_\A} : F^{(n)} \to T\A\otimes U\A\otimes T\A
$$
such that for any $(\us,\ut,\phi)\in P_n$, 
$$
\kappa^{(n)}_{\A,r_\A}(z(\us,\ut,\phi)) =  
\sum_{\al\in \on{Map}(I_\us,I)} \bigotimes_{i = 1}^n
\big( \prod_{\gamma\in I_{\us,i}} a(\al(i,\gamma)) \prod_{\delta\in I_{\ut,i}}
b(\al(\phi^{-1}(i,\delta)))\big);  
$$
we denote by $\on{Map}(I_\us,I)$ the set of all maps from $I_\us$ to $I$.  

\begin{prop} \label{prop:universal:fn}
$\kappa^{(n)}_{\A,r_\A}$ is a morphism of algebras. 
\end{prop}

\subsection{The normal ordering map}

When $\up\in P_n$, let us set 
$$
G^{(n)}_\up = 
\big( \bigotimes_{i = 1}^n  FA_{\sum_{j|j<i} p_{ji} 
+  \sum_{j|j>i} p_{ij}}\big)_{\prod_{1\leq i<j\leq n} \SG_{p_{ij}}} . 
$$
The generators of the $i$th factor of the tensor product are 
$u^{(i)}_1,\ldots,u^{(i)}_{\sum_{j|j<i} p_{ji} 
+  \sum_{j|j>i} p_{ij}}$.
Let us denote by $S_i$ this set of generators. We have a partition 
$S_i = \cup_{j|j \neq i} S_{ij}$, where 

-- if $j<i$, $S_{ij}$ consists of the 
$u^{(i)}_{\sum_{j'|j'<j} p_{j'i} + 1}, \ldots,  
u^{(i)}_{\sum_{j'|j'\leq j} p_{j'i}}$

-- if $j>i$, $S_{ij}$ consists of the 
$u^{(i)}_{\sum_{j'|j'<i} p_{j'i}  
+ \sum_{j'|i<j'<j} p_{ij'} + 1}, \ldots,  
u^{(i)}_{\sum_{j'|j'<i} p_{j'i} + \sum_{j'|i<j'\leq j} p_{ij'}}$. 

When $i<j$, there is a bijection from $S_{ij}$ to $S_{ji}$, 
sending each $u^{(i)}_{\sum_{j'|j'<i} p_{j'i}
+ \sum_{j'|i<j'<j} p_{ij'} + \al}$ to $u^{(j)}_{\sum_{k'|k'<j} p_{k'j}
+ \al}$, for 
$\al = 1,\ldots,p_{ij}$. Then $\SG_{p_{ij}}$ acts by simultaneously
permuting the elements of $S_{ij}$ and of $S_{ji}$.

When $(\sigma_i)_{i = 1,\ldots,n}$ is a collection of permutations
of $\prod_{i = 1}^n \SG_{\sum_{j|j<i} p_{ji} + \sum_{j|j>i} p_{ij}}$, 
let us set 
$$
u\big((\sigma_i)_{i = 1,\ldots,n} \big)
= \bigotimes_{i = 1}^n \big( u^{(i)}_{\sigma_i(1)} 
\cdots u^{(i)}_{\sum_{j|j<i} p_{ji} + \sum_{j|j>i} p_{ij}}  \big) . 
$$
Then for any $\up\in P_n$, the family $\{ 
u\big((\sigma_i)_{i = 1,\ldots,n} \big) | 
(\sigma_i)_{i = 1,\ldots,n} \in 
\prod_{i = 1}^n \SG_{\sum_{j|j<i} p_{ji} + \sum_{j|j>i} p_{ij}}\}$
is a generating family of $G^{(n)}_\up$.   

When $\up\in P_n$, let us set $\SG(\up) = \prod_{i = 1}^n 
\SG_{\sum_{j|j<i} p_{ji} + \sum_{j|j>i} p_{ij}}$. If $\up$
and $\uq$ belong to $P_n$ and $\usigma = (\sigma_i)_{i = 1,\ldots,n}$
(resp., $\utau = (\tau_i)_{i = 1,\ldots,n}$) belongs to 
$\SG(\up)$ (resp., to $\SG(\uq)$), define 
$\usigma*\utau$ as the element of $\SG(\up + \uq)$ such that 
for any $i = 1,\ldots,n$, $((\usigma*\utau)_i)(\al) = \sigma_i(\al)$ if $\al \leq 
\sum_{j|j>i} p_{ij} + \sum_{j|j<i} p_{ji}
$, and  $((\usigma*\utau)_i)(\al) = 
\tau_i(\al - (\sum_{j|j>i} p_{ij} + \sum_{j|j<i} p_{ji}))
+ \sum_{j|j>i} p_{ij} + \sum_{j|j<i} p_{ji}$ if $\al \leq \sum_{j|j>i}
p_{ij} + \sum_{j|j<i} p_{ji}$. 

Then there is a unique linear map $\conc_{G^{(n)}}$ from 
$G^{(n)} \otimes G^{(n)}$ to $G^{(n)}$, such that if $\up,\uq
\in P_n$ and $\usigma\in \SG(\up)$, $\utau\in \SG(\uq)$, then 
$$
\conc_{G^{(n)}} (u(\usigma)\otimes u(\utau)) = u (\usigma*\utau) .   
$$

When $k = 1,\ldots,n-1$, define $G^{k,(n)}$ as  the subspace of 
$G^{(n)}$ spanned by all $u(\usigma)$, where $\usigma$
is such that when $k'\leq k$, the generators of $\cup_{k'|k'<k} S_{kk'}$  
occur before the generators of $\cup_{k'|k'>k} S_{k'k}$ in the $k$th
factor. Then $F^{(n)}$ and $G^{(n)}$ may be identified with 
$G^{n-1,(n)}$ and $G^{1,(n)}$. 
When $k = 1,\ldots,n-2$, define $\mu_k$ as the following 
partial ordering map. $\mu_k$ maps $G^{k,(n)}$  to $G^{k+1,(n)}$. 
If $\up\in P_n$, set 
$p(i,j) = p_{ij}$ if $i<j$ and $p(i,j) = p_{ji}$ if $i>j$. Then 
if $\usigma\in \SG(\up)$, then one may find elements of free algebras
$P_j$, $j\in \{1,\ldots,n\} - \{k\}$ and partitions $\{1,\ldots,p(k,k')\}
= \cup_{r = 1}^s I_{rk'}$, such that $u(\usigma)$ has the form  
\begin{align*}
& \bigotimes_{j|j<k} P_j(u_\al^{(jj')}|j'\neq j, \al= 1,\ldots,p(j,j')) 
\otimes
\\ & 
\prod_{r = 1}^s \big(\prod_{i|i\neq k} (\prod_{\al\in I_{ri}} 
u^{(ki)}_\al) \big) 
\otimes 
\bigotimes_{j|j>k} P_j(u_\al^{(jj')}|j'\neq j, \al= 1,\ldots,p(j,j')). 
\end{align*}
Let us set $q_t = \sum_{t'<k}\card(I_{tt'})$ and  
$p_t = \sum_{t'>k}\card(I_{tt'})$ for $t = 1,\ldots,s$. 
When $u<t$ and $\al = 1,\ldots,\card(I_{tu})$, let us set 
$w'_{q_1 + \cdots + q_{t-1} + \card(I_{t1}) + \cdots + \card(I_{t,u-1})
+ \al} = u^{(ku)}_{n_\al}$, where $n_\al$ is the $\al$th element of $I_{tu}$. 
When $u>t$ and  $\al = 1,\ldots,\card(I_{tu})$, let us set 
$v_{p_1 + \cdots + p_{t-1} + \card(I_{t,k+1}) + \cdots + 
\card(I_{t,u-1})
+ \al} = u^{(ku)}_{n_\al}$. 
Let us also set $v'_p = u^{(uk)}_{\nu}$ whenever $w'_p = u^{(ku)}_{\nu}$
and $w_p = u^{(ku)}_{\nu}$  whenever $v_p = u^{(uk)}_{\nu}$. 

Let us set 
$A_i = \cup_{t = 1}^s (q_1 + \cdots + q_{t-1} + 
\card(I_{t1}) + \cdots + \card(I_{t,i-1}) 
+ \{1,\ldots,\card(I_{ti})\})$ and $B_j = 
\cup_{t = 1}^s (p_1 + \cdots + p_{t-1} + 
\card(I_{t,k+1}) + \cdots + \card(I_{t,j-1}) 
+ \{1,\ldots,\card(I_{tj})\})$. 

Let us set $\up = (p_1,\ldots,p_s)$ and $\up = (q_1,\ldots,q_s)$. 
For $(c,c')\in M_{\uq,\up}$, let us define $(p_{ij}(c,c'))_{1\leq i<j\leq n}$
as follows. If $i<j<k$ or $k<i<j$, then $p_{ij}(c,c') = p_{ij}$. 
If $i<k<j$, then $p_{ij}(c,c') = p_{ij} 
+ \card(A_j\cap (c')^{-1}(B_i))$. If $i<k$, then $p_{ik} = 
\card(c^{-1}(-\infty)\cap A_i)$.  
If $k<j$, then $p_{kj} = \card((c')^{-1}(\infty)\cap B_j)$. 

When $(c,c')\in M_{\uq,\up}$, define $u(k,\usigma),c,c') $ as the element 
of $G^{(n)}_{(p_{ij}(c,c'))_{1\leq i < j \leq n}}$ equal to 
\begin{align*}
&
\bigotimes_{j|j<k} P_j(u_\al^{(jj')}| j'\neq j,k; \al= 1,\ldots,p(j,j'); 
\prod^{\searrow}_{j\in (c')^{-1}(i')} \ad'(v_j)(v'_{i'}), i'\in A_j)
\\ & 
\otimes 
\prod^{\nearrow}_{i'\in c^{-1}(-\infty)} w'_{i'} 
\prod^{\nearrow}_{i\in (c')^{-1}(\infty)} v_i
\\ & 
\otimes
\bigotimes_{j|j>k} P_j(u_\al^{(jj')}| j'\neq j,k; \al= 1,\ldots,p(j,j'); 
\prod^{\nearrow}_{j'\in c^{-1}(i)} \ad'(w'_{j'})(w_{i}), i\in B_j))
\end{align*}

Then there is a unique map $\mu_k : G^{k,(n)} \to G^{k+1,(n)}$
such that $\mu_k(u(\usigma)) = \sum_{(c,c')\in M_{\uq,\up}}
u(k,\usigma),c,c')$.

Define now the linear map $\mu^{(n)} : F^{(n)} \to G^{(n)}$
as the composition $\mu_{n-2}\circ\cdots \circ \mu_2$.  

Then $\mu^{(n)}$ is an algebra morphism from $(F^{(n)},m_{F^{(n)}})$
to $(G^{(n)},\conc_{G^{(n)}})$. $\mu^{(n)}$ is a normal ordering 
map is a sense generalizing Proposition \ref{prop:normal:order}.

Moreover, normal ordering commutes with the operation of 
``mixing together'' elements of $G^{(n)}$. If $N$ is an 
integer and $(N_1,\ldots,N_n)$
is a partition of $N$, let us define $G^{(N)}(N_1,\ldots,N_n)$
as the direct sum of all $G^{(N)}_\up$, where $\up$ is such that 
for any $k = 1,\ldots,n$ and any pair $i,j$ of elements of 
$N_1 + \cdots + N_{k-1} + \{1,\ldots,N_k\}$, $p_{ij} = 0$. 
Assume that $(N^{(\beta)})_{\beta = 1,\ldots,\al}$ is a 
family of integers and for each $\beta = 1,\ldots,\al$, 
$(N^{(\beta)}_\al)$ is a partition of $N^{(\beta)}$. 
Assume that for $\beta = 1,\ldots,\al$, 
$x_\beta\in G^{(N^{(\beta)})}(N^{(\beta)}_1,\ldots,N^{(\beta)}_n)$,  
and let us fix a sequence 
of bijections $b_k : \{1,\ldots, \sum_\beta N^{(\beta)}_k\} 
\to I_{N^{(1)}_k,\ldots,N^{(\beta)}_k}$. To the sequence
$\ub = (b_k)_{k = 1,\ldots,n}$, we associate a map 
$$
\conc_{\ub} : \bigotimes_{\beta = 1}^{\al} 
G^{(N^{(\beta)})}(N^{(\beta)}_1,\ldots,N^{(\beta)}_N)
\to G^{(n)}, 
$$
such that if $(b'_k,b''_k)$ are the components of the map 
$b_k$, then the first factor of $\conc_\ub(
\otimes_{\beta = 1}^{\al} x_\beta)$ is the concatenation of
the $b''_1(1)$th factor of $x_{b'(1)}$, the 
$b''_1(2)$th factor of $x_{b'(2)}$, etc, the $i$th 
factor of  $\conc_\ub(
\otimes_{\beta = 1}^{\al} x_\beta)$ is the concatenation of
$N^{(b'(i))}_1 + \cdots + N^{(b'(i))}_{i-1} + 
b''_i(1)$th factor of $x_{b'(i)}$, the 
$N^{(b'(i))}_1 + \cdots + N^{(b'(i))}_{i-1} + 
b''_i(2)$th factor of $x_{b'(i)}$, etc. 

Then we have 
\begin{prop} \label{prop:normal:ordering} 
 $\mu^{(n)}\circ \conc_\ub  = \mu^{(n)}\circ \conc_\ub\circ 
(\otimes_{\beta = 1}^\al \mu^{(N^{(\beta)})})$. 
\end{prop}

{\em Proof.} We have seen that the image by $\mu^{(3)}$
of an element of $G^{(3)}$ is not changed if we apply 
the CYBE identity within in this element. Therefore the 
same is true for any map $\mu^{(k)}$. Now each map  
$\mu^{(N^{(\al)})}$ consists of repeated applications of the 
CYBE identity, so the image by $\mu^{(N)}$ of 
$\conc_\ub(\otimes_{\beta = 1}^\al  x_\beta)$ and $\conc_\ub 
(\otimes_{\beta = 1}^\al \mu^{(N^{(\beta)})}(x_\beta))$ will 
be the same. 
\hfill \qed\medskip 

Using $\mu^{(n)}$, one may prove in the same way as Corollary 
\ref{cor:Lie}: 

\begin{prop}  \label{prop:Lie:n}
Assume that for $i = 1,\ldots,n$, $L_i$ is an element of 
$FL_{\sum_{j|j<i} p_{ji} + \sum_{j|j>i} p_{ij}}$. 
Then the element of $F^{(n)}$ equal to 
$$
\mu^{(n)} \big(\bigotimes_{i = 1}^n L_i(u^{(i)}_1,\ldots,
u^{(i)}_{\sum_{j|j<i} p_{ji} + \sum_{j|j>i} p_{ij}} ) \big) 
$$
belongs to $\bigoplus_{x_2,\ldots,x_{n-2}\in \{a,b\}} F^{(ax_2
\cdots x_{n-2}b)}$. 
\end{prop}

As it happens in the case $n = 3$,  
the restriction on $\up$ (see equation (\ref{def:lie:space})) 
comes from the fact that the reorderings due to CYBE 
tranfer the $x^{(i)}_\al$ to the left and the $y^{(i)}_\al$ to the
right, so the fact that each $x^{(i)}_\al$ is  at the left of the
corresponding $y^{(i)}_\al$ is not changed.

\subsection{Construction and properties of $\delta_4^{(F)}$}

\subsubsection{The maps $x\mapsto x^{(i_1\ldots i_p)}$}

Let $i_1,\ldots,i_p$ be integers such that $1\leq i_1<\cdots < i_p\leq n$. 
Then there is a unique linear map $\us\mapsto \us^{(i_1\ldots i_p)}$
from $\NN^p$ to $\NN^n$, such that $\us^{(i_1\ldots i_p)}_i = 0$
if $i\notin\{i_1,\ldots,i_p\}$ and $\us^{(i_1\ldots i_p)}_{i_j} = 
s_j$ for $j = 1,\ldots,p$. If $(\us,\ut,\phi)$ belongs to $\Phi_p$, 
define $\phi^{(i_1 \ldots i_p)}$ as the map from 
$I_{\us^{(i_1\ldots i_p)}}$ to $I_{\ut^{(i_1\ldots i_p)}}$
such that  if $(k,\al)\in I_{\us}$ and $(\phi(k,\al)) = 
(l,\beta)$, then $\phi^{(i_1\ldots i_p)}(i_k,\al) = (i_l,\beta)$. 
Denote also by $x\mapsto x^{(i_1\ldots i_p)}$ the map from 
$\Phi_p$ to $\Phi_n$, such that $(\us,\ut,\phi)^{(i_1\ldots i_p)} 
= (\us^{(i_1\ldots i_p)},\ut^{(i_1\ldots i_p)},
\phi^{(i_1\ldots i_p)})$. Let us also denote by $x\mapsto 
x^{(i_1 \ldots x_p)}$ the linear map from $F^{(p)}$ to $F^{(n)}$, 
such that for any $(\us,\ut,\phi)\in\Phi_p$, 
$z(\us,\ut,\phi)^{(i_1\ldots i_p)} = z((\us,\ut,\phi)^{(i_1\ldots i_p)})$. 
Then   $x\mapsto 
x^{(i_1 \ldots x_p)}$ is an algebra morphism from 
$(F^{(p)},m_{F^{(p)}})$ to $(F^{(n)},m_{F^{(n)}})$.

\subsubsection{}

Let us denote by $r$ the element of $F^{(2)}_1$ equal to 
$x^{(12)}_1\otimes y^{(12)}_1$.

\begin{lemma} \label{entretien}
1) Let $\rho$ belong to $\bigoplus_{x_1,\ldots,x_{n-3}\in \{a,b\}}
F^{(ax_1\ldots x_{n-3}b)}$. Then if $i<j$ (resp., if $j<i$), 
the commutators
$[r^{(ij)},\rho^{1\ldots i-1,i+1\ldots n}]$  
(resp., $[r^{(ji)},\rho^{1\ldots i-1,i+1\ldots n}]$)
belongs to $\bigoplus_{x_1,\ldots, x_{n-3}\in \{a,b\}}
F^{(ax_1\ldots x_{n-2}b)}$. 

2) If $i<j<k$, we have $[r^{(ij)},r^{(ik)}] + [r^{(ij)},r^{(jk)}] 
+ [r^{(ik)},r^{(jk)}] = 0$ and if $i,j,k,l$ are all distinct and 
$i<j$, and $k<l$, then $[r^{(ij)},r^{(kl)}] = 0$. 
\end{lemma}

{\em Proof.} It suffices to prove 1) when $x$ belongs to 
some $F^{(ax_1\ldots x_{n-3}b)}$ and is the tensor
product of Lie polynomials. There is nothing to 
do when $i<j$ and $x_{j-1} = b$, or when $i>j$ and $x_j = a$. 
In the two other cases, one applies   
Proposition \ref{prop:Lie:n} to a suitable family 
of Lie polynomials. 

The first part of 2) is a consequence of the fact that 
the map $x\mapsto x^{(ijk)}$ from $F^{(3)}$ to $F^{(n)}$
is an algebra morphism. The second part is immediate. 
\hfill \qed\medskip 

When $x\in F^{(3)}$ define $\delta^{F^{(3)} \to F^{(4)}}(x)$
as follows 
\begin{align*}
&
\delta^{F^{(3)}\to F^{(4)}}(x) = [r^{(12)} + r^{(13)} + r^{(14)}, x^{(234)}] 
-  [  - r^{(12)} + r^{(23)} + r^{(24)}, x^{(134)}]
\\ & 
+ [  - r^{(13)} - r^{(23)} + r^{(34)}, x^{(124)}] 
-  [  - r^{(14)} - r^{(24)} - r^{(34)}, x^{(123)}] . 
\end{align*}  
$\delta^{F^{(3)} \to F^{(4)}}$ maps $F^{(3)}$ to $F^{(4)}$. 
It follows from Lemma \ref{entretien}, 1) that if 
$x\in F^{(aab)}\oplus F^{(abb)}$,  $\delta^{F^{(3)} \to 
F^{(4)}}$ maps 
maps $F^{(aab)}\oplus F^{(abb)}$ to 
$\oplus_{x,y\in \{a,b\}} F^{(axyb)}$. We define $\delta_4^{(F)}$
as the resulting map  from $F^{(aab)}\oplus F^{(abb)}$ to 
$\oplus_{x,y\in \{a,b\}} F^{(axyb)}$. 

Then it follows from Proposition \ref{prop:universal:fn} that 
$\delta_4^{(F)}$ satisfies the conclusion of \ref{prop:delta:4}. 

Recall that $\delta_3^{(F)}$ is the restriction to $\prod_{n\geq 0}
F_n$ of the map $\delta^{F^{(2)}\to F^{(3)}} : 
F^{(2)} \to F^{(3)}$ such that 
$$
\delta^{F^{(2)}\to F^{(3)}}(x) = 
[r^{(12)},x^{(13)}]  + [r^{(12)},x^{(23)}]  
+ [r^{(13)},x^{(23)}]  + 
[x^{(12)},r^{(13)}]  + [x^{(12)},r^{(23)}]  
+ [x^{(13)},r^{(23)}]  . 
$$
Then it follows from Lemma \ref{entretien}, 2) that 
the composition $\delta^{F^{(3)}\to F^{(4)}} 
\circ \delta^{F^{(2)}\to F^{(3)}}$, 
which is a map from $F^{(2)}$ to $F^{(4)}$, is zero. 
It follows that the composition $\delta_4^{(F)}\circ \delta_3^{(F)}$
is also zero.   

This proves Proposition \ref{prop:delta:diff}, 2). 

\begin{remark} It is now clear how to define $\delta^{(F)}_n
: \oplus_{x_1,\ldots,x_{n-3}\in\{a,b\}} F^{(a x_1\cdots x_{n-3}b)} 
\to \oplus_{y_1,\ldots,y_{n-2}\in\{a,b\}}F^{(ay_1\cdots y_{n-2}b)}$, 
such that $\delta^{(F)}_n \circ \delta^{(F)}_{n-1} = 0$. 
\end{remark}

\newpage
\section{Computation of cohomology groups} 
\label{app:cohomologies}

\subsection{Computation of $H^2_n$} \label{sect:H2}

Clearly, $\Ker(\delta_3^{(F)}) \cap F_1 = \KK r$. Let us assume 
now that $n>1$. We want to prove that $\Ker(\delta_3^{(F)})\cap F_n
= 0$. 

Let $x$ belong to  $\Ker(\delta_3^{(F)})\cap F_n$. Since the 
family $(x^{(12)}_1\cdots x^{(12)}_n \otimes 
y^{(12)}_{\sigma(1)} \cdots y^{(12)}_{\sigma(n)})_{\sigma\in\SG_n}$
is a basis of $F^{(2)}_n = (FA_n\otimes FA_n)_{\SG_n}$, 
there exists a unique family $(A_\sigma)_{\sigma\in \SG_n}$
in $\KK^{\SG_n}$ such that 
$$
x = \sum_{\sigma\in\SG_n} A_\sigma
x^{(12)}_1\cdots x^{(12)}_n \otimes 
y^{(12)}_{\sigma(1)} \cdots y^{(12)}_{\sigma(n)} . 
$$  
Moreover: 

\begin{lemma} \label{observation} 
The condition that $x$ belong to $(FL_n\otimes FL_n)_{\SG_n}$
is equivalent to the condition that 
$\sum_{\sigma\in \SG_n} A_\sigma 
x_{\sigma(1)}\cdots x_{\sigma(n)}$  
and  $\sum_{\sigma\in \SG_n} A_{\sigma^{-1}} 
x_{\sigma(1)}\cdots x_{\sigma(n)}$ are both Lie polynomials in 
the free algebra with generators $x_1,\ldots,x_n$. 
\end{lemma}

Let $\pi_{01n}$ denote the projection of $F^{(3)}$ on 
$F_{01n}$ parallel to $\bigoplus_{(p,q,r)|(p,q,r)\neq (0,1,n)}
F_{pqr}$. Let us apply $\pi_{01n}$ to the identity 
$\delta_3^{(F)}(x) = 0$.  

Since $[x^{12},r^{13}]\in F_{0n1}$, 
$[r^{13},x^{23}]\in F_{n01}$ and $[x^{13},r^{23}]\in F_{10n}$ 
and since $n\neq 1$, 
$\pi_{01n}(\delta_3^{(F)}(x))$ 
is the same as the image by $\pi_{01n}$ of 
$[r^{12},x^{23}] + [x^{12},r^{13}] + [r^{12},x^{13}]$. 
This is also the image by $\pi_{01n}$ of 
\begin{align} \label{pobeda}
&  \sum_{\sigma\in\SG_n} A_\sigma \mu^{(3)} 
( x^{(2)}_1 \otimes y^{(2)}_1 x^{(1)}_{\sigma(1)} 
\cdots x^{(1)}_{\sigma(n)} \otimes y^{(1)}_1 \cdots y^{(1)}_n) 
\\ & \nonumber 
+ \sum_{\sigma\in\SG_n} A_\sigma \mu^{(3)} 
( x^{(2)}_{\sigma(1)}\cdots x^{(2)}_{\sigma(n)} 
\otimes y^{(2)}_1 \cdots y^{(2)}_n x^{(1)}_{1} 
\otimes y^{(1)}_1 ) 
\\ & \nonumber 
+ \sum_{\sigma\in \SG_n} A_\sigma 
[x^{(2)}_1,
x^{(3)}_{\sigma(1)} \cdots x^{(3)}_{\sigma(n)}]\otimes y^{(2)}_1 
\otimes y^{(3)}_1\cdots y^{(3)}_n , 
\end{align} 
where the argument of the first (resp., second) expression 
in $\mu^{(3)}$ belongs to $G^{(3)}_{n10}$  (resp., $G^{(3)}_{1n0}$). 

The image of (\ref{pobeda}) by $\pi_{01n}$ is then equal to 
\begin{align} \label{cassirer} 
& \sum_{\sigma\in\SG_n} A_\sigma [x^{(3)}_{\sigma(1)},\cdots
,[x^{(3)}_{\sigma(n)},x^{(2)}_1]] \otimes y^{(2)}_1\otimes
y^{(3)}_1 \cdots y^{(3)}_n
\\ & \nonumber 
+  \sum_{\sigma\in\SG_n} A_\sigma x^{(3)}_{\sigma(1)} 
\cdots [x^{(3)}_1,x^{(2)}_1]\cdots x^{(3)}_{\sigma(n)}
\otimes y^{(2)}_1 \otimes [[y^{(3)}_1,y^{(3)}_2],\cdots,y^{(3)}_n]
\\ & \nonumber
+ \sum_{\sigma\in\SG_n} A_\sigma [x^{(2)}_1,
x^{(3)}_{\sigma(1)} \cdots x^{(3)}_{\sigma(n)}]\otimes y^{(2)}_1 
\otimes y^{(3)}_1\cdots y^{(3)}_n . 
\end{align}

\begin{lemma} \label{chrono}
Let $(X_\sigma)_{\sigma\in\SG_n}$ belong to $\KK^{\SG_n}$ and 
assume that  $X = \sum_{\sigma\in\SG_n} X_\sigma x_{\sigma(1)}
\cdots x_{\sigma(n)}$ is a Lie element of the free algebra
with generators $x,x_1,\dots,x_n$. Then we have the identity
\begin{equation} \label{basic:Lie}
[X,x] = \sum_{\sigma\in\SG_n} X_\sigma [x_{\sigma(1)},\cdots,
[x_{\sigma(n)},x]] . 
\end{equation}  
\end{lemma}

{\em Proof of Lemma.}
Let us denote by $\ad$ the adjoint action of the free algebra with 
$n+1$ generators on itself by the rule $\ad(a)(b) = [a,b]$. 
Then we should show the identity $[X,x] = \ad(X)(x)$. 
$X$ is a linear combination of Lie polynomials of the form 
$[x_{\tau(1)},\cdots,[x_{\sigma(n-1)},x_{\sigma(n)}]]$
(see e.g.\ \cite{Bbk}), so it suffices to check (\ref{basic:Lie})
when $X$ is such a polynomial. Moreover, we may assume after relabeling 
indices that $\tau$ is the inversion $\tau(i) = n+1 - i$. 
We have now to prove that if $X_n = [x_n,\cdots,[x_2,x_1]]$, then 
\begin{equation} \label{basic:Lie:n}
[X_n,x] = \ad(X_n)(x). 
\end{equation} 
Let us prove (\ref{basic:Lie:n}) by induction. The case $n = 1$
is obvious.   Assume that (\ref{basic:Lie:n}) holds at step $n$. 
Then  
\begin{align*}
& \ad(X_n)(x) = \ad([x_{n+1},X_n])(x)  
= (\ad(x_{n+1})\ad(X_n) - \ad(X_n)\ad(x_{n+1}))(x)
\\ & 
= [x_{n+1},[X_n,x]] - [X_n,[x_{n+1},x]] = [X_{n+1},x] , 
\end{align*}
where the first equality of the second line uses the twice the 
induction hypothesis, first with $x$, then with $x$ replaced 
by $[x_{n+1},x]$. 
\hfill \qed\medskip 

Lemma \ref{chrono} implies that the first and last terms of 
(\ref{cassirer}) cancel out, so that 
\begin{align} \label{cass:2}
& \sum_{\sigma\in\SG_n}
A_\sigma x^{(3)}_{\sigma(1)} \cdots
x^{(3)}_{\sigma(\sigma^{-1}(1) - 1)} 
\big( x^{(3)}_1 x^{(2)}_1 - x^{(2)}_1 x^{(3)}_1 \big) 
x^{(3)}_{\sigma(\sigma^{-1}(1) + 1)}\cdots x^{(3)}_{\sigma(n)}
\\ & \nonumber 
\otimes y^{(2)}_1 \otimes [[y^{(3)}_1,y^{(3)}_2],\cdots,y^{(3)}_n]
\end{align}
is zero. Decompose (\ref{cass:2}) as a sum $\sum_{i = 1}^{n+1} A_i$, 
where $A_i$ contains $x^{(2)}_1$ at the $i$th position. Then 
each $A_i$ is zero. Replacing $x^{(2)}_1$ by $1$ in each of these 
equalities, we find  the following set of equalities in 
$(FA_n\otimes FA_n)_{\SG_n}$
$$
\sum_{\sigma\in \SG_n|\sigma(1) = 1}
A_\sigma x^{(3)}_{\sigma(1)} \cdots x^{(3)}_{\sigma(n)}
\otimes [[y^{(3)}_1,y^{(3)}_2],\cdots,y^{(3)}_n] = 0, 
$$
\begin{align*}
& \sum_{\sigma\in \SG_n|\sigma(k) = 1}
A_\sigma x^{(3)}_{\sigma(1)} \cdots x^{(3)}_{\sigma(n)}
\otimes [[y^{(3)}_1,y^{(3)}_2],\cdots,y^{(3)}_n] 
\\ & =  
\sum_{\sigma\in \SG_n|\sigma(k+1) = 1}
A_\sigma x^{(3)}_{\sigma(1)} \cdots x^{(3)}_{\sigma(n)}
\otimes [[y^{(3)}_1,y^{(3)}_2],\cdots,y^{(3)}_n] , 
\end{align*}
for $k = 1,\ldots,n-1$
$$
\sum_{\sigma\in \SG_n|\sigma(n) = 1}
A_\sigma x^{(3)}_{\sigma(1)} \cdots x^{(3)}_{\sigma(n)}
\otimes [[y^{(3)}_1,y^{(3)}_2],\cdots,y^{(3)}_n] = 0  . 
$$
We have therefore for any $k = 1,\ldots,n$
$$
\sum_{\sigma\in \SG_n|\sigma(k) = 1}
A_\sigma x^{(3)}_{\sigma(1)} \cdots x^{(3)}_{\sigma(n)}
\otimes [[y^{(3)}_1,y^{(3)}_2],\cdots,y^{(3)}_n] = 0  . 
$$
Adding up these equalities, and using the identity  
$x^{(3)}_{\sigma(1)}\cdots x^{(3)}_{\sigma(n)}\otimes 
y^{(1)}_1\cdots y^{(1)}_n =  
x^{(3)}_{1}\cdots x^{(3)}_{n}\otimes 
y^{(1)}_{\sigma^{-1}(1)}\cdots y^{(1)}_{\sigma^{-1}(n)}$ 
in $(FA_n\otimes FA_n)_{\SG_n}$ we get 
\begin{equation} \label{cass:3}
x^{(3)}_1\cdots x^{(3)}_n\otimes \sum_{\sigma\in\SG_n}
A_\sigma [[y^{(3)}_{\sigma^{-1}(1)},y^{(3)}_{\sigma^{-1}(2)}],
\cdots, y^{(3)}_{\sigma^{-1}(n)}] = 0. 
\end{equation}

\begin{prop} \label{reut} (see \cite{Reut})
If $X = \sum_{\sigma\in\SG_n} X_\sigma x_{\sigma(1)}\cdots 
x_{\sigma(n)}$ is a Lie polynomial in the free algebra 
with generators $x_1,\ldots,x_n$, then 
$$
\sum_{\sigma\in\SG_n} X_\sigma [[x_{\sigma(1)},x_{\sigma(2)}],
\cdots,x_{\sigma(n)}] = n X . 
$$
\end{prop}

We have seen that $\sum_{\sigma\in\SG_n} 
A_\sigma x_{\sigma^{-1}(1)} \cdots x_{\sigma^{-1}(n)}$
is a Lie polynomial, therefore (\ref{cass:3}) is equal to 
$$
n \cdot x^{(3)}_1\cdots x^{(3)}_n \otimes
\sum_{\sigma\in\SG_n} A_\sigma y^{(3)}_{\sigma^{-1}(1)}  
\cdots y^{(3)}_{\sigma^{-1}(n)} = nx .  
$$
Therefore $x = 0$. 

\subsection{Computation of $H^3_n$}

\subsubsection{Form of the elements of $\Imm(\delta_3^{(F)})$}

\begin{prop} \label{prereut}
Let $(A_\sigma)_{\sigma\in\SG_n}\in \KK^{\SG_n}$ be such that
$\sum_{\sigma\in\SG_n} A_\sigma x_{\sigma(1)}\cdots x_{\sigma(n)}$
is a Lie polynomial. Then for any $k = 1,\ldots,n$, we have
$$
\sum_{\sigma\in\SG_n} A_\sigma x_{\sigma(1)}\cdots x_{\sigma(n)}
=  \sum_{\sigma\in\SG_n} 
A_\sigma [x_{\sigma(1)},\cdots [x_{\sigma(n-1)}, x_{\sigma(n)}]]. 
$$  
\end{prop}

{\em Proof.} We may assume that 
$\sum_{\sigma\in\SG_n} A_\sigma x_{\sigma(1)}\cdots x_{\sigma(n)}
= [x_1,\ldots, [x_{n-1},x_n]]$. Let us then prove the result by induction
on $n$. Assume that we proved the result up to order $n-1$
and let us treat the case of order $n$. Let us define  
$(A_\sigma(n))_{\sigma\in\SG_n}$ as the elements of $\KK$ such that 
$$
\sum_{\sigma\in\SG_n} A_\sigma(n) x_{\sigma(1)} \cdots x_{\sigma(n)}
= [x_1,\ldots, [x_{n-1},x_n]] . 
$$
Then if $k = 2,\ldots,n$ and if $\sigma\in\SG_n$ is such that 
$\sigma(n) = k$ and $A_\sigma(n)\neq 0$, then $\sigma(1) = 1$. 

For any $\sigma\in\SG_n$ such that $\sigma(1) = 1$, let us denote by 
$\sigma'$ the element of $\SG_{n-1}$ such that $\sigma'(k) = 
\sigma(k+1) - 1$, for $k = 1,\ldots,n-1$. Then when $\sigma\in \SG_n$
and $\sigma(1) = 1$, we have $A_\sigma(n) = A_{\sigma'}(n-1)$. 

Then if $k = 2,\ldots,n$, 
\begin{align} \label{partial}
& \sum_{\sigma\in\SG_n|\sigma(n) = k} A_\sigma(n) [x_{\sigma(1)},
\ldots, [x_{\sigma(n-1)},x_{\sigma(n)}]]
\\ & \nonumber 
= \sum_{\tau\in\SG_{n-1}} A_\tau(n-1) [x_{1},
[x_{\sigma'(1)+1},\ldots, [x_{\sigma'(n-2)+1},x_{\sigma'(n-1)+1}]]
\\ & \nonumber 
= [x_1,[x_2,\ldots,[x_{n-1},x_n]]] ,  
\end{align}  
where the second equality follows from the induction 
hypothesis applied to variables $x_2,\ldots,x_n$. 
This proves the result at order $n$, when $k = 2,\ldots,n$. 
Then 
\begin{align*}
& 
\sum_{\sigma\in\SG_n|\sigma(n) = 1} A_\sigma(n)
[x_{\sigma(1)},\ldots,[x_{\sigma(n-1)},x_{\sigma(n)}]]
\\ & 
= \sum_{\sigma\in\SG_n} A_\sigma(n)
[x_{\sigma(1)},\ldots,[x_{\sigma(n-1)},x_{\sigma(n)}]]
- \sum_{k = 2}^n \sum_{\sigma\in\SG_n|\sigma(n) = k} A_\sigma(n)
[x_{\sigma(1)},\ldots,[x_{\sigma(n-1)},x_{\sigma(n)}]]
\\ & 
= 
n \sum_{\sigma\in\SG_n} A_\sigma(n) x_{\sigma(1)}\cdots x_{\sigma(n)} 
- (n-1) \sum_{\sigma\in\SG_n} A_\sigma(n) x_{\sigma(1)}\cdots x_{\sigma(n)} , 
\end{align*}
where the second equality 
follows from Proposition \ref{reut} and from equalities 
(\ref{partial}). 
So 
$$
\sum_{\sigma\in\SG_n|\sigma(n) = 1} A_\sigma(n)
[x_{\sigma(1)},\ldots,[x_{\sigma(n-1)},x_{\sigma(n)}]]
= \sum_{\sigma\in\SG_n} A_\sigma(n) x_{\sigma(1)}\cdots x_{\sigma(n)} 
= [x_1,\ldots,[x_{n-1},x_n]], 
$$
which proves the remaining case of the result at order $n$. 
\hfill \qed\medskip

\begin{prop} \label{thl}
Let $n>1$ and $y\in F_n$. Set 
$x = \delta_3^{(F)}(y)$, and let us decompose $x$ as 
$\sum_{p,q|p+q = n+1} x^{(aab)}_{p,q,n+1} 
+ \sum_{p,q|p+q = n+1} x^{(abb)}_{n+1,p,q}$, with  
$$
x^{(aab)}_{p,q,n+1} \in (FL_p\otimes FL_q\otimes FL_{n+1})_{\SG_p\times\SG_q}
\quad 
\on{and}\quad 
x^{(abb)}_{n+1,p,q} \in 
(FL_{n+1} \otimes FL_p\otimes FL_q)_{\SG_p\times\SG_q} .  
$$
Then 

1) $x^{(aab)}_{n,1,n+1} = [r^{(23)},y^{(13)}]$ and
$x^{(abb)}_{n+1,1,n} = - [r^{(12)},y^{(13)}]$. 
 
2)  $x^{(aab)}_{1,n,n+1} = 0$ and 
$x^{(abb)}_{n+1,1,n} = 0$.  
\end{prop}  

{\em Proof.} Proposition 
 \ref{prereut} allows some simplifications in the computations of 
the end of  Section \ref{sect:H2}, which imply the two first 
equalities. To prove the two last equalities, let us proceed as 
in Section \ref{sect:H2}. For example, if $y = 
\sum_{\sigma\in\SG_n} A_\sigma x^{(12)}_1\cdots x^{(12)}_n\otimes 
y^{(12)}_{\sigma(1)} \cdots y^{(12)}_{\sigma(n)}$, then 
the nonzero contributions to 
$x^{(aab)}_{1,n,n+1}$ are those of $[r^{(12)},y^{(23)}]$ 
and of $[r^{(13)},y^{(23)}]$, which are respectively 
$$
\sum_{k = 1}^n\sum_{\sigma\in \SG_n} A_\sigma
x^{(13)}_1\cdots x^{(13)}_n\otimes x^{(23)}
y^{(13)}_{\sigma(1)} \cdots 
[y^{(13)}_{\sigma(k)}, y^{(23)}]
\cdots y^{(13)}_{\sigma(n)} 
$$
and
$$
\sum_{k = 1}^n\sum_{\sigma\in \SG_n} A_\sigma
x^{(13)}_1\cdots x^{(13)}_n\otimes x^{(23)}
[y^{(23)} , y^{(13)}_{\sigma(1)} \cdots 
y^{(13)}_{\sigma(n)}] 
$$
and cancel out (here we do not use the fact that $y$ belongs to 
$(FL_n\otimes FL_n)_{\SG_n}$). 
\hfill \qed\medskip

\subsubsection{Computation of $H^3_n$}

The result is obvious in the case $n = 2$. 

Let us assume that $n>2$ and let $x$ belong to $F^{Lie,(3)}_{n}$
be such that $\delta^{(F)}_4(x) = 0$. We want to show the 
existence of $y\in F^{Lie,(2)}_{n-1}$ such that 
$x = \delta^{(F)}_3(y)$. 

\begin{prop}
We have $x^{(aab)}_{1,n-1,n} = 0$ and $x^{(abb)}_{n,n-1,1} = 0$. 
\end{prop}

{\em Proof.} When $z\in F^{Lie,(4)}_n$, $x,y\in\{a,b\}$
and $p,q,r,s$ are integers, we denote by $z^{(axyb)}_{pqrs}$
the projection of $z$ on $F^{(axyb)}_{pqrs}$ parallel to the
direct sum of all other $F^{(abx'y'b)}_{p'q'r's'}$. 
Then $(\delta_4^{(F)}(x))^{(aaab)}_{1,1,n-1,n+1}$ is equal to 
$$
\big( [r^{(12)} + r^{(13)} + r^{(14)}, (x^{(aab)}_{1,n-1,n})^{(234)} 
] - [r^{(23)} + r^{(24)} , 
(x^{(aab)}_{1,n-1,n})^{(134)} ] \big)^{(aaab)}_{1,1,n-1,n+1} , 
$$
because all other summands of $\delta_4^{(F)}(x)$ project to 
zero on  $F^{(aaab)}_{1,1,n-1,n+1}$, either obviously, or in 
the case of $[r^{(34)},(x^{(aab)}_{112})^{(123)} + 
(x^{(aab)}_{112})^{(124)}]$, because we assumed that $n>2$. 

\begin{lemma}
$$
\big( [r^{(12)} + r^{(13)} + r^{(14)}, (x^{(aab)}_{1,n-1,n})^{(234)} 
] \big)^{(aaab)}_{1,1,n-1,n+1} = 0
$$
\end{lemma}

{\em Proof of Lemma.} Let $(Z_{k,\sigma})_{(k,\sigma)\in\{0,\ldots,n-1\}
\times\SG_{n-1}}$
be the elements of $\KK$ such that 
$$
x^{(aab)}_{1,n-1,n} = \sum_{k = 0}^{n-1} \sum_{\sigma\in \SG_{n-1}}
Z_{k,\sigma} x_1^{(3)} \otimes x_1^{(1)} \cdots x_{n-1}^{(1)}
\otimes y^{(1)}_{\sigma(1)} \cdots y^{(1)}_{\sigma(k)}
y^{(3)}_1 y^{(1)}_{\sigma(k+1)} \cdots y^{(1)}_{\sigma(n-1)} . 
$$ 
Then 
\begin{align*}
& [r^{(12)},(x^{(aab)}_{1,n-1,n})^{(234)}]
\\ & 
= \sum_{k = 0}^{n-1} \sum_{\sigma\in\SG_{n-1}}
Z_{k,\sigma} x_1^{(14)} \otimes x^{(24)}_1
\otimes x^{(34)}_1 \cdots x^{(34)}_{n-1} 
\\ & 
\otimes 
y^{(34)}_{\sigma(1)} \cdots y^{(34)}_{\sigma(k)} 
[y^{(24)}_1,y^{(14)}_1]
y^{(34)}_{\sigma(k+1)} \cdots y^{(34)}_{\sigma(n-1)} 
\end{align*}
and
\begin{align*}
& [r^{(13)},(x^{(aab)}_{1,n-1,n})^{(234)}]
\\ & 
= \sum_{k = 0}^{n-1} \sum_{\sigma\in\SG_{n-1}}
Z_{k,\sigma} \sum_{s = 1}^{n-1}
x_1^{(14)} \otimes x^{(24)}_1
\otimes x^{(34)}_1 \cdots x^{(34)}_{n-1} 
\\ & 
\otimes 
y^{(34)}_{\sigma(1)} \cdots [y^{(34)}_{\sigma(s)},y^{(14)}_{1}]
\cdots y^{(34)}_{\sigma(k)}  y^{(24)}_1
y^{(34)}_{\sigma(k+1)} \cdots y^{(34)}_{\sigma(n-1)} ,  
\end{align*}
so that the sum of these terms is equal to $ - [r^{(14)},
(x^{(aab)}_{1,n-1,n})^{(234)}]$. 
\hfill \qed\medskip 

{\em End of proof of Proposition.}
In the same way, one proves 
\begin{align*}
& \big( [r^{(23)} + r^{(24)},(x^{(aab)}_{1,n-1,n})^{(234)}]
\big)^{(aaab)}_{1,1,n-1,n+1}
\\ & 
= \sum_{k = 0}^{n-1}\sum_{\sigma\in\SG_{n-1}}
Z_{k,\sigma} x^{(14)}_1\otimes x^{(24)}_1\otimes 
x^{(34)}_1 \cdots x^{(34)}_{n-1}
\\ & \otimes 
y^{(34)}_{\sigma(1)} \cdots y^{(34)}_{\sigma(k)}
[y_1^{(14)},y^{(24)}_1]
y^{(34)}_{\sigma(k+1)} \cdots y^{(34)}_{\sigma(n-1)} . 
\end{align*}
Therefore $\sum_{k = 0}^{n-1}\sum_{\sigma\in\SG_{n-1}}
Z_{k,\sigma} y^{(34)}_{\sigma(1)} \cdots y^{(34)}_{\sigma(k)}
[y_1^{(14)},y^{(24)}_1]
y^{(34)}_{\sigma(k+1)} \cdots y^{(34)}_{\sigma(n-1)}  = 0$. 
For any $k$, the sum of all the terms in this sum in which $y_1^{(14)}$ 
and $y_1^{(24)}$ appear in the $k$th and $(k+1)$st position 
is also zero, so for any $k$ we have 
$\sum_{\sigma\in\SG_{n-1}} Z_{k,\sigma} y^{(34)}_{\sigma(1)} 
\cdots y^{(34)}_{\sigma(n-1)} = 0$. So the $Z_{k,\sigma}$ are
all zero and $x^{(aab)}_{1,n-1,n} = 0$. 

The proof of $x^{(abb)}_{n,n-1,1} = 0$ is similar. 
\hfill \qed\medskip 

Recall that $F^{Lie,(n)}$ is the direct sum 
$\bigoplus_{x_1,\ldots,x_{n-2}\in\{a,b\}} 
F^{(ax_1 \ldots x_{n-2}b)}$. If $z\in F^{Lie,(n)}$, let us denote by 
$x^{(ax_1\ldots x_{n-2}b)}$ the projection of $z$ on 
$F^{(ax_1 \ldots x_{n-2}b)}$ parallel to the direct sum of all 
other $F^{(ax'_1 \ldots x'_{n-2}b)}$. 

\begin{prop} \label{garges}
If $z$ belongs to $\oplus_{x\in \{a,b\}} F^{(axb)}_n$, then 
$$
(\delta_4(z))^{(aabb)} = 
- [r^{(13)} , (z^{(aab)})^{(124)} + (z^{(aab)})^{(214)}] 
- [r^{(24)}, (z^{(abb)})^{(134)} + (z^{(abb)})^{(314)}] . 
$$
\end{prop}

{\em Proof.} If $n$ is an integer $\geq 2$, 
$x_1,\ldots,x_{n-3}$ belong to $\{a,b\}$, and if
$1\leq i<j\leq n$, then 
$$
[r^{(ij)}, (F^{(ax_1\cdots x_{n-3}b)})^{1,\ldots,i-1,i+1,\ldots,n} 
]\subset 
\oplus_{x\in\{a,b\}} F^{(ax_1\cdots x_{i-2}a x_{i-1}\cdots x_{j-3}x x_{j-2}\cdots x_{n-3}b)}
$$
and 
$$
[r^{(ij)}, (F^{(ax_1\cdots x_{n-3}b)})^{1,\ldots,j-1,j+1,\ldots,n}  ]\subset 
\oplus_{x\in\{a,b\}} F^{(ax_1\cdots x_{i-2}a x_{i-1}\cdots x_{j-3}x x_{j-2}\cdots x_{n-3}b)}. 
$$
This implies that 
$(\delta^{(F)}_4(z))^{(aabb)}$ is equal to 
\begin{align} \label{expr:delta:4}
& \Big( [r^{(12)} + r^{(13)} + r^{(14)}, (z^{(abb)})^{(234)}]
+ [r^{(13)},(z^{(aab)})^{(234)}] 
- [r^{(23)},(z^{(aab)})^{(134)}] 
\\ & \nonumber 
- [r^{(23)} + r^{(24)},(z^{(abb)})^{(134)}] 
- [r^{(13)} + r^{(23)},(z^{(aab)})^{(124)}] 
- [r^{(23)},(z^{(abb)})^{124}]
\\ & \nonumber 
+ [r^{(14)} + r^{(24)} + r^{(34)}, (z^{(abb)})^{(123)}]
+ [r^{(14)},(z^{(abb)})^{(123)}]
\Big)^{(aabb)} . 
\end{align}
The reasoning of Proposition \ref{thl}, 2) implies that 
$$
\big( [r^{(12)} + r^{(13)} + r^{(14)}, (z^{(abb)})^{(234)}] \big)^{(aabb)}
\quad \on{and} \quad 
\big( [r^{(14)} + r^{(24)} + r^{(34)}, (z^{(abb)})^{(123)}]
\big)^{(aabb)}
$$
are zero (here we no not use the fact that the components of 
$z^{(aab)}$ and $z^{(abb)}$ are Lie polynomials)
and the reasoning of Proposition \ref{thl}, 1) relying on 
Proposition \ref{prereut} and the fact that the components of 
$z^{(aab)}$ and $z^{(abb)}$, shows that  
$$
([r^{(23)},(z^{(aab)})^{(134)}])^{(aabb)} 
= - [r^{(23)},(z^{(aab)})^{(124)}]
$$
and 
$$
([r^{(23)},(z^{(abb)})^{(124)}])^{(aabb)} 
= - [r^{(23)},(z^{(abb)})^{(134)}] . 
$$
Permuting the two first tensor factors of 
these relations, we find  
$$
([r^{(13)},(z^{(aab)})^{(134)}])^{(aabb)} 
= - [r^{(13)},(z^{(aab)})^{(214)}]
$$
and
$$
([r^{(23)},(z^{(abb)})^{(123)}])^{(aabb)} 
= - [r^{(23)},(z^{(abb)})^{(143)}] . 
$$
Substituting these expressions in (\ref{expr:delta:4}) gives 
the result. 
\hfill \qed\medskip 

\begin{remark} One proves in the same way that 
$(\delta^{(F)}_4(z))^{(abab)}$ is identically zero. 
\hfill \qed\medskip 
\end{remark}
 
\begin{cor} \label{9:ab}
Recall that $x$ belongs to $F^{Lie,(3)}_n$ and is such that 
$\delta^{(F)}_4(x) = 0$. Then there exists $y\in F_{n-1}$
such that 
$$
x^{(aab)}_{n-1,1,n} = [r^{(23)},y^{(13)}] 
\quad \on{and} \quad   
x^{(abb)}_{n,1,n-1} = - [r^{(12)},y^{(13)}] . 
$$
Moreover, for any integers $p,q$ such that $p+q = n$ and $p,q>1$, 
we have 
$$
x^{(aab)}_{p,q,n} +  (x^{(aab)}_{q,p,n})^{(213)}  = 0,  \quad  
x^{(abb)}_{n,p,q} +  (x^{(abb)}_{n,q,p})^{(132)} = 0 . 
$$
\end{cor}

{\em Proof.}
It follows from Proposition \ref{garges} that 
\begin{equation} \label{delta:0}
[r^{(13)} , (x^{(aab)})^{(124)} + (x^{(aab)})^{(214)}] +  
[r^{(24)}, (x^{(abb)})^{(134)} + (x^{(abb)})^{(143)}] = 0. 
\end{equation}
Let us project this equation on $F^{(aabb)}_{n,1,1,n}$ parallel to 
the sum of all other $F^{(aabb)}_{pqrs}$. Since $x^{(aab)}_{1,n-1,n} = 0$
and $x^{(abb)}_{n,n-1,1} = 0$, we get 
$$
[r^{(13)},(x^{(aab)}_{n-1,1,n})^{(124)}] +  
[r^{(24)},(x^{(abb)}_{n,1,n-1})^{(134)}] = 0. 
$$ 
Let us set 
$$
x^{(aab)}_{n-1,1,n} = \sum_{k = 0}^{n-1} \sum_{\sigma\in\SG_{n-1}}
A_{k,\sigma} x^{(3)}_1\cdots x^{(3)}_{n-1} \otimes x^{(1)}_1
\otimes y^{(3)}_{\sigma(1)}\cdots y^{(3)}_{\sigma(k)} 
y^{(1)}_1
y^{(3)}_{\sigma(k+1)}\cdots y^{(3)}_{\sigma(n-1)}
$$
and 
$$
x^{(abb)}_{n,1,n-1} = \sum_{k = 0}^{n-1} \sum_{\sigma\in\SG_{n-1}} 
B_{k,\sigma} 
x^{(3)}_{\sigma(1)}\cdots x^{(3)}_{\sigma(k)} 
x^{(2)}_1
x^{(3)}_{\sigma(k+1)}\cdots x^{(3)}_{\sigma(n-1)} 
\otimes y^{(2)}_1 \otimes 
y^{(3)}_1\cdots y^{(3)}_{n-1} . 
$$
Then we find 
\begin{align*}
& \sum_{k = 0}^{n-1} \sum_{\sigma\in\SG_{n-1}} 
A_{k,\sigma} [x^{(13)}_1,x^{(14)}_1\cdots x^{(14)}_{n-1}]
\otimes x^{(24)}_1\otimes y^{(13)}_1 \otimes y^{(14)}_{\sigma(1)}
\cdots y^{(14)}_{\sigma(k)}
y^{(24)}_1 y^{(14)}_{\sigma(k+1)}\cdots y^{(14)}_{\sigma(n-1)} 
\\ & 
+  \sum_{k = 0}^{n-1} \sum_{\sigma\in\SG_{n-1}}
x^{(14)}_{\sigma(1)} \cdots x^{(14)}_{\sigma(k)} 
x^{(13)}_1
x^{(14)}_{\sigma(k+1)} \cdots x^{(14)}_{\sigma(n-1)}
\otimes x^{(24)}_1 \otimes y^{(13)}_1 \otimes
[y^{(24)}_1,y^{(14)}_1 \cdots y^{(14)}_{n-1}] = 0 . 
\end{align*}
Identifying terms, we find that $A_{k,\sigma} = B_{k,\sigma} = 0$
whenever $k\notin\{0,n-1\}$. Moreover, for any $\sigma\in\SG_{n-1}$, 
we get  
$$
A_{0,\sigma} = - B_{0,\sigma^{-1}}, \
A_{n-1,\sigma} =  B_{0,\sigma^{-1}}, \
A_{0,\sigma} = B_{n-1,\sigma^{-1}}, \
- A_{n-1,\sigma} = - B_{n-1,\sigma^{-1}} . 
$$
Let us set $C_\sigma = A_{0,\sigma}$, then $A_{0,\sigma} = 
- A_{n-1,\sigma} = C_\sigma$ and $B_{0,\sigma^{-1}} = 
- B_{n-1,\sigma^{-1}} = - C_\sigma$, so if we set 
$$
y = \sum_{\sigma\in\SG_{n-1}} x^{(12)}_1\cdots x^{(12)}_{n-1}
\otimes y^{(12)}_{\sigma(1)} \cdots y^{(12)}_{\sigma(n-1)}, 
$$
$y\in F^{(2)}$ and we get 
$$
x^{(aab)}_{n-1,1,n} = [r^{(23)},y^{(13)}] \quad  
\on{and} \quad 
x^{(abb)}_{n,1,n-1} = - [r^{(12)},y^{(13)}].  
$$
 
Let us show that $y \in F_{n-1} = 
(FL_{n-1}\otimes FL_{n-1})_{\SG_{n-1}}$. 
Since $x^{(aab)}_{n-1,1,n}$ and $x^{(abb)}_{n,1,n-1}$
belong to $F^{(aab)}_{n-1,n,1}$  and $F^{(abb)}_{1,n-1,n}$, 
the commutators $[y^{(23)}_1,
\sum_{\sigma\in\SG_{n-1}} 
C_\sigma y^{(13)}_{\sigma(1)}\cdots y^{(13)}_{\sigma(n-1)}]$
and $[x^{(12)}_1,
\sum_{\sigma\in\SG_{n-1}} 
C_{\sigma^{-1}} x^{(13)}_{\sigma(1)}\cdots x^{(13)}_{\sigma(n-1)}]$
are Lie polynomials. 
The first statement implies that  
$\sum_{\sigma\in\SG_{n-1}} 
C_\sigma y^{(13)}_{\sigma(1)}\cdots y^{(13)}_{\sigma(n-1)}$
is a Lie polynomial, and the second statement implies that 
$\sum_{\sigma\in\SG_{n-1}} 
C_{\sigma^{-1}} x^{(13)}_{\sigma(1)}\cdots x^{(13)}_{\sigma(n-1)}$
is also a Lie polynomial; by virtue of Lemma \ref{observation}, 
this implies that $y$ belongs to $F_{n-1}$.

To prove the second part of the Proposition, let us 
now project equation (\ref{delta:0}) on $F^{(aabb)}_{p+1,q,1,n}$. 
Since $q\neq 1$, the contribution of the second term of (\ref{delta:0})
is zero, so that $[r^{(13)},(x^{(aab)})^{(124)}_{p,q,n} 
+ (x^{(aab)})^{(214)}_{p,q,n}] = 0$, therefore 
$$
x^{(aab)}_{p,q,n} 
+ (x^{(aab)})^{(213)}_{p,q,n} = 0 . 
$$
\hfill \qed\medskip 

Let us set now $x' = x - \delta^{(F)}_3(x)$. 

\begin{prop} \label{isr}
We have $x' = 0$. 
\end{prop} 

{\em Proof.} Let us summarize the properties of $x'$. We have
$\delta^{(4)}_F(x') = 0$, $(x')^{(aab)}_{n-1,1,n} =  
(x')^{(aab)}_{1,n-1,n} = 0$  , $(x')^{(abb)}_{n,1,n-1} =  
(x')^{(abb)}_{n,n-1,1} = 0$ and for any pair of integers 
$p,q$ such that $p+q = n$, $(x')^{(aab)}_{p,q,n} + 
((x')^{(aab)})^{(213)}_{p,q,n} = 0$ and    
$(x')^{(abb)}_{n,p,q} + ((x')^{(abb)})^{(132)}_{n,p,q} = 0$. 

The first property follows from $\delta^{(F)}_4 \circ \delta^{(F)}_3 = 0$, 
the second property follows from Proposition \ref{thl}, 1), and the third
property is a consequence of the first property and Proposition \ref{9:ab}. 

Let $k$ be an integer such that $1<k<n$, and let us project the equality 
$\delta^{(F)}_4(x') = 0$ to $F^{(aaab)}_{k,n-k,1,n+1}$ parallel to all other 
$F^{(axyb)}_{pqrs}$. We have 
$$
[r^{(12)},(F^{(abb)})^{(234)}]\subset F^{(aabb)} \oplus F^{(abbb)}, 
$$
and
$$ 
[r^{(12)},(F_{pqn}^{(aab)})^{(234)}]\subset \bigoplus_{k = 1}^p 
F^{(aaab)}_{k,p+1-k,q,n+1} \oplus F^{(abab)}_{p+1,1,q,n}, 
$$
so that for any $z\in F^{Lie,(4)}_n$, 
\begin{align*}
& (\delta^{(F)}_4(z))^{(aaab)}_{k,n-k,1,n+1} = 
\big( [r^{(12)},(z^{(aab)}_{n-1,1,n})^{(234)}] 
+ [r^{(13)},(z^{(aab)}_{n-k,k,n})^{(234)}] 
- [r^{(23)},(z^{(aab)}_{k,n-k,n})^{(134)}] 
\\ & + [r^{(34)},(z^{(aab)}_{k,n-k,n})^{(124)}]
+ [r^{(34)},(z^{(aab)}_{k,n-k,n})^{(123)}]
\big)^{(aaab)}_{k,n-k,1,n+1} . 
\end{align*}

\begin{lemma}
For any $w\in F^{(aab)}_{k,n-k,n}$, we have 
$$
\big( [r^{(13)},w^{(324)}] + [r^{(23)},w^{(134)}]
\big)^{(aaab)}_{k,n-k,1,n+1} = [r^{(34)},w^{(124)}] . 
$$
\end{lemma}

{\em Proof.} There exists a unique family $(A_{\sigma,\tau,\tau'})$
in $\KK^{\SG_{k,n-k}\times\SG_k \times\SG_{n-k}}$, such that 
$$
w = \sum_{(\sigma,\tau,\tau')\in \SG_{k,n-k}\times\SG_k \times\SG_{n-k}}
A_{\sigma,\tau,\tau'} x^{(13)}_{\tau(1)} \cdots x^{(13)}_{\tau(k)} 
\otimes x^{(23)}_{\tau'(1)} \cdots x^{(23)}_{\tau'(n-k)}
\otimes y_{\sigma(1)} \cdots y_{\sigma(n)},  
$$
where $(y_1,\ldots,y_n) = 
(y^{(13)}_1,\ldots,y^{(13)}_k,y^{(23)}_1,\ldots,y^{(23)}_{n-k})$. Then 
for each $(\sigma,\tau')\in \SG_{k,n-k}\times\SG_{n-k}$, 
$\sum_{\tau\in\SG_k} A_{\sigma,\tau,\tau'} x^{(13)}_{\tau(1)}\cdots x^{(13)}_{\tau(k)}$
is a Lie polynomial, and 
for each $(\sigma,\tau)\in \SG_{k,n-k}\times\SG_{k}$, 
$\sum_{\tau'\in\SG_k} A_{\sigma,\tau,\tau'} x^{(23)}_{\tau'(1)}\cdots x^{(23)}_{\tau'(n-k)}$
is also a Lie polynomial.  
Then 
\begin{align*}
& ([r^{(23)},w^{(134)}])^{(aaab)}_{k,n-k,1,n+1} 
\\ & 
= \sum_{(\sigma,\tau,\tau')\in
\SG_{k,n-k}\times\SG_k\times\SG_{n-k}}
A_{\sigma,\tau,\tau'} \mu^{(4)} \big( 
x^{(14)}_{\tau(1)} \cdots x^{(14)}_{\tau(k)}
\otimes x^{(23)}_1 \otimes [y^{(23)}_1, 
x^{(34)}_{\tau'(1)}  \cdots x^{(34)}_{\tau'(n-k)}]
\\ & 
\otimes y'_{\sigma(1)}  \cdots y'_{\sigma(n)} 
\big)^{(aaab)}_{k,n-k,1,n+1}
\\ & 
= 
\sum_{(\sigma,\tau,\tau')\in
\SG_{k,n-k}\times\SG_k\times\SG_{n-k}}
A_{\sigma,\tau,\tau'}  
x^{(14)}_{\tau(1)} \cdots x^{(14)}_{\tau(k)}
\otimes 
[x^{(24)}_{\tau'(1)},  \cdots [x^{(24)}_{\tau'(n-k-1)},
x^{(24)}_1]] \otimes y^{(23)}_1 
\\ & 
\otimes y''_{\sigma(1)}  \cdots [y''_{k+\tau'(n-k)},y^{(24)}_1]
\cdots y''_{\sigma(n)} ,  
\end{align*}
where $(y'_1,\ldots,y'_n) = 
(y^{(14)}_1,\ldots,y^{(14)}_k,y^{(34)}_1,\ldots,y^{(34)}_{n-k})$
and 
$$
(y''_1,\ldots,y''_n) = 
(y^{(13)}_1,\ldots,y^{(13)}_k,y^{(23)}_1,\ldots,y^{(23)}_{n-k}).
$$

Now Lemma \ref{prereut} and the fact that 
$\sum_{\tau'\in\SG_k} A_{\sigma,\tau,\tau'} 
x^{(23)}_{\tau'(1)}\cdots x^{(23)}_{\tau'(n-k)}$ is a Lie 
polynomial implies that 
\begin{align} \label{edgar}
& ([r^{(23)},w^{(134)}])^{(aaab)}_{k,n-k,1,n+1}
\\ & \nonumber 
=  \sum_{(\sigma,\tau,\tau')\in\SG_{k,n-k}\times\SG_{k}\times\SG_{n-k}}
A_{\sigma,\tau,\tau'}
\sum_{i|\sigma(i)\in\{k+1,\ldots,n\}}
x^{(14)}_{\tau(1)} \cdots x^{(14)}_{\tau(k)}
\otimes 
x^{(24)}_{\tau'(1)} \cdots x^{(24)}_{\tau'(n-k)}
\otimes y^{(23)}_1 
\\ & \nonumber 
\otimes y''_{\sigma(1)}  \cdots [y^{(24)}_1, y''_{\sigma(i)}]
\cdots y''_{\sigma(n)} . 
\end{align}

Applying $x\mapsto x^{(2134)}$ to the equality (\ref{edgar}), 
where $k$ and $w$ are replaced by $n-k$ and $w^{(213)}$
yields
\begin{align} \label{jacobs}
& ([r^{(23)},w^{(134)}])^{(aaab)}_{k,n-k,1,n+1}
\\ & \nonumber 
= 
\sum_{(\sigma,\tau,\tau')\in\SG_{k,n-k}\times\SG_{k}\times\SG_{n-k}}
A_{\sigma,\tau,\tau'}
\sum_{i|\sigma(i)\in\{1,\ldots,k\}}
x^{(14)}_{\tau(1)} \cdots x^{(14)}_{\tau(k)}
\otimes 
x^{(24)}_{\tau'(1)} \cdots x^{(24)}_{\tau'(n-k)}
\otimes y^{(23)}_1 
\\ & \nonumber 
\otimes y''_{\sigma(1)}  \cdots [y^{(24)}_1, y''_{\sigma(i)}]
\cdots y''_{\sigma(n)} . 
\end{align}
The result now follows from the addition of (\ref{edgar}) and 
(\ref{jacobs}). 
\hfill \qed\medskip 

\begin{lemma}
For any $w\in F^{(aab)}_{k,n-k,n}$, we have 
$$
\big( [r^{(34)},w^{(123)} + w^{(124)}]\big)^{(aaab)}_{k,n-k,1,n+1} = 0. 
$$
\end{lemma}

{\em Proof.} There exists a unique family $(A_\sigma)_{\sigma\in\SG_n}\in \KK^{\SG_n}$, 
such that 
$$
w = \sum_{\sigma\in\SG_n} A_\sigma x^{(13)}_1 \cdots x^{(13)}_{k}  \otimes
x^{(13)}_{k+1} \cdots x^{(13)}_n \otimes y_{\sigma(1)} \cdots y_{\sigma(n)} , 
$$
where $(y_1,\ldots,y_n) = (y^{(13)}_1,\ldots,y^{(13)}_k,y^{(23)}_1,\ldots,y^{(23)}_n)$. 
Then $\sum_{\sigma\in\SG_n} A_\sigma y_{\sigma(1)}\cdots y_{\sigma(n)}$ is a Lie
polynomial. 
We have 
\begin{align*}
& [r^{(34)},w^{(123)}] = \mu^{(4)} \big( \sum_{\sigma\in\SG_n}
A_\sigma x^{(13)}_1\cdots x^{(13)}_k  \otimes x^{(23)}_{k+1}\cdots x^{(23)}_n
\otimes [x^{(34)}_1,y_{\sigma(1)}\cdots y_{\sigma(n)}] \otimes y^{(34)}_1 
\big) 
\\ & 
= \mu^{(4)} \big( \sum_{\sigma\in\SG_n}
A_\sigma x^{(13)}_1\cdots x^{(13)}_k  \otimes x^{(23)}_{k+1}\cdots x^{(23)}_n
\otimes x^{(34)}_1 \otimes [y'_{\sigma(1)},\cdots, [y'_{\sigma(n)}, y^{(34)}_1]] 
\big) 
\\ & = - [r^{(34)},w^{(124)}]
\end{align*}
where $(y'_1,\ldots,y'_n) = (y^{(14)}_1,\ldots,y^{(14)}_k,y^{(24)}_1,\ldots,y^{(24)}_n)$. 
The last equality follows from Lemma \ref{chrono} and the fact that 
$\sum_{\sigma\in\SG_n} A_\sigma y_{\sigma(1)}\cdots y_{\sigma(n)}$ is a Lie
polynomial. 
\hfill \qed\medskip 

{\em End of proof of Proposition \ref{isr}.}
By virtue of the two Lemmas above, and 
since $(x')^{(aab)}_{n-1,1,n} = 0$, and 
$(x')^{(aab)}_{k,n-k,n} + ((x')^{(aab)}_{n-k,k,n})^{(213)} = 0$, 
$(\delta^{(F)}_4(x'))^{(aaab)}_{k,n-k,1,n+1} = 0$ yields 
$$
[r^{(34)},(x')^{(aab)}_{k,n-k,n} ] = 0, 
$$
which implies that $(x')^{(aab)}_{k,n-k,n} = 0$, so $(x')^{(aab)} = 0$.
One shows in the same way that $(x')^{(abb)} = 0$, so $x' = 0$.
\hfill \qed\medskip 

\begin{remark}
The projection of 
$\delta^{(F)}_4(x') = 0$ to $F^{(aaab)}_{k,1,n-k,n+1}$ also 
yields the result. On the other hand, the projection of
$\delta^{(F)}_4(x)$ to $F^{(aaab)}_{1,k,n-k,n+1}$ is identically
zero. \hfill \qed\medskip 
\end{remark}

{\em End of the computation of $H^3_n$.} We have shown that for any 
$n>2$ and $x\in F^{Lie,(3)}_n$ such that $\delta^{(F)}_4(x) = 0$, 
there exists $y\in F_{n-1}$ such that $x = \delta^{(F)}_3(y)$. 
Therefore $H^3_n = 0$. \hfill \qed\medskip

\newpage

\section{Universal shuffle algebras (proof of Theorem \ref{thm:univ})}
\label{sect:proof:thm:univ}

In this Section, we define universal shuffle algebras $\Sh^{(F)}_k$. 
These algebras have universal properties with respect to the tensor 
powers $\Sh(\A)^{\otimes k}$, where $\A$ is a Lie algebra endowed with 
a solution $r_\A\in\A\otimes\A$ of CYBE. 

\subsection{Definition of $\Sh^{(F)}_k$}

Let $k$ be an integer $\geq 0$. We set $\Sh^{(F)}_k = \KK$ for 
$k = 0$ and $k = 1$. When $k\geq 2$, we put the following definitions. 
If $\ualpha = (\al_1,\ldots,\al_k)\in \NN^k$, we define
$X(\ual)$ as the set of all maps $x : I_{\ual} \to \{a,b\}$
(recall that $I_\ual = \{(i,\beta) | i \in\{1,\ldots,n\}, \beta\in 
\{1,\ldots,\al_i\}\}$). If $x\in X(\ual)$, let us denote by 
$P(x,\ual)$ the set of all maps $\up : x^{-1}(a)\times 
x^{-1}(b) \to\NN$, such that 
for any $(i,\beta)\in x^{-1}(a)$ and $(i',\beta')\in x^{-1}(b)$, 
$\up((i,\beta),(i',\beta')) = 0$ whenever $i = i'$ or 
$(i,\beta)>(i',\beta')$ in the lexicographical order.  

We then set 
\begin{equation} \label{def:univ:sh}
\Sh^{(F)}_k 
= \bigoplus_{\ual\in\NN^k} \big(
\bigoplus_{x\in X(\ual)} \big( \bigoplus_{\up\in P(x,\ual)}
F^{(x(1,1)\ldots x(1,\al_1)x(2,1)\ldots x(k,\al_k))}_{
\big( \up((i,\beta),(i',\beta')) 
\big)_{((i,\beta),(i',\beta'))\in x^{-1}(a) \times x^{-1}(b)}}
\big) \big) , 
\end{equation}
where the space $F^{(x_1\ldots x_n)}_{(p(k,l))_{(k,l)\in K\times L}}$ 
is defined by 
(\ref{def:lie:space}). The summand corresponding by $\ual = 0$
is $\KK$.  
(For each $i$, one should think 
of the tensor product of all the tensor factors  
indexed by $(i,\al)$ as of the analogue of 
the $i$th factor of $\Sh(\A)^{\otimes k}$.) 

We denote by $\Sh^{(F)}_k(\ual)$ 
(resp., $\Sh^{(F)}_k(\ual,x)$, 
$\Sh^{(F)}_k(\ual,x,\up)$) the graded component of 
$\Sh^{(F)}_k$ corresponding to $\ual$ (resp., to $(\ual,x)$,
$(\ual,x,\up)$). 

\subsection{Operations of $\Sh^{(F)}_k$}

\subsubsection{The multiplication $m_{\Sh^{(F)}_k}$}

Let us first define a bilinear map 
$$
\ell : ( \otimes_{i = 1}^l FL_{\al_i + \al'_i}) 
\times (\Sh^{(F)}_k(\ual) \otimes \Sh^{(F)}_k(\ual') )
\to \Sh^{(F)}_k(\underline{1})  , 
$$
where $\underline{1}$ is the element of $\NN^k$ with all 
components equal to $1$. 

Assume that $u$ and $v$ are decomposed as $u = \otimes_{i = 1}^k 
(\otimes_{\beta = 1}^{\al_i} u_{i,\beta})$ and
$v = \otimes_{i = 1}^k 
(\otimes_{\beta = 1}^{\al'_i} v_{i,\beta})$, then 
we view $\otimes_{i = 1}^k L_i(u_{i,1},\ldots,u_{i,\beta_i},
v_{i,1},\ldots,v_{i,\beta'_i})$ as an element of $G^{(k)}$. 
Then it follows from Corollary \ref{prop:Lie:n} that 
$\mu^{(k)} (\otimes_{i = 1}^k L_i(u_{i,1},\ldots,u_{i,\beta_i},
v_{i,1},\ldots,v_{i,\beta'_i}))$ belongs to 
$\oplus_{x\in \on{Map}(\{1,\ldots,k\},\{a,b\}) }
F^{(x(1)\ldots x(k))}_{\up}$.  
The latter space is exactly 
$\Sh^{(F)}_k(\underline{1})$, and we set 
$$
\ell(\otimes_{i = 1}^k L_i,u\otimes v) = \mu^{(k)}
(\otimes_{i = 1}^k L_i(u_{i,1},\ldots,u_{i,\beta_i},
v_{i,1},\ldots,v_{i,\beta'_i})) . 
$$

Let us fix $\ugamma\in \NN_{>0}^k$ and for each $i$, let us 
fix $\gamma_i$-partitions of $\al_i$ and $\al'_i$, namely 
$\al_i = \al_{i1} + \cdots + \al_{i\gamma_i}$ 
and $\al'_i = \al'_{i1} + \cdots + \al'_{i\gamma_i}$. 
Then we define a bilinear map 
$$
\ell_{\ual,\ual',\ugamma,(\al_{ij}),(\al'_{ij})} : 
\big(  \otimes_{i = 1}^l (\otimes_{\beta = 1}^{\gamma_i}
FL_{\al_{i\beta} + \al'_{i\beta}}) \big) 
\times (\Sh^{(F)}_k(\ual) \otimes \Sh^{(F)}_k(\ual') )
\to \Sh^{(F)}_k(\ugamma)   
$$
as follows. Let us define $\Delta_{\ugamma,(\al_{ij})}$
as the linear map 
$$
\Delta_{\ugamma,(\al_{ij})} : 
\Sh^{(F)}_k(\ual) \to 
\Sh^{(F)}_{\sum_i\gamma_i}(\al_11,\ldots,\al_{1,\gamma_1},
\al_{21},\ldots,\al_{k\gamma_k})
$$
defined as the canonical injection of 
$\Sh^{(F)}_k(\ual,x,\up)$ into $\Sh^{(F)}_{\sum_i\gamma_i}
( (\al_{ij})_{ (i,j)| 1\leq i\leq n, 1\leq j\leq\gamma_i} , 
x',\up')$, where $x'$ is the composition of $x$ with the 
lexicographical bijection 
$I_{\al_{11},\ldots,\al_{k\gamma_k}} \to I_\ual$, 
and $\up'$  is the composition of $\up$ with the square of this
bijection. 

In the same way, $\Delta_{1,\underline{1}}$ is an injection 
of $\Sh^{(F)}_k(\ugamma)$ in $\Sh^{(F)}_{\sum_i\gamma_i}(\underline{1})$, 
associated with the partitions $(1,\ldots,1)$ ($\gamma_i$ times $1$)
of each $\gamma_i$.

Then if $u\in \Sh^{(F)}_k(\ual)$  and $v\in \Sh^{(F)}_k(\ual')$, 
and if for each $(i,\beta)\in I_\ugamma$, $L_{i,\beta}\in 
FL_{\al_{i\beta} + \al'_{i\beta}}$, then 
$$
\ell \big(
\otimes_{i = 1}^k (\otimes_{\beta = 1}^{\gamma_i} L_{i\beta}) , 
\Delta_{\ugamma,(\al_{ij})}(u) \otimes 
\Delta_{\ugamma,(\al'_{ij})}(v) \big)
$$ 
is in the image of $\Delta_{1,\underline{1}}$, and we denote by 
$$
\ell_{\ual,\ual',\ugamma,(\al_{ij}),(\al'_{ij})} 
\big(
\otimes_{i = 1}^k (\otimes_{\beta = 1}^{\gamma_i} L_{i\beta}) , 
u \otimes v \big)
$$
as the preimage of this element by 
$\Delta_{1,\underline{1}}$. 

Then there exists a unique linear map $m_{\Sh^{(F)}_k} 
: \Sh^{(F)}_k \otimes \Sh^{(F)}_k  
\to \Sh^{(F)}_k$, such that  
if $\ual\in\NN^k$ and $\ual'\in\NN^{k'}$, and if 
$u\in\Sh^{(F)}_k(\ual)$ and $u'\in\Sh^{(F)}_k(\ual')$, 
\begin{align*}
& m_{\Sh^{(F)}_k} (u\otimes v) = 
\sum_{\ugamma\in\NN^k} \sum_{(\al_{ij})\in P(\ual,\ugamma),
(\al'_{ij})\in P(\ual',\ugamma)}
\\ & 
\ell_{\ual,\ual',\ugamma,(\al_{ij}),(\al'_{ij})} 
\big(
\otimes_{i = 1}^k (\otimes_{\beta = 1}^{\gamma_i} B_{\al_{i\beta},
\al'_{i\beta}}) ,  u \otimes v \big) , 
\end{align*}
where $P(\ual,\ugamma)$ is the set of collections 
$\big( (\al_{1j})_{j = 1,\ldots,\gamma_1}, \ldots, 
(\al_{kj})_{j = 1,\ldots,\gamma_k}\big)$, where for each $i$,
$(\al_{ij})_{j = 1,\ldots,\gamma_i}$ is a $\gamma_i$-partition 
of $\al_i$. 

\begin{prop}
 $m_{\Sh^{(F)}_k}$ is associative. 
\end{prop}

{\em Proof.} This follows from the fact that the $B_{pq}$
satisfy the identities (\ref{cond:Bpq}),
 and from Proposition \ref{prop:normal:ordering}
with $\al = 3$. 
\hfill \qed\medskip 

\subsubsection{The maps $x\mapsto x^{(i_1\ldots i_k)}$}

Let $k$ and $l$ be integers such that $k\leq l$, and 
let $(i_1,\ldots,i_k)$ be integers in $\{1,\ldots,l\}$ such that 
$i_1<i_2\ldots <i_k$. If $\ual\in\NN^k$, define 
$\ual^{(i_1\ldots i_k)}$ by $(\ual^{(i_1\ldots i_k)})_{i_s} = \al_s$
and $(\ual^{(i_1\ldots i_k)})_{t} = 0$ if $t\notin\{i_1,\ldots,i_k\}$. 
If $x\in X(\ual)$, define $x^{(i_1\ldots i_k)}$
as the element of $X(\ual^{(i_1\ldots i_k)})$ equal to  
the composition of $x$ with the lexicographical bijection 
between $I_{\ual^{(i_1 \ldots i_k)}}$ and $I_\ual$.  
If $\up\in P(x,\ual)$, define $\up^{(i_1\ldots i_l)}$
as the element of $P(x^{(i_1\ldots i_k)},\ual^{(i_1\ldots i_k)})$ 
given by the compisition of $\up$ with the square of this
bijection. 

Let $x\mapsto x^{(i_1\ldots i_k)}$ 
be the linear map from $\Sh^{(F)}_k$  to $\Sh^{(F)}_l$
defined as the direct sum of all canonical injections
of $\Sh^{(F)}_k(\ual,x,\up)$  into $\Sh^{(F)}_l(\ual^{(i_1\ldots i_k)}
,x^{(i_1\ldots i_k)},\up^{(i_1\ldots i_k)})$.  

Then $x\mapsto x^{(i_1\ldots i_k)}$ is an algebra morphism, 
and if $j_1,\ldots,j_{k'}$ are such that $1\leq j_1 < \ldots < j_{k'}
\leq l$ and that $\{i_1,\ldots,i_{k}\}$  and $\{j_1,\ldots,j_{k'}\}$
are disjoint, then $[x^{(i_1\ldots i_k)},y^{(j_1\ldots j_{k'})}]
= 0$ for any $x\in \Sh^{(F)}_k$  and $y\in \Sh^{(F)}_l$.

\subsubsection{The morphisms $\Delta_{k,i}$}

Assume that $k$ is an integer and let $i$ be an integer
such that $1\leq i \leq k$. 

If $\ual\in\NN^k$ and $\beta = 1,\ldots,\al_i$, define 
$\ual(i,\beta)$ as the element of $\NN^{k+1}$ equal to
$(\al_1,\ldots,\al_{i-1},\beta,\al_i - \beta,\al_{i+1},\ldots,\al_k)$. 
If $x\in X(\ual)$, define $x(i,\beta) = x'$ as the element of 
$X(\ual(i,\beta))$ such that $x'(j,\gamma) = x(j,\gamma)$ if
$j\leq i$, $x'(j,\gamma) = x(j-1,\gamma)$ if
$j\geq i+2$, and $x'(i+1,\gamma) = x(i,\beta + \gamma)$; 
this is the composition of $x$ with the lexicographical 
bijection between $I_{\ual}$ and $I_{\ual(i,\beta)}$.
If $\up\in P(x,\ual)$, define $\up(i,\beta)$ as the 
composition of $\up$ with the square of this bijection.

Define $\Delta_{k,\ual,i,\beta}$ as the canonical injection map from 
$\Sh^{(F)}_k(\ual,x,\up)$ to 
$$
\Sh^{(F)}_{k+1}(\ual(i,\beta),x(i,\beta),
\up(i,\beta))
$$ 
and $\Delta_{k,i}$ as 
$\oplus_{\ual\in\NN^k, \beta = 1,\ldots,\al_i} \Delta_{k,\ual,i,\beta}$. 
Then $\Delta_{k,i}$ is a linear 
map from $\Sh^{(F)}_k$  to $\Sh^{(F)}_{k+1}$. 

Moreover,  $\Delta_{k,i}$ is an algebra morphism, and  
the maps $\Delta_{k,i}$ and $x\mapsto x^{(i_1\ldots i_k)}$
satisfy coassociativity and compatibility rules
$$
\Delta_{k+1,j}\circ\Delta_{k,i} = \Delta_{k+1,i+1}\circ \Delta_{k,j}
$$
if $j\leq i$,  
$$
\Delta_{l,i}(x^{(i_1\ldots i_k)}) =x^{(i_1,\ldots, i_s, i_{s+1} + 1,
\ldots i_{k}+1)}
$$ 
if $s\neq 0,k$ and $i_s < i < i_{s+1}$, or if 
$s = 0$ and $i<i_1$, or if $s = k$ and $i>i_k$, 
and 
$$
\Delta_{l,i_s}(x^{(i_1\ldots i_k)}) =
(\Delta_{k,s}(x))^{(i_1,\ldots, i_s, i_{s} + 1,
\ldots i_{k}+1)}. 
$$

\subsection{Universal properties of $\Sh^{(F)}_k$}

Let $(\A,r_\A)$ be the pair of a Lie algebra and a solution
$r_\A\in\A\otimes\A$ of CYBE. If $k$ is an integer 
$\geq 0$, $\ual\in \NN^k$, $x\in X(\ual)$,  
let us denote by $\kappa_{\A,r_\A}(k,\ual,x)$ the linear
map from $\Sh^{(F)}_k$ to $\Sh(\A)^{\otimes k}$
given by the composition of 
$\kappa^{(x(1,1),\ldots,x(k,\al_k))}_{\A,r_\A}$, the 
canonical isomorphism $\A^{\otimes (\al_{1} + \cdots + \al_{k})}
\to \otimes_{i = 1}^{k} \A^{\otimes \al_i}$ and the tensor 
product $\otimes_{i = 1}^k \iota_{\al_i}$, where $\iota_s$
is the canonical injection of $\A^{\otimes s}$ in $\Sh(\A)$ as
its part of degree $s$. 

Define $\kappa_{\A,r_\A}^{\Sh^{(F)}_k}$ as the linear map 
from $\Sh^{(F)}_k$ to $\Sh(\A)^{\otimes k}$ equal to the 
sum 
$$
\sum_{\ual\in\NN^k}\sum_{x\in X(\ual)}
\kappa_{\A,r_\A}(k,\ual,x).
$$ 

\begin{prop} 
$\kappa_{\A,r_\A}^{\Sh^{(F)}_k}$ is an algebra morphism. 
Moreover, we have 
$$
\kappa_{\A,r_\A}^{\Sh^{(F)}_{k+1}} \circ 
\Delta_{k,i} = 
\id_{\Sh(\A)}^{\otimes (i-1)}\otimes\Delta_{\Sh(\A)}
\otimes \id_{\Sh(\A)}^{\otimes (k-i-1)}\circ 
\kappa_{\A,r_\A}^{\Sh^{(F)}_{k}} 
$$
and
$$
\kappa_{\A,r_\A}^{\Sh^{(F)}_{l}}(x^{(i_1\ldots i_k)})  
= (\kappa_{\A,r_\A}^{\Sh^{(F)}_{k}}(x))^{(i_1\ldots i_k)}   
$$
for any $x\in \Sh^{(F)}_k$. 
\end{prop}

The proof is straightforward. This Proposition explains 
why $\Sh^{(F)}_k,m_k,\Delta_{k,i}$ and $x\mapsto x^{(i_1\ldots i_k)}$
should be viewed as universal versions of $\Sh(\A)^{\otimes k},
m_{\Sh(\A)}^{\otimes k}, \id_{\Sh(\A)}^{\otimes (i-1)}\otimes\Delta_{\Sh(\A)}
\otimes \id_{\Sh(\A)}^{\otimes (k-i-1)}$ and the map
$x\mapsto x^{(i_1\ldots i_k)}$. 

\subsection{Proof of Theorem \ref{thm:univ}}

In the statement of Proposition \ref{prop:leeds:first}, 
we may replace $A^{\otimes k}$
($k = 2,3,4$) by $\Sh^{(F)}_k$ and $r_A$ by the element 
$r_{\Sh_2^{(F)}}\in \Sh^{(F)}_2((1,1),(1\mapsto a,2\mapsto b),
(1,2)\mapsto 1)
\subset \Sh^{(F)}_2$ equal to 
$$
r\in F_1 = F^{(ab)}_1 . 
$$
The proof is a direct 
transposition of the proof of Proposition \ref{prop:leeds:first}, 
which is in Appendix \ref{proof:problem:3}.  

Then using the maps $\Delta_{k,i}$ and $x\mapsto x^{(i_1\ldots i_k)}$, 
we may reproduce step by step the proof of Theorem \ref{training}. 
This proves Theorem \ref{thm:univ}. 
\hfill \qed\medskip   

\end{appendix}

\end{document}